\documentclass{amsart}
\pdfoutput=1

\usepackage[utf8]{inputenc}
\usepackage{amssymb}
\usepackage{amsmath}
\usepackage{amsthm}
\usepackage{bm}
\usepackage{graphicx}
\usepackage[square,numbers,sort]{natbib}
\usepackage{xspace}
\usepackage{overpic}
\usepackage{algorithm2e}
\usepackage{enumerate}  % For lower-case Roman numbers
\usepackage{enumitem}  % For lower-case Roman numbers
\usepackage{subfigure}
\usepackage{todonotes}
\definecolor{lightred}{rgb}{1,.4,.4}

\usepackage{colonequals}
\usepackage{hyperref}

\newtheorem{theorem}{Theorem}[section]

\newtheorem{assumption}[theorem]{Assumption}
\newtheorem{lemma}[theorem]{Lemma}
\newtheorem{proposition}[theorem]{Proposition}
\theoremstyle{remark}
\newtheorem{remark}[theorem]{Remark}

\newcommand{\Lcal}{\mathcal{L}}
\newcommand{\Jcal}{\mathcal{J}}
\newcommand{\Ccal}{\mathcal{C}}
\newcommand{\Vcal}{\mathcal{V}}

\newcommand{\Mcal}{\mathcal{M}}
\newcommand{\be}{\mathbf{e}}
\newcommand{\Eig}{\bm{\lambda}}
\newcommand{\Eigfunc}{\operatorname{Eig}}

\newcommand{\bvarepsilon}{{\bm{\varepsilon}}}
\newcommand{\bve}{{\bm{\varepsilon}}}

\DeclareRobustCommand{\Bu}{{\boldsymbol{\mathnormal u}}}
\DeclareRobustCommand{\Bv}{{\boldsymbol{\mathnormal v}}}
\DeclareRobustCommand{\BU}{{\boldsymbol{\mathnormal U}}}
\DeclareRobustCommand{\BV}{{\boldsymbol{\mathnormal V}}}
\DeclareRobustCommand{\BW}{{\boldsymbol{\mathnormal W}}}
\newcommand{\bu}{\Bu}
\newcommand{\bv}{\Bv}

\newcommand{\tr}{\operatorname{tr}}
\newcommand{\R}{\mathbb{R}}   % Real numbers
\renewcommand{\S}{\mathbb{S}}   % Symmetric matrices
\newcommand{\Dim}{m}
\newcommand{\TNNMGEX}{TNNMG-EX\xspace}
\newcommand{\TNNMGPRE}{TNNMG-PRE\xspace}

\providecommand{\inner}[2]{\left\langle #1, #2 \right\rangle}

\providecommand{\abs}[1]{\lvert#1\rvert}

\DeclareMathOperator*{\argmin}{arg\,min}

\graphicspath{{gfx/}}

\makeatletter
\AtBeginDocument{
  \hypersetup{
    pdftitle = {\@title},
    pdfauthor = {\@author}
  }
}
\makeatother

\title[Nonsmooth multigrid for phase-field brittle-fracture problems]{Truncated Nonsmooth Newton Multigrid for phase-field brittle-fracture problems, with analysis}

\thanks{%
  The authors gratefully acknowledge the financial support by the
  German Federal Ministry of Education and Research through the
  ParaPhase project within the framework
  ``IKT 2020 -- Forschung für Innovationen''
  (project numbers 01-H15005C, 01-15005D, 01-15005E)%
}

\author{Carsten Gräser}
\address{Carsten Gräser\\
Friedrich-Alexander-Universität Erlangen--Nürnberg\\
Department Mathematik\\
Cauerstr. 11\\
91058 Erlangen\\
Germany}
\email{graeser@math.fau.de}

\author{Daniel Kienle}
\address{Daniel Kienle\\
Universität Stuttgart\\
Institut für Angewandte Mechanik\\
Pfaffenwaldring 7\\
70569 Stuttgart\\
Germany}
\email{kienle@mechbau.uni-stuttgart.de}

\author{Oliver Sander}
\address{Oliver Sander\\
Technische Universität Dresden\\
Institut für Numerische Mathematik\\
Zellescher Weg 12--14\\
01069 Dresden\\
Germany\\
ORCID: \href{https://orcid.org/0000-0003-1093-6374}{0000-0003-1093-6374}
}
\email{oliver.sander@tu-dresden.de}

\begin{document}

\begin{abstract}
  We propose the Truncated Nonsmooth Newton Multigrid Method (TNNMG)
  as a solver for the spatial problems of the small-strain brittle-fracture phase-field equations.
  TNNMG is a nonsmooth multigrid method that can
  solve biconvex, block-separably nonsmooth minimization problems with linear time complexity.
  It exploits the variational structure inherent in the problem,
  and handles the pointwise irreversibility
  constraint on the damage variable directly, without regularization or the
  introduction of a local history field.
  In the paper we introduce the method and show how it can be applied to several established
  models of phase-field brittle fracture.  We then prove convergence of the solver
  to a solution of the nonsmooth Euler--Lagrange equations of the spatial problem for any load and initial iterate.
  On the way, we show several crucial convexity and regularity properties of the models considered here.
  Numerical comparisons to an operator-splitting algorithm show a
  considerable speed increase, without loss of robustness.
\end{abstract}

\maketitle

%\bigskip


\section{Introduction}
\label{sec:introduction}

The equations of phase-field models of brittle fracture present a number of
challenges to the designers of numerical solution algorithms~\cite{ambati_gerasimov_delorenzis:2015}.
Even in the small-strain case the equations are nonlinear, due to the multiplicative
coupling of the mechanical stresses to the degradation function.
At the same time the non-healing condition introduces an inequality constraint.
Finally, eigenvalue-based splittings of the energy density as in~\cite{miehe+welschinger+hofacker10a}
make the equations nondifferentiable.

In this paper we focus on the spatial problems of small-strain brittle-fracture phase-field models
obtained by a suitable time discretization.
The standard approach to solving these spatial problems is based on operator
splitting.  Algorithms based on this approach, also known as staggered schemes, alternate between solving a displacement
problem with fixed damage and a damage problem with fixed displacement.
Both subproblems are elliptic and well-understood, and such methods are therefore
straightforward to implement.  The method can be interpreted as a nonlinear
Gauß--Seidel method~\cite{farrell+maurini17}, which provides a natural framework
for convergence proofs.
Applications of the operator splitting scheme and
its extensions appear, e.g., in \cite{bourdin07,burke+ortner+sueli10}.
With particular semi-implicit time discretizations it is also possible
to solve the spatial problems by solving only one damage problem and one
displacement problem~\cite{miehe+hofacker+welschinger10}.  This is very fast,
but works only if the load steps are small enough.

In contrast, other works propose monolithic solution schemes based on Newton's
method \cite{gerasimov+delorenzis16,farrell+maurini17,wick17modified,wick17error,wambacq_ulloa_lombaert_francois:2021}.
For the unmodified Newton method only local convergence can be shown,
and failure to converge for large load steps is readily observed in practice~\cite{wu_huang_nguyen:2020}.
Therefore, various authors have proposed extensions or modifications of the
Newton idea to stabilize the method.
In \cite{gerasimov+delorenzis16} a line search strategy is
applied to enlarge the domain of convergence of Newton's method.
\citet{wu_huang_nguyen:2020} and \citet{kristensen_martinez:2020}
propose to use the BFGS quasi-Newton algorithm, claiming that it is
more stable than Newton's method and more efficient than operator splitting.
\citet{heister_wheeler_wick:2015} proposed a modified Newton scheme
which was later improved by \citet{wick17modified} with an adaptive transition from Newton’s method to the modified Newton
scheme.
Recent results in~\cite{kopanicakova_kothari_krause:2022} suggest that nonlinear preconditioning
can speed up the solution process tremendously.
Finally, the authors of \cite{singh_verhoosel_deborst_vanbrummelen:2016,may_vignollet_deborst:2016}
suggest an arc-length method based on the fracture surface, and an adaptive time stepping
scheme to enhance the robustness.
% However, these strategies are not always effective and their applicability to non-standard phase-field damage
% models, i.e., the AT1 and PF-CZM, needs further investigation.''
In summary, while monolithic Newton-type methods are reported to be faster than
operator-splitting algorithms, the latter ones are more robust~\cite{ambati_gerasimov_delorenzis:2015,wu_huang_nguyen:2020}.

Various approaches are used in the literature to deal with the damage irreversibility.
A natural approach is to regularize the constraint, as investigated
in \cite{miehe+welschinger+hofacker10a,wheeler+wick+wollner14,kopanicakova_kothari_krause:2022}.
This leads to an additional parameter, and to ill-conditioned tangent matrices~\cite{gerasimov_delorenzis:2019}.
Interior-point solvers implement an automatic control of the new parameter;
they are investigated in~\cite{wambacq_ulloa_lombaert_francois:2021}.
An alternative approach considers
the thermodynamic driving force of the fracture phase-field as a
global unknown yielding a three-field formulation which results in a saddle-point
principle~\cite{miehe+welschinger+hofacker10a}.
A third formulation
considers the Karush--Kuhn--Tucker conditions and shifts the
thermodynamic driving force of the fracture phase-field into a local
history field representing the maximum over time of the elastic
energy~\cite{miehe+hofacker+welschinger10}.
This approach, frequently known as the $\mathcal{H}$-field technique,
therefore trades the inequality constraint for the nondifferentiable
maximum function.
  Alternatively, the time discrete
  problems can be reformulated as semismooth systems by means of
  so-called complementarity functions. This strategy can be combined
  with monolithic semismooth Newton techniques~\cite{MangWickWollner2020}
  or nested Newton and active set methods~\cite{heister_wheeler_wick:2015}.
Unfortunately, these approaches spoil the variational structure of the spatial problems.
Augmented Lagrangian solvers as in \cite{wheeler+wick+wollner14} introduce extra variables.
Closest in spirit to the present manuscript is the use of bound-constrained
optimization solvers, used, e.g.,  in \cite{amor_marigo_maurini:2009,farrell+maurini17,wu:2018,burke+ortner+sueli10}.
None of these approaches are fully satisfactory.

The effect of the nondifferentiable terms caused by anisotropic splittings of the mechanical
energy density as in~\cite{miehe+welschinger+hofacker10a,steinke+kaliske18,amor_marigo_maurini:2009} is rarely discussed
in the literature.  Hybrid formulations like the one proposed in \cite{ambati_gerasimov_delorenzis:2015}
try to overcome the additional computational difficulties of these splittings by further changes to the model,
again at the cost of sacrificing the variational structure.

All these approaches are slow in the sense that they have to solve global partial differential equations
at each Newton or operator-splitting iteration.
This is expensive, even if efficient multigrid methods are used for the linear
subproblems (as, e.g., in~\cite{jodlbauer_langer_wick:2020}).
When the methods use direct sparse solvers for the linear tangent problems, memory
consumption can become problematic, too.
At the same time, the problem of
small-strain phase-field brittle fracture has a lot of elegant variational structure; in particular, it fits
directly into the rate-independent framework of \citet{mielke_roubicek:2015}.
As a consequence, implicit time discretization leads to a sequence of coercive minimization problems
for the displacement and damage fields together.  These problems are not convex, but they are
biconvex, i.e., convex (even strongly convex) in each variable separately.
Pointwise inequality restrictions $\dot{d} \ge 0$ to handle the irreversibility
of the fracture process as proposed in \cite{miehe+welschinger+hofacker10a} reduce the smoothness of the
objective functional, but do not influence its convexity or coercivity properties.
The same holds for anisotropic energy splittings based on linear quantities or the eigenvalues of the mechanical strain.

Recent years have shown that nonsmooth multigrid methods
are able to solve variational nonsmooth problems from mechanics efficiently
without the need for solving global linear systems of
equations~\cite{GraeserSander2019,sander_jaap:2020,graeser_sack_sander:2009,GraeserKornhuber2009b,kornhuber:1997, KopanicakovaKrause2020}.
This can make them vastly more efficient than operator-splitting or Newton-based methods.
As there are no sparse matrix factorizations, memory consumption remains
linear in the number of unknowns.
In addition, these multigrid methods can be shown to converge globally (i.e., from any initial iterate
and for any load step) to a stationary point of the objective functional.
The proof exploits the above-mentioned variational structure,
together with certain separability properties.
As one such method, the Truncated Nonsmooth Newton Multigrid method (TNNMG) can treat the pointwise constraints of the
increment problems directly, i.e., without artificial regularization or tricks like
the $\mathcal{H}$-field technique~\cite{miehe+hofacker+welschinger10}. The idea is that TNNMG only needs to handle these constraints
in a series of low-dimensional subproblems, each of which is easy to solve
by itself. As a consequence, solving the problems with constraints is not
appreciably slower than solving the corresponding unconstrained problem.

In this paper we show how the TNNMG method can be used to solve small-strain brittle-fracture
problems. This involves in particular verifying that the increment functionals have the
required convexity and smoothness properties.  We do this for a range of different degradation functions
and local crack surface densities (including the standard Ambrosio--Tortorelli functionals
of type~1 and~2), closely related to the family of models considered in~\cite{burke_ortner_sueli:2013}.
We cover elastic energies with various types of anisotropic splittings,
including the splitting based on strain eigenvalues of \citet{miehe+welschinger+hofacker10a}.
For the proofs we use results from the convex analysis
of spectral functions~\cite{QiYang2003,baes:2007}.
Extension to the slightly more general damage models of~\cite{marigo_maurini_pham:2016}
is straightforward, provided the stored elastic energy has the required convexity
and smoothness properties.
In contrast to the multilevel trust region method proposed in \cite{KopanicakovaKrause2020},
the TNNMG method presented here relies on nonsmooth Newton techniques
leading to linear subproblems and thus gives more flexibility
in the selection of coarse grid solvers.

%\medskip

The paper is organized as follows: Section~\ref{sec:phase_field_models} discusses
a framework of small-strain phase-field models for brittle fracture,
and shows the range of applicability of the TNNMG solver.
Section~\ref{sec:time_discrete_model} introduces the natural fully implicit time discretization,
and proves existence of solutions for the spatial problems.
In both sections we pay particular attention to the mathematical properties
of the energy functionals.  In Section~\ref{sec:tnnmg_method},
finally, we introduce the TNNMG method. We explain its construction,
discuss various algorithmic options, and prove that it converges globally
to stationary points of the increment energy functional.
The numerical efficiency
is then demonstrated in Section~\ref{sec:numerical_examples}.
Our reference for comparison, briefly revisited in Section~\ref{sec:operator_splitting},
is the operator-splitting iteration
proposed in~\cite{burke_ortner_sueli:2013}, which uses a projected Newton method
\cite{bertsekas:1982} for the constrained damage problems.
We compare the solvers for
two- and three-dimensional example problems with different forms
of the local crack surface density, and with and without spectral splittings.
We observe a noticeable performance increase, without loss of robustness.

\section{Phase-field models of brittle fracture}
\label{sec:phase_field_models}

This section presents a range of phase-field models for brittle fracture, and discusses
its smoothness and convexity properties.

Consider a deformable $\Dim$-dimensional object represented by a domain $\Omega \in \R^\Dim$.
The deformation
of such an object is characterized by a displacement field $\bu : \Omega \to \R^\Dim$. The object
is supposed to exhibit small-strain deformations and elastic material behavior only,
and we therefore introduce the linearized strain tensor
$\bve(\Bu) \colonequals \frac{1}{2} (\nabla \bu + \nabla \bu^T)$.
Following~\cite{miehe+welschinger+hofacker10a}, we model the fracturing by a scalar damage
field $d : \Omega \to [0,1]$, where $d=0$ signifies intact material, and $d=1$ a fully broken one.
Dirichlet boundary conditions can be posed both for the displacement and for the damage field.
For this we select two not necessarily equal subsets $\Gamma_{D,\bu},\Gamma_{D,d} \subset \partial \Omega$ of
the domain boundary, and require
\begin{equation*}
 \bu = \bu_0 \qquad \text{on $\Gamma_{D,\bu}$},
 \qquad
 d = d_0 \qquad \text{on $\Gamma_{D,d}$},
\end{equation*}
where $\bu_0$ and $d_0$ are two given functions.

Displacement and damage field evolve together, governed by a system of coupled nonsmooth
partial differential equations.  Disregarding inertia effects, we obtain a rate-independent system
in the sense of \citet{mielke_roubicek:2015}.  Such a system can be written using the Biot equation
\begin{equation}
\label{eq:biot_equation}
 D_{(\bu,d)} \mathcal{E}(t,\bu,d) + \partial_{\dot{d}} \mathcal{R}(d,\dot{d}) \ni 0,
\end{equation}
where $D_{(\bu,d)} \mathcal{E}(t,\bu,d)$ means the G\^ateaux derivative with respect to the second and third arguments
of $\mathcal{E}$,
and $\partial_{\dot{d}} \mathcal{R}(d,\dot{d})$ is the convex subdifferential with respect to the second argument
of the dissipation potential $\mathcal{R}$.

In this equation, $\mathcal{E}$ is a potential energy, which we assume to be of the form
\begin{align}
\label{eq:potential_energy}
 \mathcal{E}(t,\Bu, d)
 & =
 \int_\Omega \psi(\bve(\bu),d)\,dV + \int_\Omega g_c\gamma(d,\nabla d)\,dV
  + P_\text{ext}(t,\bu) + \int_\Omega I_{[0,1]} (d)\,dV.
\end{align}
The term $\psi$ is a degraded elastic energy density, and will be discussed in detail
in Section~\ref{sec:elastic_energy}. The term $\gamma$ models the local crack surface density,
and will be discussed in Section~\ref{sec:crack_surface_density}.
The number $g_c$ is Griffith's critical energy release rate, a material parameter.
$P_\text{ext}$ represents
time-dependent volume and surface forces, which drive the evolution.
We assume that $P_\text{ext}$ is linear and $H^1(\Omega)$-continuous in~$\bu$, and differentiable in $t$
with bounded time derivative.

The last term of~\eqref{eq:potential_energy} implements the restriction that the damage field can only assume values
between $0$ and~$1$.  For a set $\mathcal{K} \subset \R$
we define the indicator functional
\begin{equation*}
 I_{\mathcal{K}} : \R \to \R \cup \{\infty\},
 \qquad
 I_\mathcal{K}(x)
 \colonequals
 \begin{cases}
  0 & \text{if $x \in \mathcal{K}$}, \\
  \infty & \text{otherwise}.
 \end{cases}
\end{equation*}
For a closed, convex, nonempty set $\mathcal{K}$,
the functional $I_\mathcal{K}$ is convex, lower semicontinuous, and proper.
Adding the constraint $d \in [0,1]$ explicitly is not always necessary, as some fracture models lead to evolutions
that satisfy the constraints implicitly.  However, as pointwise bounds come with practically
no cost when using the TNNMG solver, we do include them to extend our range of models.

To make the potential energy $\mathcal{E}$ well defined, we will in general consider
it on the first-order Sobolev space $H^1(\Omega, \R^\Dim) \times H^1(\Omega, \R)$.
Incorporating the boundary conditions leads to the affine subspace
\begin{align*}
  \mathbf{H}_{\bu_0}^1 \times H_{d_0}^1
    \colonequals \Bigl\{\bv \in H^1(\Omega, \R^\Dim) \,\big|\, \bv|_{\Gamma_{D,\bu}} = \bu_0\Bigr\}
    \times \Bigl\{v \in H^1(\Omega) \,\big|\, v|_{\Gamma_{D,d}} = d_0\Bigr\}.
\end{align*}

The second term of the Biot equation~\eqref{eq:biot_equation} is $\partial_{\dot{d}} \mathcal{R}(d,\dot{d})$,
where $\mathcal{R}$ is the dissipation potential
\begin{align}
 \label{eq:dissipation_potential_mielke1}
  \mathcal{R}(d,\dot{d}) & = \int_\Omega I_{[0,\infty)} (\dot{d})\,dV.
\end{align}
It implements the pointwise non-healing condition $\dot{d} \ge 0$ as proposed by \cite{miehe+welschinger+hofacker10a}.
Note that $\mathcal{R}(d,\cdot) : H^1(\Omega) \to [0,\infty]$
is convex and lower semicontinuous, and $\mathcal{R}(d,0) = 0$.
The fact that $\mathcal{R}$ is positively 1-homogeneous in $\dot d$ implies the rate-independence
of the system.
Since the particular functional~\eqref{eq:dissipation_potential_mielke1} only depends on $\dot{d}$
but not on $d$, we will also write $\mathcal{R}(\dot{d}) \colonequals \mathcal{R}(d,\dot{d})$
and $\partial \mathcal{R}(\dot{d}) \colonequals \partial_{\dot{d}}\mathcal{R}(d,\dot{d})$.

\begin{remark}
  \label{rem:variational_formulation}
  Using the definition of the subdifferential $\partial\mathcal{R}$
  it is straightforward to see that the
  Biot equation~\eqref{eq:biot_equation}
  is equivalent to the variational inequality
  \begin{align}
    \label{eq:biot_equation_vi}
    \big \langle D_{(\bu,d)}\mathcal{E}(t,\bu,d),\, (\bv, e)-(\bu,\dot{d})\big \rangle
    + \mathcal{R}(e)
    \geq \mathcal{R}(\dot{d})
    \qquad \forall (\bv,e) \in \mathbf{H}_{\bu_0}^1 \times H_{0}^1.
  \end{align}
  Likewise, it is equivalent to the coupled system
  \begin{alignat*}{2}
    \inner{D_{\bu} \mathcal{E}(t,\bu,d)}{\bv} &= 0
      &\forall \bv &\in \mathbf{H}_{0}^1,\\
    \big \langle D_{d} \mathcal{E}(t,\bu,d), \, e-\dot{d} \big \rangle
      + \mathcal{R}(e)
      &\geq \mathcal{R}(\dot{d}),
    & \qquad \forall e & \in H_{0}^1,
  \end{alignat*}
  where we have denoted by
  $\mathbf{H}_{0}^1 \times H_0^1 \colonequals \mathbf{H}_{\bu_0}^1 \times H_{d_0}^1 - (\bu_0,d_0)$
  the homogeneous space corresponding to $\mathbf{H}_{\bu_0}^1 \times H_{d_0}^1$.
  Since these two notions of solutions are based on first-order
  derivatives of $\mathcal{E}$ and $\mathcal{R}$, they may lead to local minimizers or saddle points of the
  functional $\mathcal{E}$. To overcome this, there are alternative
  so-called \emph{energetic formulations} of the problem.
  In general energetic solutions are solutions of \eqref{eq:biot_equation_vi},
  while the converse is only true for convex $\mathcal{E}$.
  We will not discuss such solutions here,
  but refer to \cite{mielke_roubicek:2015,Thomas_phd:2010},
  where they are discussed for damage-related and more general rate-independent processes.
\end{remark}

\begin{remark}
 In the engineering literature, the same problem is sometimes formulated as
 \begin{equation*}
  \partial_{(\dot{\bu},\dot{d})} \Pi(\dot{\bu},\dot{d};\bu,d) \ni 0,
 \end{equation*}
 with the rate potential
\begin{align*}
  \Pi(\dot\Bu, \dot d; \Bu, d)
  & \colonequals
  \frac{d}{dt} \mathcal{E}_0(\bu,d) + P_\text{ext}(t,\dot{\bu}) + \mathcal{R}(d,\dot{d}),
\end{align*}
(see, e.g., \cite[Section~4.2]{miehe+welschinger+hofacker10a}), where $\mathcal{E}_0(\bu,d)$
denotes the parts of the energy $\mathcal{E}(t,\bu,d)$ that do not explicitly depend on $t$.
If $P_\text{ext}(t,\bu)$ is linear in $\bu$, and if the problem is sufficiently smooth,
then this formulation is equivalent to the Biot equation~\eqref{eq:biot_equation}.
Indeed, we then have
\begin{equation*}
  \Pi(\dot\Bu, \dot d; \Bu, d)
  \colonequals
  D_{(\bu,d)} \mathcal{E}_0(\bu,d)(\dot{\bu}, \dot{d}) + P_\text{ext}(t,\dot{\bu}) + \mathcal{R}(d,\dot{d}),
\end{equation*}
and
\begin{align*}
 \partial_{(\dot{\bu},\dot{d})} \Pi(\dot{\bu},\dot{d};\bu,d)
 & =
 D_{(\bu,d)} \mathcal{E}_0 (\bu,d) + P_\text{ext}(t,\cdot) + \partial_{\dot{d}} \mathcal{R}(d,\dot{d}) \\
 & =
 D_{(\bu,d)} \big[ \mathcal{E}_0 (\bu,d) + P_\text{ext}(t,\bu) \big] + \partial_{\dot{d}} \mathcal{R}(d,\dot{d}).
\end{align*}
Requiring this to contain $0$ is the Biot equation~\eqref{eq:biot_equation}.
\end{remark}

We now discuss the potential energy and the dissipation potential.

\subsection{Degraded elastic energy density}
\label{sec:elastic_energy}

We consider models that behave linearly elastic and isotropic if the material is in an undamaged state.
That is, for the undamaged stored energy density we use the St.\,Venant--Kirchhoff material law,
whose energy density is given by
\begin{equation*}
\psi_0(\bve) = \frac{\lambda}{2}\tr[\bve]^2 + \mu\tr[\bve^2],
\end{equation*}
with Lamé parameters $\mu>0$ and $\lambda > -\tfrac{2}{3}\mu$.
With this choice of parameters the quadratic functional $\psi_0$ is
strongly convex on $\S^\Dim$, the set of real symmetric $m \times m$ matrices.

The undamaged energy density is splitting as
\begin{align*}
  \psi_0(\bve) = \psi_0^+(\bve) + \psi_0^-(\bve)
\end{align*}
into a part $\psi_0^+$ that produces damage
and another part $\psi_0^-$ that does not. The damage-producing part is then scaled by a so-called
degradation function
\begin{equation*}
 g : [0,1] \to [0,1],
\end{equation*}
and the energy density $\psi: \S^\Dim \times [0,1] \to \R$ takes the form
\begin{equation}
\label{eq:damaged_energy_density}
\psi(\bve,d) = [g(d)+k]  \psi^+_{0}(\bve) + (1+k) \psi^-_0(\bve).
\end{equation}
The residual stiffness $k>0$ guarantees a well-posed problem in case of fracture.

Various different degradation functions have appeared in the literature~\cite{steinke+kaliske18,kuhn+schuelter+mueller15,burke_ortner_sueli:2013}.
While the details vary, there appears to be agreement on the following properties:

\begin{assumption}
\label{ass:degradation_function}
The degradation function $g :[0,1] \to [0,1]$ is differentiable, monotone decreasing,
and fulfills $g(0) = 1$ and $g(1) = 0$.
\end{assumption}

\begin{remark}
As an alternative to this assumption, one could consider degradation functions
with $g(0)=1-k$ such that the original energy density $\psi_0$ is exactly
recovered in the undamaged case (see, e.g., \cite{wick17error}).
However, in the typical regime of $k\ll 1$
both will essentially yield the same result. Thus we will follow the more common
approach $g(0)=1$ here, although the proposed methods could equally be applied
to the alternative one.
\end{remark}

Note that several authors require $g'(1)=0$ in order to ensure that the evolution does not
lead to values of $d$ larger than~1. We do not need this assumption here, because the
pointwise constraint $d \le 1$ is enforced explicitly by the energy term~\eqref{eq:potential_energy}.

The following specific degradation functions all fulfill Assumption~\ref{ass:degradation_function}:
\begin{alignat*}{2}
 g_a(d) & = (1-d)^2  & \qquad & \text{(from \cite{bourdin_francfort_marigo:2000})} \\
 g_b(d) & = (1-d)^2 \cdot (2d+1)  && \text{(from \cite{kuhn+schuelter+mueller15})} \\
 g_c(d) & = (1-d)^3 \cdot (3d+1) && \text{(from \cite{kuhn+schuelter+mueller15})} \\
 g_d(d) & = \frac{\exp(b d) - (b ( d - 1) + 1) \exp(b)}{(b - 1) \exp(b) + 1},
  \qquad
  b>0  && \text{(from \cite{steinke+kaliske18})}.
\end{alignat*}
Note that the functions $g_a$ and $g_d$ are strictly convex, but $g_b$ and $g_c$ are not
even convex.
For the rest of the paper we will restrict our considerations
to convex twice continuously differentiable degradation functions $g$.

%\bigskip

Various splittings of $\psi_0$ have been proposed in the literature.
We cover four common strain-based splittings taking the
form \eqref{eq:damaged_energy_density}.%
\footnote{The stress-based splitting of \citet{steinke+kaliske18} is left for future work.}
All those splittings have the property that $\psi_0 = \psi^+_0 + \psi^-_0$,
and we will show that all have the following essential properties:
\begin{enumerate}[label=(P$_{\arabic{*}}$), ref=(P$_{\arabic{*}}$)]
  \item\label{item:density_differentiable}
    $\psi(\bve,\cdot) \in C^2$ for all $\bve \in \S^\Dim$
    and $\psi(\cdot,d) \in C^{1,1}$ for all  $d \in [0,1]$,
    i.e., $\psi(\cdot,d)$ is differentiable with locally Lipschitz continuous derivative.
  \item\label{item:density_semismooth}
    The gradient $\nabla \psi(\cdot,d)$ is semismooth for all $d \in [0,1]$.
  \item\label{item:density_lipschitz}
    The gradient $\nabla \psi(\cdot,d)$ is globally Lipschitz continuous uniformly in $d$,
    i.e., there exists $L \geq 0$ independent of $d$
    such that for all matrices $A, B \in \S^\Dim$ we have
    \begin{equation*}
      |\nabla\psi(A,d)-\nabla\psi(B,d)| \leq L |A-B|_F,
    \end{equation*}
    where $|\cdot|_F$ denotes the Frobenius norm.%
    \footnote{While we only require $L\geq 0$ in~\ref{item:density_lipschitz},
    the strong convexity~\ref{item:density_strongly_convex} in fact implies $L>0$.}
  \item\label{item:density_strongly_convex}
    $\psi(\cdot,d) : \S^\Dim \to \R$ is strongly convex uniformly in $d$,
    i.e., there exists $\eta > 0$ independent of $d$
    such that for all matrices $A, B \in \S^\Dim$ we have
    \begin{equation*}
      \psi\big(tA+(1-t)B,d\big)\leq t \psi(A,d) + (1-t) \psi(B,d)-{\frac {1}{2}}  \eta t(1-t)|A-B|_F^2.
    \end{equation*}
  \item\label{item:density_coercive}
    $\psi(\cdot,d)$ is coercive uniformly in $d$ in the sense
    that there exists $C > 0$ independent of $d$ such that
    $\psi(\bve,d) \geq C |\bve|_F^2$.% The comment here prevents spurious whitespace in front of the footnote mark.
    \footnote{Under the additional assumption that $0 = \psi(0,d) \leq \psi(\bve,d)$ holds for all $\bve$ and $d$
    one can show that \ref{item:density_strongly_convex} implies \ref{item:density_coercive}.}
\end{enumerate}
We remind the readers that the gradient $\nabla \psi(\cdot,d)$ is called semismooth if
for any point $A \in \S^\Dim$ and any direction $V \in \S^\Dim$ the limit
\begin{align*}
  \lim_{n \to \infty}G_n V_n
\end{align*}
exists and is unique for all sequences $V_n$ and $G_n$ with
$V_n \to V$ and $G_n \in \partial (\nabla \psi(\cdot,d))(A+t_n V_n)$
for a sequence $t_n \searrow 0$.
The set $\partial (\nabla \psi(\cdot,d))(A)$ denotes
Clarke's generalized Jacobian of the locally Lipschitz
continuous map
$\nabla \psi(\cdot,d) : \S^\Dim \to \S^\Dim$ at $A \in \S^\Dim$
(cf.\,\cite{QiSun:1993}).
Notice that the strong convexity~\ref{item:density_strongly_convex}
implies strong monotonicity of $\nabla \psi$, i.e.,
\begin{align*}
  \inner{\nabla \psi(A)-\nabla \psi(B)}{A-B} \geq \eta |A-B|_F^2.
\end{align*}
Here and in the following we denote by
$\inner{\cdot}{\cdot}:\Vcal^* \times \Vcal \to \R$ the duality pairing of a vector space $\Vcal$.

For the splittings considered in the following
we will only prove
\ref{item:density_differentiable} and \ref{item:density_semismooth}
directly,
and show that the simplified assumptions of the following
lemma hold true.  This then implies
\ref{item:density_lipschitz}, \ref{item:density_strongly_convex},
and \ref{item:density_coercive}.
\begin{lemma}
  \label{lemma:density_simplified}
  Let $\psi_0^+$ and $\psi_0^-$ be convex, non-negative,
  and differentiable with Lipschitz continuous gradients
  $\nabla\psi_0^+$ and $\nabla\psi_0^-$.
  Then $\psi$ satisfies
  \ref{item:density_lipschitz}, \ref{item:density_strongly_convex},
  and \ref{item:density_coercive}.
\end{lemma}
\begin{proof}
  Let $L^+$ and $L^-$ be the Lipschitz constants of
  $\nabla \psi_0^+$ and $\nabla \psi_0^-$, respectively.
  Then $\nabla\psi(\cdot,d)$ is Lipschitz continuous
  with uniform Lipschitz constant
  $(1+k)(L^++L^-)$, because $g(d)+k \leq 1+k$.

  To show strong convexity, we first note that $\psi_0 = \psi_0^+ + \psi_0^-$
  is strongly convex on $\S^\Dim$
  with a modulus $\eta>0$ independent of $d$.
  Now consider the function
  \begin{align*}
    \bve \mapsto
    &
        \psi(\bve,d) - k \psi_0(\bve) = g(d)\psi_0^+(\bve) + \psi_0^-(\bve).
  \end{align*}
  Since this is a weighted sum of two convex functions
  $\psi_0^+$ and $\psi_0^-$ with non-negative weights
  $g(d) \geq 0$ and $1$, it is itself convex.
  Thus, as a sum of this convex function
  and the strongly convex functions $C\psi_0$,
  the function $\psi(\cdot,d)$ is itself strongly
  convex and inherits the convexity modulus $k\eta$ of $k\psi_0$.
  Finally, we note that with the same $\eta$ we have
  \begin{equation*}
    \psi(\bve,d)  \geq
    k \psi_0(\bve) \geq k\frac{\eta}2|\bve|_F^2.
  \qedhere
  \end{equation*}
\end{proof}

Despite those strong properties of $\psi(\cdot,d)$ we note that
$\psi(\bve,d)$ is not convex in $d$ and $\bve$ together for any of
the splittings considered below.

\subsubsection{Isotropic splitting}

In this model, any strain will lead to damage.  The splitting is therefore
\begin{equation}
 \label{eq:symmetric_damaged_energy}
 \psi_0^+(\bve) = \psi_0(\bve),
 \qquad
 \psi_0^-(\bve) = 0.
\end{equation}
Without proof, we note the following simple properties of the
energy density $\psi$ defined by \eqref{eq:damaged_energy_density} and this splitting:

\begin{lemma}
  The energy density $\psi$ defined in \eqref{eq:damaged_energy_density}
  with the isotropic splitting~\eqref{eq:symmetric_damaged_energy} has the
  properties \ref{item:density_differentiable}--\ref{item:density_coercive}.
  Furthermore $\psi(\cdot,d)$ has the stronger property
  that it is in $C^\infty$ and quadratic for all $d \in [0,1]$.
\end{lemma}

\subsubsection{Volumetric decompositions}
\label{sec:volumetric_decompositions}

The isotropic model is unphysical, because it produces fracturing for all kinds of
strain.  In~\cite{lancioni_royer_carfagni:2009}, \citeauthor{lancioni_royer_carfagni:2009}
obtained better results by letting only the deviatoric strain contribute to the
degradation.  They introduced the split
\begin{equation*}
 \psi_0^+(\bve) = \psi_0(\operatorname{dev} \bve),
 \qquad
 \psi_0^-(\bve) = \psi_0(\operatorname{vol} \bve),
\end{equation*}
with the deviatoric--volumetric strain splitting
\begin{equation*}
 \operatorname{vol} \bve \colonequals \frac{\tr \bve}{\Dim} I,
 \qquad
 \operatorname{dev} \bve \colonequals \bve - \operatorname{vol} \bve.
\end{equation*}
With these definitions, the energies are
\begin{equation}
\label{eq:isotropic_volumetric_split}
 \psi_0^+(\bve) = \Big( \frac{\mu}{\Dim} + \frac{\lambda}{2} \Big)(\tr \bve)^2,
 \qquad
 \psi_0^-(\bve) = \mu \Big(\bve^2 - \frac{1}{\Dim} (\tr \bve)^2\Big) = \mu \operatorname{dev} \bve : \operatorname{dev} \bve.
\end{equation}

\begin{lemma}
  The energy density $\psi$ defined in \eqref{eq:damaged_energy_density}
  with the isotropic volumetric splitting~\eqref{eq:isotropic_volumetric_split} has the
  properties \ref{item:density_differentiable}--\ref{item:density_coercive}.
  Furthermore $\psi(\cdot,d)$ has the stronger property
  that it is in $C^\infty$ and quadratic for all $d \in [0,1]$.
\end{lemma}

\begin{proof}
  $C^\infty$-smoothness and thus \ref{item:density_differentiable} and \ref{item:density_semismooth}
  are straightforward.
  The fact that $\psi_0^+$ and $\psi_0^-$ are quadratic, convex, and non-negative
  allows to derive \ref{item:density_lipschitz}, \ref{item:density_strongly_convex},
  and \ref{item:density_coercive} from Lemma~\ref{lemma:density_simplified}
  and implies that $\psi(\cdot,d)$ is also quadratic.
\end{proof}

The decomposition of \citeauthor{lancioni_royer_carfagni:2009} is still isotropic. \citet{amor_marigo_maurini:2009} proposed to
only degrade the expansive part of the volumetric strain.
Using the ramp functions
\begin{align*}
  \langle x\rangle_+ \colonequals \max \{0,x\}, \qquad
  \langle x\rangle_- \colonequals \min\{0,x\}
\end{align*}
that provide the decompositions
$x = \langle x \rangle_+ + \langle x \rangle_-$
and $x^2 = \langle x \rangle_+^2 + \langle x \rangle_-^2$,
they proposed the energy split
\begin{equation}
\label{eq:anisotropic_volumetric_split}
 \psi_0^+(\bve) = \Big( \frac{\mu}{\Dim} + \frac{\lambda}{2} \Big)\langle \tr \bve\rangle_+^2,
 \qquad
 \psi_0^-(\bve) = \Big( \frac{\mu}{\Dim} + \frac{\lambda}{2} \Big)\langle \tr \bve\rangle_-^2 + \mu \operatorname{dev} \bve : \operatorname{dev} \bve,
\end{equation}
where only the tensile volumetric strain contributes to damage.

\begin{lemma}
  The energy density $\psi$ defined in \eqref{eq:damaged_energy_density}
  with the anisotropic volumetric splitting~\eqref{eq:anisotropic_volumetric_split} has the
  properties \ref{item:density_differentiable}--\ref{item:density_coercive}.
  Furthermore
  $\psi(\cdot,d)$ is not $C^2$, unless  $g(d) = 1$.
\end{lemma}

\begin{proof}
  We first note that the squared ramp functions $\langle \cdot \rangle_{\pm}^2$ are convex, $C^{1,1}$
  with derivatives having a global Lipschitz constant $2$, and piecewise $C^2$
  (in the sense of \cite[Definition 2.19]{Ulbrich2002}).
  Hence the functions $\psi_0^{\pm}$ are also $C^{1,1}$
  with globally Lipschitz gradients and piecewise $C^2$,
  which shows \ref{item:density_differentiable} and
  (using Lemma~\ref{lemma:density_simplified}) \ref{item:density_lipschitz}.
  Being piecewise $C^2$ implies semismoothness
  \ref{item:density_semismooth} of $\nabla \psi(\cdot,d)$ \cite[Proposition~2.26]{Ulbrich2002}.
  Noting that $\mu/\Dim + \lambda/2 >0$,
  convexity of the squared ramp functions
  furthermore implies that the functions $\psi_0^{\pm}$ are also convex
  and non-negative, which by
  Lemma~\ref{lemma:density_simplified}
  provides \ref{item:density_strongly_convex}
  and \ref{item:density_coercive}.

  For $g(d)+k=1$ the functional $\psi(\cdot,d)$ is quadratic and thus $C^2$.
  In the case $g(d)\neq1$, if $\psi(\cdot,d)$ would be $C^2$, then
  the function $t \mapsto \psi(tI,d)$ 
  would also be~$C^2$. However,
  this function takes the form
  \begin{align*}
    \psi(tI,d)
    = \Big( \frac{\mu}{\Dim} + \frac{\lambda}{2} \Big) \Dim^2 t^2
    \begin{cases}
      g(d)+k &\text{ if $t\geq 0$},\\
      1 +k &\text{ if $t<0$}
    \end{cases}
  \end{align*}
  and is thus piecewise quadratic but not $C^2$ in $t=0$. \qedhere
\end{proof}

\subsubsection{Spectral decomposition}

A more elaborate nonlinear splitting
separating the tensile and compressive parts of the elastic energy
was introduced in~\cite{miehe+welschinger+hofacker10a}.
To define this splitting it is convenient to introduce the ordered
eigenvalue function $\Eigfunc: \S^\Dim \to \R^\Dim$
on the space $\S^\Dim$ of symmetric $\Dim \times \Dim$ matrices,
mapping any symmetric matrix $M$
to the vector $\Eigfunc(M) \in \R^\Dim$ containing its
eigenvalues in ascending order.
Using the ramp functions
the tensile and compressive energies $\psi_0^+$
and $\psi_0^-$ are then defined as
\begin{equation}
\label{eq:spectral_energy_split}
\psi^{\pm}_0(\bve)
\colonequals
  \frac{\lambda}{2} \Big\langle \sum_{i=1}^\Dim \Eigfunc(\bve)_i \Big\rangle_\pm^2
  + \mu \sum_{i=1}^\Dim\langle \Eigfunc(\bve)_i\rangle_\pm^2.
\end{equation}
Note that this indeed defines a splitting $\psi_0 = \psi^+_0 + \psi^-_0$.
For this splitting we will make the additional assumption
that $\lambda\geq 0$.

To quantify the properties of $\psi(\bve,d)$ with respect
to the strain tensor $\bve$ we use the theory of spectral functions~\cite{lewis:1996}.
To this end we note that we can write $\psi^\pm_0$ as
\begin{align*}
  \psi^\pm_0  = \widehat{\psi}^\pm_0 \circ \Eigfunc: \S^\Dim \to \R
\end{align*}
with
\begin{equation*}
  \widehat{\psi}^{\pm}_0(\Eig)
  \colonequals
  \frac{\lambda}{2}\Big\langle \sum_{i=1}^\Dim \Eig_i \Big\rangle_\pm^2
  + \mu \sum_{i=1}^\Dim \langle \Eig_i \rangle_\pm^2.
\end{equation*}
The functions
$\widehat{\psi}^\pm_0$ are \emph{symmetric} in the sense
that $\widehat{\psi}^\pm_0(\Eig)$ does not depend on the order
of the entries of $\Eig \in \R^\Dim$.
Having this form we can infer properties of
the functions $\psi_0^\pm = \widehat{\psi}^\pm_0 \circ \Eigfunc$
from properties of the symmetric functions
$\widehat{\psi}^\pm_0$.

\begin{lemma}
\label{lem:smoothness_of_elastic_energy_density}
  Let $\lambda \geq 0$.
  Then the energy density $\psi$ defined in \eqref{eq:damaged_energy_density}
  with the spectral splitting~\eqref{eq:spectral_energy_split} has the
  properties \ref{item:density_differentiable}--\ref{item:density_coercive}.
  Furthermore
  $\psi(\cdot,d)$ is not $C^2$, unless $g(d) = 1$.
\end{lemma}
\begin{proof}
  We will first show \ref{item:density_differentiable}--\ref{item:density_coercive}.
  An essential ingredient is
  that the squared ramp functions $\langle \cdot \rangle_\pm^2$
  are non-negative, piecewise quadratic, and convex.

  \ref{item:density_differentiable}
  The squared ramp functions $\langle \cdot \rangle_\pm^2$
  and thus $\widehat{\psi}^\pm_0$ are $C^{1,1}$.
  Now \cite[Proposition 4.3]{QiYang2003} shows that
  the spectral functions
  $\psi^\pm_0 =\widehat{\psi}^\pm_0 \circ \Eigfunc$ are
  also $C^{1,1}$. Hence the same applies to $\psi(\cdot,d)$.

  \ref{item:density_semismooth}
  The squared ramp functions $\langle \cdot \rangle_\pm^2$
  are piecewise $C^2$ functions.
  Hence the gradients $\nabla \widehat{\psi}^\pm_0$ are piecewise $C^1$
  functions (in the sense of \cite[Definition 2.19]{Ulbrich2002})
  and thus semismooth \cite[Proposition 2.26]{Ulbrich2002}.
  Now \cite[Proposition 4.5]{QiYang2003} provides semismoothness
  of $\nabla \psi^\pm_0$ and thus of $\nabla\psi(\cdot,d)$.

  \ref{item:density_lipschitz}
  Since the functions $\widehat{\psi}^\pm_0$ are
  piecewise quadratic and $C^{1,1}$ the gradients
  $\nabla \widehat{\psi}^\pm_0$ are globally Lipschitz continuous.
  Now Corollary~43 of~\cite{baes:2007} provides
  global Lipschitz continuity of the
  gradients $\nabla \psi^\pm_0$ of the
  spectral functions $\psi^\pm_0$
  in the more general context of Euclidean Jordan algebras
  (which includes the special case of symmetric matrices).
  In fact, the Lipschitz constant of $\nabla \widehat{\psi}^\pm_0$
  equals the one for $\nabla \psi^\pm_0$
  if $\S^\Dim$ is equipped with the Frobenius norm.
  Using Lemma~\ref{lemma:density_simplified}
  this implies uniform Lipschitz continuity of $\psi(\cdot,d)$.

  \ref{item:density_strongly_convex},\ref{item:density_coercive}
  %    As a composition of sum functions and of the
  %    convex continuous function $\langle \cdot \rangle_\pm^2$
  %    the functions $\widehat{\psi}^\pm_0$ are convex
  %    and lower semicontinuous on $\R^\Dim$.
  %    By~\cite[Corollary\,2.4]{lewis:1996} the functions
  %    $\psi^\pm_0$ are then convex and lower semicontinuous
  %    on $\S^\Dim$.
  %    To show even strict convexity, observe that $\psi$ as a function on $\R^3$ is
  %    strictly convex on $\R^3$.  In particular, it is therefore \emph{essentially strictly convex}.
  %    Therefore (essential) strict convexity of $\psi$ as a function on $\S^d$ follows from
  %    \cite[Corollary\,3.3]{lewis:1996}.
  Since the functions $\widehat{\psi}^\pm_0$ are
  weighted sums of convex, non-negative squared ramp functions
  with non\-negative weights, they are convex
  and non-negative themselves.
  Convexity of the functions $\psi_0^\pm$
  then follows from~\cite[Theorem\,41]{baes:2007}
  while non-negativity of those functions is trivial.
  Now Lemma~\ref{lemma:density_simplified} provides
  \ref{item:density_strongly_convex} and \ref{item:density_coercive}.

  To characterize second order differentiability
  of $\psi(\cdot,d)$ we first consider $g(d)=1$.
  Then $\psi(\cdot,d)$ coincides with the
  quadratic function $(1+k)\psi_0 = (1+k)\psi^+_0 + (1+k)\psi^-_0$ and is thus~$C^2$.
  In the case $g(d)\neq 1$, if $\psi(\cdot,d)$ would be $C^2$, then
  the function $t \mapsto \psi(tE,d)$ for the fixed matrix $E$ with
  $E_{ij} = \delta_{1i}\delta_{1j}$ would also be $C^2$. However,
  this function takes the form
  \begin{align*}
    \psi(tE,d)
    = \Bigl(\frac{\lambda}{2}+1\Bigr)t^2
    \begin{cases}
      g(d)+k &\text{ if $t\geq 0$},\\
      1+k &\text{ if $t<0$},
    \end{cases}
    %      = \Bigl(\frac{\lambda}{2}+\mu \Bigr)
    %        \Bigl(
    %          (g(d)+k) \langle t \rangle_+^2 + \langle t \rangle_-^2 \Bigr)
  \end{align*}
  and is thus piecewise quadratic but not $C^2$ in $t=0$. \qedhere
\end{proof}

\begin{remark}
  \label{rem:generalized_jacobian}
  One can show that $\S^\Dim$ decomposes into finitely many disjoint subsets~$\mathcal{A}_i$
  such that $\psi(\cdot,d)$ is twice continuously differentiable
  in the interior of each of these sets.
  A matrix $\bve \in \S^\Dim$ is in the intersection of several $\overline{\mathcal{A}}_i$
  if it either has an eigenvalue $\Eigfunc(\bve)_i = 0$ or if $\tr \bve=0$.
  While $\nabla \psi(\cdot,d)$ is not differentiable at those points,
  there are still generalized second-order derivatives.
  For example, the generalized Jacobian in the sense of Clarke contains
  the derivatives of $\nabla \psi(\cdot,d)$ with respect to all the adjacent
  sets $\mathcal{A}_i$.
  Semismoothness essentially means that such generalized derivatives
  provide an approximation that can be exploited in a generalized Newton method.
\end{remark}

\begin{remark}
  The additional assumption $\lambda \geq 0$ is essential for
  convexity of $\psi(\cdot,d)$.
  To see this we consider for $\Dim=2$
  the line segment
  \begin{align*}
    \bigl\{D(t) = \operatorname{diag}(-1,t) \, \big| \, t \in (0,1) \bigr\} \subset \S^2.
  \end{align*}
  Then, along this line segment, $\psi(\cdot,1)$ is quadratic and takes the form
  \begin{align*}
    \psi(D(t),1)
    &= k \mu t^2 + (1+k)\Bigl(\frac{\lambda}{2}(t-1)^2 + \mu\Bigr)\\
    &= \Bigl(k \mu + (1+k)\frac{\lambda}{2}\Bigr) t^2 + (1+k)\Bigl(-\lambda t + \frac{\lambda}{2} +  \mu \Bigr)
  \end{align*}
  which is strictly concave for $\lambda<0$ and sufficiently small $k\ll 1$.
\end{remark}

\subsection{Crack surface density}
\label{sec:crack_surface_density}
The crack surface density function per unit volume of the solid is typically of the form \cite{pham_amor_marigo_maurini:2011}
\begin{equation*}
 \gamma(d,\nabla d)
 \colonequals
  \frac{1}{4 c_w}\Big(\frac{w(d)}{l} + l \abs{\nabla{d}}^2 \Big),
\end{equation*}
with parameters $c_w$ and $l$, and a parameter function $w : [0,1] \to [0,1]$.
Motivated by the seminal work of Modica and Mortola~\cite{ModicaMortola77a, ModicaMortola77b},
the associated crack surface functional
\begin{align*}
  \Ccal(d) \colonequals \int_\Omega \gamma(d,\nabla d) \,dx
\end{align*}
is called a Modica--Mortola functional.
The internal length scale parameter $l$ controls the size of the diffusive zone between a completely intact and a completely damaged material.
For $l \to 0$ the regularized crack surface yields a sharp crack topology in the sense of $\Gamma$-convergence~{\cite{ModicaMortola77a, ModicaMortola77b}.
For a given function $w$, the normalization constant $c_w$
must be chosen such that the integral of $\gamma(d,\nabla d)$ over the fractured domain converges to the surface measure of the crack set as $l \to 0$.

The function $w(d)$ models the local fracture energy. Two types of local crack density functions appear in the literature. Double-well potentials (as briefly reviewed in~\cite{ambati_gerasimov_delorenzis:2015})
provide an energy barrier between broken and unbroken state, but will be disregarded here.
Instead, we only consider single-well potentials $w$, which grow monotonically
from the intact state $w(0)=0$ to the damaged state $w(1)=1$.
For such potentials the normalization constant is given by
\begin{align}\label{eq:cw_value}
  c_w \colonequals \int_0^1 \sqrt{w(s)} \,ds.
\end{align}
While this follows from the $\Gamma$-convergence proof
(see, e.g., \cite{Alberti2000}), the proper constant can also be computed
by elementary means:
Consider an open set $\Gamma \subset \R^{\Dim-1}$ of measure $|\Gamma|$
embedded in a domain $\Omega = \R \times \Gamma$.  We identify $\Gamma$ with
the set $\{x \in \Omega \; :\; x_1 = 0 \}$, and interpret it as a crack in $\Omega$.
To approximate $\Gamma$ by a damage field $d : \Omega \to [0,1]$,
assume that $d$ minimizes the crack surface energy
\begin{align*}
  \Ccal(d) \colonequals \int_\Omega \gamma(d,\nabla d) \,dx,
\end{align*}
subject to the crack conditions $d(0,\xi)=1$ and $d(\pm\infty,\xi)=0$
for all $\xi \in \Gamma$.
We find that $d$ is given by $d(x) = \hat{d}(|x_1|/l)$ where
$\hat{d}: [0,\infty)$ solves
the normalized scalar Euler--Lagrange equation
\begin{align*}
  w'(\hat{d}) - 2 \hat{d}''=0 \qquad \text{ in } (0,\infty),
\end{align*}
with boundary condition $\hat{d}(0) = 1$, or, equivalently,
the initial value problem,
\begin{align}\label{eq:crack_profile_ode}
  \hat{d}'(s) = - \sqrt{w(\hat{d}(s))} \qquad s>0, \qquad \hat{d}(0)=1.
\end{align}
Using \eqref{eq:crack_profile_ode}
we find that the crack surface energy of the minimizer $d$ is
\begin{align*}
  \Ccal(d)
%  = \int_\Omega \gamma(d,\nabla d) dx for a vector space $\Vcal$
%    &= 2|\Gamma| \int_0^\infty \gamma(d(s,0),\nabla d(s,0) ds
    = \frac{2|\Gamma|}{4 c_w} \int_0^\infty w(\hat{d}) + |\hat{d}'|^2 ds 
%    \\
%    &= \frac{|\Gamma|}{2c_w} \int_0^\infty \Bigl(\sqrt{w(\hat{d})} + \hat{d}'\Bigr)^2 - 2\sqrt{w(\hat{d})}\hat{d}'\,ds
%    \\
    &= -\frac{2|\Gamma|}{4 c_w} \int_0^\infty 2\sqrt{w(\hat{d})}\hat{d}'\,ds
    = \frac{|\Gamma|}{c_w} \int_0^1 \sqrt{w(s)}\,ds.
\end{align*}
Thus $c_w$ has to be selected as in \eqref{eq:cw_value} to ensure that
the crack surface energy scales like the limit crack surface measure $|\Gamma|$.

Note that the function $\hat{d}$ is the $w$-dependent crack profile,
rescaled by the length parameter $l$.
In order to relate the length parameter to the actual crack width,
we use the cone construction from~\cite{miehe+hofacker+welschinger10}
and define the crack width to be the average width of a tangential finite cone
fitted to the crack profile and the zero function (Figure~\ref{fig:crack_cone_width}).
From the initial condition $\hat{d}(0) = 1$,
the crack profile equation \eqref{eq:crack_profile_ode},
and the normalization $w(1)=1$ it follows that $\hat{d}'(0) = -1$.
Thus---for the normalized solution $\hat{d}$---the cone has width $2$ at its base
and average width $1$.
Hence the average crack cone width of the rescaled solution $d$ is $l$.
\begin{figure}
  \begin{center}
    \input{crack_cone.pgf}
  \end{center}
  \caption{Crack width definition for the AT-2 functional}
  \label{fig:crack_cone_width}
\end{figure}

We will focus on the two widely used potentials
\begin{alignat*}{2}
 w(d) & =d,
  & \qquad &
  c_w = \frac{2}{3} \\
 \intertext{and}
 w(d) & = d^2,
  & &
  c_w = \frac{1}{2}.
\end{alignat*}
They are referred to in the literature as Ambrosio--Tortorelli (AT) functionals of type~1 and~2, respectively.
The corresponding crack profiles are given by
\begin{align*}
  \hat{d}_{\text{AT-1}}(s) =
  \begin{cases}
    \left(1-\frac{s}{2}\right)^2 & \text{if $s<2$},\\
    0 & \text{otherwise}
  \end{cases}
  \qquad \text{ and } \qquad
  \hat{d}_{\text{AT-2}}(s) = \exp(-s).
\end{align*}

% Kuhn schreibt: Most commonly, the convex quadratic function
% with $\beta = -1$ is used, see e.g. [19–21]. But also the linear case
% $\beta = 0$ is found e.g. in [22,23].}

Some authors like \cite{kuhn+schuelter+mueller15} prefer $w(d) = d^2$
because it has a local minimizer at $d=0$.
Thus, in the absence of mechanical strain,
the unfractured solution $d \equiv 0$ is a minimizer of the total energy.
As a result, no additional constraints need to be
applied to ensure that $d \ge 0$.  However, this argument becomes void when solver
technology is available that can handle the explicit constraints $0 \le d \le 1$.
In contrast, for the AT-1 functional we have $w' \neq 0$ in the intact state $d=0$.
Together with the constraint $d \in [0,1]$ this leads to a threshold,
i.e., a minimum load required to cause damage~\cite{pham_amor_marigo_maurini:2011}.
A numerical comparison of the AT-1 and AT-2 functionals can be found in~\cite{burke_ortner_sueli:2013}.

\citet{kuhn+schuelter+mueller15} proposed to regard the Ambrosio--Tortorelli functionals
as special instances of the general family defined by
\begin{equation}
\label{eq:monotone_local_crack_density}
 w(d) = (1+\beta(1-d))d,
\end{equation}
with $\beta \in [-1,1]$.  The Ambrosio--Tortorelli functionals are obtained by setting
$\beta = 0$ for AT-1 and $\beta = -1$ for AT-2.
Further choices of $w$ are proposed in \cite{pham_amor_marigo_maurini:2011},
which also contains a detailed stability analysis for one-dimensional problems.

We note the following properties of the functional $w$ in~\eqref{eq:monotone_local_crack_density}:
\begin{lemma}
 The function $w$ given in~\eqref{eq:monotone_local_crack_density} has the following properties:
 \begin{enumerate}
  \item It fulfills $w(0)=0$ and $w(1)=1$.
  \item It is strictly monotone increasing on $[0,1]$ for all $\beta \in [-1,1]$.
  \item It is convex for all $\beta \le 0$, and strictly convex for all $\beta < 0$.
 \end{enumerate}
\end{lemma}

For the rest of the paper we will assume the $w(\cdot)$ takes the form \eqref{eq:monotone_local_crack_density}
with $\beta\leq 0$ such that $w(\cdot)$ is guaranteed to be convex and quadratic.

\section{Discretization and the algebraic increment potential}
\label{sec:time_discrete_model}

We use a fully implicit discretization in time, and Lagrange finite elements for discretization
in space. By using a fully implicit
time discretization we retain the variational structure of the problem.
Most of this section is spent investigating the properties of the increment functional.

\subsection{Time discretization}
\label{sec:time_discretization}

It is shown in~\cite{mielke_roubicek:2015} that there is a natural time discretization for
\eqref{eq:biot_equation} that consists of sequences of minimization problems.
To simplify the presentation we will not derive the time discretization
from an energetic formulation of the time-dependent problem
(cf.\ Remark~\ref{rem:variational_formulation}) but first discretize the variational
formulation \eqref{eq:biot_equation_vi} and then reformulate
the time-discrete problem as a sequence of minimization problems.
Let the time interval $[0,T]$ be subdivided by time points $t_n$, $n=0,1,2,\dots$
and denote by $(\bu_n, d_n) \in \mathbf{H}_{\bu_0} \times H_{d_0}$
the discrete approximation of $(\bu(t_n), d(t_n))$.
Approximating the time derivative $\dot{d}(t_{n+1})$ by the backward difference quotient
$(d_{n+1} - d_n)/\tau_n$ for the time step size $\tau_n \colonequals t_{n+1} - t_n$,
and inserting this into the time-continuous variational inequality~\eqref{eq:biot_equation_vi}
we obtain the time-discrete variational inequality for $(\bu_{n+1}, d_{n+1})$
\begin{multline*}
  \Big \langle D_{(\bu_{n+1},d_{n+1})}
    \mathcal{E}(t_{n+1},\bu_{n+1},d_{n+1}),\,(\bv, e)-\big(\bu_{n+1},(d_{n+1} - d_n)/\tau_n \big) \Big \rangle
  \\
  + \mathcal{R}(e)
  - \mathcal{R}\big((d_{n+1} - d_n)/\tau_n \big)
  \geq 0
  \qquad \forall (\bv,e) \in \mathbf{H}_{\bu_0}^1 \times H_{0}^1.
\end{multline*}
Testing with
\begin{align*}
  (\bv,e)=\Bigl(
    \bu_{n+1} + \frac{1}{\tau_n}(\hat{\bv}-\bu_{n+1}),
    \frac{1}{\tau_n}(\hat{e} - d_n)\Bigr)
\end{align*}
for $(\hat{\bv}, \hat{e}) \in \mathbf{H}_{\bu_0}^1 \times H_{d_0}^1$,
using the fact that $\mathcal{R}$ is $1$-homogeneous,
and relabeling $(\hat{\bv}, \hat{e})$ to $(\bv,e)$
yields
\begin{multline}
  \label{eq:time_discret_vi}
  \Big \langle D_{(\bu_{n+1},d_{n+1})} \mathcal{E}(t_{n+1},\bu_{n+1},d_{n+1}),\,(\bv, e)-(\bu_{n+1},d_{n+1}) \Big \rangle
  \\
  + \mathcal{R}(e - d_n)
  - \mathcal{R}(d_{n+1} - d_n)
  \geq 0
  \qquad \forall (\bv,e) \in \mathbf{H}_{\bu_0}^1 \times H_{d_0}^1.
\end{multline}
This is the first-order optimality system for the minimization problem
\begin{equation}
\label{eq:increment_problem}
 (\Bu_{n+1}, d_{n+1})
 \colonequals
  \argmin_{(\tilde{\Bu},\tilde{d}) \in \mathbf{H}^1_{\bu_0}\times H^1_{d_0}} \Pi^\tau_{n+1}(\tilde{\Bu},\tilde{d}),
\end{equation}
with the increment potential
\begin{align*}
 \Pi^\tau_{n+1}(\Bu,d)
 &\colonequals \mathcal{E}(t_{n+1},\Bu, d) + \mathcal{R}(d - d_n)
 \\
 &=
    \int_\Omega \psi(\bve(\Bu),d)\,dV + \int_\Omega g_c\gamma(d,\nabla d)\,dV
    + P_\text{ext}(t_{n+1},\bu)
    + \int_\Omega I_{[d_n,1]} (d) \,dV.
\end{align*}
  Although the variational inequality~\eqref{eq:time_discret_vi} is not equivalent to
  the minimization problem~\eqref{eq:increment_problem} because $\mathcal{E}$ is not convex,
  we will use the minimization formulation
  in the following.

Note that the time step size does not appear in this functional, which means that the model
is rate-independent. Note also that the increment potential depends on the previous time step
only through the indicator functional.

\begin{lemma}
  \label{lemma:time_discrete_coercivity}
  Assume that $\Gamma_{D,\bu}$ is non-trivial in the sense
  that its $\Dim-1$-dimensional Hausdorff-measure is positive.
  Then the functional $\Pi^\tau_{n+1}$ is coercive on $\mathbf{H}_{\bu_0}^1 \times H_{d_0}^1$.
\end{lemma}
\begin{proof}
  Using the uniform coercivity \ref{item:density_coercive}
  of $\psi(\cdot,d)$, $w(1)>0$, and $w(d)\geq 0$ for $d \in [0,1]$ we get
  \begin{align*}
    \int_\Omega \psi(\bve(\Bu),d)\,dV
    + \int_\Omega g_c\gamma(d,\nabla d)\,dV
    \geq C \int_\Omega |\bve(\Bu)|_F^2 + |\nabla d|^2 \,dV
  \end{align*}
  for some constant $C>0$.
  Using Korn's inequality for $\Bu$,
  the Poincar\'e inequality for $d$,
  and the fact that
  $P_\text{ext}(t_{n+1},\bu)$ grows at most linearly we get
  for another constant $C>0$
  \begin{align*}
    \Pi^\tau_{n+1}(\Bu,d)
    &\geq C \biggl(\|\bu\|_1^2 + \|d\|_1^2 - 1 - \Big(\int_\Omega d \,dV\Bigr)^2\biggr)
            + \int_\Omega I_{[d_n,1]} (d) \,dV\\
    &\geq C \biggl(\|\bu\|_1^2 + \|d\|_1^2 - 1 - |\Omega|^2\biggr),
  \end{align*}
  where we have used that the constraint $d\in [0,1]$ implies
  $|\int_\Omega d \,dV|\leq |\Omega|$ in the second inequality.
\end{proof}

\begin{lemma}
  The functional $\Pi^\tau_{n+1}$ is weakly lower semicontinuous on
  $\mathbf{H}^1_{\bu_0} \times H^1_{d_0}$.
\end{lemma}
\begin{proof}
  Since weak lower semicontinuity of the other terms in $\Pi^\tau_{n+1}$ follows
  from convexity and lower semicontinuity of the integrands,
  we only need to consider the non-convex term
  \begin{align}
  \label{eq:wlsc_nonconvex_term}
    \int_\Omega \psi(\bve(\Bu),d) + I_{[d_n,1]}(d)\,dV.
  \end{align}
  To this end we note that \eqref{eq:wlsc_nonconvex_term} can be written as
  $J(\Bu,d,\nabla \Bu)$ for
  \begin{align*}
    J(\Bu,d,\xi) \colonequals \int_\Omega F\big(x,(\Bu(x),d(x)),\xi\big)\,dV
  \end{align*}
  and the density
  $F: \Omega \times (\R^\Dim \times \R) \times \R^{\Dim \times \Dim} \to \R \cup \{ \infty \}$ given by
  \begin{align*}
    F\big(x,(\Bu,d),\xi\big) = \psi(\tfrac{1}{2}(\xi + \xi^T),d) + I_{[d_n,1]}(d).
  \end{align*}
  Since $F$ is a Carathéodory function, non-negative (and thus uniformly bounded from below),
  and convex in~$\xi$ for all $(x,(\Bu,d)) \in \Omega \times (\R^\Dim \times \R)$,
  it satisfies the assumptions of Theorem~3.4 in~\cite{Dacorogna1989}.

  Now let $(\Bu^\nu,d^\nu) \rightharpoonup (\bu,d)$ be a weakly convergent sequence
  in $\mathbf{H}^1_{\bu_0} \times H^1_{d_0}$.
  Then, by the compactness of the embedding into $L^2(\Omega,\R^\Dim \times \R)$ we get
  \begin{align*}
    (\Bu^\nu,d^\nu) \to (\bu,d) \qquad \text{ in } L^2(\Omega, \R^\Dim \times \R).
  \end{align*}
  Furthermore, the $H^1(\Omega,\R^\Dim \times \R)$-weak convergence of $(\Bu^\nu,d^\nu)$ implies
  $L^2(\Omega, \R^{\Dim\times \Dim})$-weak convergence of $\nabla \Bu^\nu$
  \begin{align*}
    \xi^\nu \colonequals \nabla \Bu^\nu \rightharpoonup \nabla \Bu \equalscolon \xi \qquad \text{ in } L^2(\Omega, \R^{\Dim\times \Dim}),
  \end{align*}
  because $(\bu,d) \mapsto \eta(\nabla \Bu)$ is in $H^1(\Omega,\R^\Dim \times \R)'$
  for each $\eta \in L^2(\Omega, \R^{\Dim\times \Dim})'$.
  Now Theorem~3.4 of~\cite{Dacorogna1989} provides
  \begin{equation*}
    \liminf_{\nu \to \infty} J(\Bu^\nu,d^\nu,\xi^\nu) \ge J(\Bu,d, \xi) = J(\Bu, d,\nabla \Bu).
    \qedhere
  \end{equation*}
\end{proof}

As a direct consequence of coercivity and weak lower semicontinuity
we get existence of a minimizer of the increment functional:
\begin{theorem}
  There is a solution to the minimization problem
  \eqref{eq:increment_problem}, i.e., there exists
  a global minimizer
  $(\Bu_{n+1}, d_{n+1}) \in \mathbf{H}^1_{\bu_0}\times H^1_{d_0}$
  of $\Pi^\tau_{n+1}$.
\end{theorem}

\subsection{Finite element discretization}

The increment problem~\eqref{eq:increment_problem} of the previous section is posed on the pair of spaces $\mathbf{H}^1_{\Bu_0}$
for the displacements and $H^1_{d_0}$ for the damage variable.  Let $\mathcal{G}$ be a conforming
finite element grid for $\Omega$. We discretize the function spaces by standard
first-order Lagrangian finite elements.
In order to derive an algebraic form of the discretized increment functional
we make use of the standard scalar nodal basis $\{\theta_i\}_{i=1}^M$
associated to the grid nodes $\{p_1,\dots,p_M\} \equalscolon \mathcal{N} \subset \Omega$.
Identifying the $\R^\Dim$-valued and scalar finite element functions
$\bu$ and $d$ with their coefficient vectors $\bu \in \R^{M,\Dim}$ and $d \in \R^M$, respectively,
we write
\begin{align*}
  \bu_j = \sum_{i=1}^M u_{i,j}\theta_i, \qquad d = \sum_{i=1}^M d_i \theta_i,
\end{align*}
where $u_{i,j} = \bu_j(p_i)$ and $d_i=d(p_i)$.
For the integration we use three kinds of quadrature rules:
Integrals of smooth nonlinear terms over a grid element $e$ are
approximated using a higher-order quadrature rule $\int_{e,h}\,dV$,
while the
integral over the nonsmooth term $I_{[d^n,1]} (d)$ is approximated
using the grid nodes $p_i$ as quadrature points, which is often referred to
as \emph{lumping}.
All polynomial terms are integrated exactly using appropriate quadrature rules.
Using these approximations we obtain
the algebraic increment functional
$\Jcal \colonequals \Pi^{\tau,\mathcal{G}}_{n+1}$
given by
\begin{equation}
  \label{eq:lumped_algebraic_increment_functional}
 \Jcal(\bu,d)
 \colonequals
    \underbrace{
      \int_{\Omega,h}(\psi(\bvarepsilon(\Bu),d))\, dV
      + \int_\Omega g_c\gamma(d,\nabla d)\,dV
      + P_\text{ext}(t_{n+1},\bu)
    }_{\equalscolon \Jcal_0(\bu,d)}
    + \underbrace{
      \sum_{i=1}^M I_{[d_n(p_i),1]}(d_i)
    }_{\equalscolon \varphi(d)}.
\end{equation}
Here the quadrature rule $\int_{\Omega,h} (\cdot)\,dV$ is given by
\begin{align*}
  \int_{\Omega,h}f\, dV \colonequals \sum_{e \in \mathcal{G}}\int_{e,h}f\, dV, \qquad
  \int_{e,h}f \,dV \colonequals \sum_{\alpha=1}^{\alpha_{\max}} f(q_{e,\alpha}) \omega_{e,\alpha}
\end{align*}
with positive weights $\omega_{e,\alpha}$ on each element $e$.
Notice that we do not need quadrature weights in the last term of $\Jcal$,
because the indicator function only takes values in $\{0,\infty\}$.

To elucidate the algebraic structure of $\Jcal$
we introduce the linear operator
$\Lcal : (\R^{M,\Dim} \times \R^M) \to ((\S^\Dim \times \R)^{\alpha_{\max}})^{\mathcal{G}}$
with
\begin{align*}
  \Lcal(\bu,d)_{e,\alpha} \colonequals ((\bve(\bu))(q_{e,\alpha}), d(q_{e,\alpha}))
  \qquad \alpha=1,\dots,\alpha_{\max}, \quad e \in \mathcal{G}.
\end{align*}
Then the first part $\Jcal_0$ of the functional can be written as
\begin{equation*}
  \Jcal_0(\bu,d) =
    \underbrace{
      \sum_{e \in \mathcal{G}}
      \sum_{\alpha=1}^{\alpha_{\max}}
      \psi(\Lcal(\bu,d)_{e,\alpha}) \omega_{e,\alpha}
    }_{\equalscolon A(\bu,d)}
    +
    \underbrace{
      \int_\Omega g_c\gamma(d,\nabla d)\,dV
      + P_\text{ext}(t_{n+1},\bu)
    }_{\equalscolon B(\bu,d)}.
\end{equation*}
Note that for the price of a more complex index notation, the linear operator
$\Lcal$ can also be written as a sparse matrix with suitable blocking structure.
In this case $\Lcal(\cdot,\cdot)_{e,\alpha} : (\R^{M,\Dim} \times \R^M) \to (\S^\Dim \times \R)$
corresponds to the $(e,\alpha)$-th sparse row of this matrix.

As an approximation of the boundary conditions from
$\mathbf{H}^1_{\bu_0}\times H^1_{d_0}$
we will consider~$\Jcal$ on the affine subspace
$H^\text{alg} = H^\text{alg}_{\Bu_0} \times H^\text{alg}_{d_0}$
where
\begin{align*}
  H^\text{alg}_{\Bu_0} & \colonequals \big\{\Bu \in \R^{M,\Dim} \, |\,
      \Bu(p) = \Bu_0(p) \quad\forall p \in \mathcal{N} \cap \Gamma_{D,\Bu} \big\},\\
  H^\text{alg}_{d_0}   & \colonequals \big\{d \in \R^M \, |\,
      d(p) = d_0(p) \quad\forall p \in \mathcal{N} \cap \Gamma_{D,d}\big\}.
\end{align*}
The associated homogeneous subspace of $H^\text{alg}$ is denoted by $H^\text{alg}_0$.
In the following we make the assumption that $\mathcal{N} \cap \Gamma_{D,\Bu}$
contains the vertices of at least one boundary grid face, which ensures
a discrete Korn inequality
such that $\|\bve(\Bu)\|_0 \geq C\|\Bu\|_1$ holds
for all $\Bu \in H^\text{alg}_{\Bu_0}$.
Furthermore we introduce the discrete feasible set
\begin{align*}
  \mathcal{K}^\text{alg}
  \colonequals
    H^\text{alg}_{\Bu_0} \times \left(H^\text{alg}_{d_0} \cap \mathcal{K}^\text{alg}_d\right), \qquad
  \mathcal{K}^\text{alg}_d
  \colonequals
    \prod_{i=1}^M [d_n(p_i),1]
\end{align*}
that additionally incorporates the pointwise irreversibility constraints.

\subsection{Properties of the discrete incremental potential}

The convergence property of the TNNMG algorithm heavily rely on the algebraic
structure of the problem. Hence we now collect the essential structural
properties of the algebraic increment functional $\Jcal$.
While stronger properties hold true for some splittings
of~$\psi$, we only note the necessary properties
shared by all of the proposed splittings.
In order to preserve the significant properties
in the presence of numerical quadrature, we assume that the
quadrature rule $\int_{e,h}f \,dV$
can at least integrate
the isotropic energy $f=|\bve(\bu)|_F^2$ exactly for
any finite element function $\bu$.

\begin{lemma}
  \label{lemma:fully_discrete_properties}
  The functional $\Jcal_0 = A+B$ has the following properties:
  \begin{enumerate}
    \item
      $\Jcal_0(\cdot,d) \in C^{1,1}$ and $\Jcal_0(\bu,\cdot) \in C^2$
      for any $d \in \mathcal{K}^\text{alg}_d$ and $\bu \in \R^{M,\Dim}$.
    \item
      The gradient $\nabla \Jcal_0(\cdot,d)$ is semismooth
      for any $d \in \mathcal{K}^\text{alg}_d$.
    \item
      The gradient $\nabla \Jcal_0(\cdot,d)$ is globally Lipschitz continuous uniformly in $d$.
    \item
      $\Jcal_0(\cdot,d)$ is strongly convex uniformly in $d$ on $H^\text{alg}_{\Bu_0}$.
    \item
      $\Jcal_0(\bu,\cdot)$ is convex on $\mathcal{K}^\text{alg}_d$
      for any $\bu \in \R^{M,\Dim}$.
  \end{enumerate}
\end{lemma}
\begin{proof}
  The smoothness properties,
  uniform global Lipschitz continuity,
  and convexity
  follow from the corresponding properties of $\psi$ and $\gamma$,
  and from linearity of~$\Lcal$.

  To see uniform strong convexity, we note that
  uniform strong convexity of $\psi(\cdot,d)$
  implies that there is some $\eta>0$ such that
  $\phi(\bve,d) \colonequals \psi(\bve,d) - \frac{\eta}{2}|\bve|_F^2$
  is convex.
  Using the exactness assumption on the quadrature rule
  we get
  \begin{align*}
    \sum_{e \in \mathcal{G}}
      \sum_{\alpha=1}^{\alpha_{\max}}
      \phi(\Lcal(\bu,d)_{e,\alpha}) \omega_{e,\alpha}
    = A(\bu,d) - \frac{\eta}{2}\|\bve(\Bu)\|_0^2.
  \end{align*}
  Since this is a weighted sum of convex
  functions $\phi(\cdot, d(q_{e,\alpha}))$
  with positive weights $\omega_{e,\alpha}$,
  it is itself convex with respect to $\Bu$.
  Thus $A(\bu,d)$ is the sum of a convex function and the
  function $\|\bve(\Bu)\|_0^2$. Since the latter
  is strongly convex on $H^\text{alg}_{\Bu_0}$ independently of $d$,
  the same applies to $A(\bu,d)$ and $(A+B)(\cdot,d)$.

  Finally we note that convexity of $g$ and $\gamma$
  imply convexity of $(A+B)(\bu,\cdot)$.
\end{proof}

The TNNMG algorithm is based on a crucial property called
block-separability, which states that the nonsmooth part
of the objective functional
can be written as a sum, such that the sets of
independent variables of the addend functionals are disjoint.
We note that $\Jcal = \Jcal_0+\varphi$ is of the desired form,
with a smooth part $\Jcal_0 = A+B$ and a block-separable nonsmooth part
\begin{align}
  \label{eq:nonsmooth_decomposition}
  \varphi(d) \colonequals \sum_{i=1}^M \varphi_i(d_i), \qquad
  \varphi_i(\xi) \colonequals I_{[d_n(p_i),1]}(\xi),
\end{align}
which can also be written as the indicator functional
$\varphi(d) = I_{\mathcal{K}^\text{alg}_d}(d)$
of the feasible set $\mathcal{K}^\text{alg}_d$
of the $n+1$-th time step.

Due to the nonsmoothness of $\varphi$,
the smoothness properties of $\Jcal_0$ do obviously not
carry over to the full functional $\Jcal$.
Furthermore $\Jcal$ is in general not
convex as a whole.
However we still have the following:

\begin{lemma}
  \label{lemma:fully_discrete_coercivity}
  The functional $\Jcal$
  is proper, lower semicontinuous, and coercive
  on~$H^\text{alg}$.
  Furthermore it is convex in $\bu$ and convex in $d$.
\end{lemma}
\begin{proof}
  Being the indicator function of the closed, nonempty, convex
  set $\mathcal{K}^\text{alg}_d$ it is clear that the separable nonsmooth
  functional $\varphi$ is convex, proper, and lower semicontinuous.
  Combining this with smoothness of $\Jcal_0$ we get that
  $\Jcal$ is proper and lower semicontinuous.
  Similarly, convexity in $\bu$ and $d$ follows from the corresponding
  properties of $\Jcal_0$ and $\varphi$.

  Using the uniform coercivity \ref{item:density_coercive}
  of $\psi(\cdot,d)$, $w(1)>0$, and $w(d)\geq 0$
  (as in the proof of Lemma~\ref{lemma:time_discrete_coercivity})
  and the exactness assumption on the quadrature rule
  (as in the proof of Lemma~\ref{lemma:fully_discrete_properties})
  we get
  \begin{align*}
    A(\Bu,d) + B(\Bu,d) - P_\text{ext}(t_{n+1},\bu)
    \geq C \int_\Omega |\bve(\Bu)|_F^2 + |\nabla d|^2 \,dV
  \end{align*}
  for some constant $C>0$.
  Now we can proceed as in the proof of
  Lemma~\ref{lemma:time_discrete_coercivity}
  to show coercivity of $\Jcal$.
\end{proof}

\section{Truncated Nonsmooth Newton Multigrid for brittle fracture}
\label{sec:tnnmg_method}

The Truncated Nonsmooth Newton Multigrid method (TNNMG) is designed to solve nonsmooth block-separable minimization problems
on Euclidean spaces.
In a nutshell, one step of the TNNMG method consists of a nonlinear
Gauß--Seidel-type smoother and a subsequent inexact
Newton-type correction in a constrained subspace.
The nonlinear smoother
computes local corrections by subsequent (possibly inexact)
solving of reduced minimization problems
in small subspaces.
As the nonlinear smoother is responsible
for ensuring convergence, while the Newton corrections
accelerate the convergence, the ingredients of the
nonlinear smoother have to be selected carefully.

It is a well known result \cite{Glowinski1984_NumNonlinVar}
that nonsmooth Gauß--Seidel-type
methods can easily get stuck if the subspace decomposition
used to construct localized minimization problems is not aligned
with the decomposition induced by the block-separable
nonsmooth term. In our case, the nonsmooth term $\varphi$
is separable with respect to the decomposition of
unknowns induced by the grid vertices.
An additional requirement is that the local minimization problems
must be uniquely solvable, which is typically
ensured by choosing the decomposition such that the
local problems are strictly convex.

In view of these requirements we first decompose the
space according to the grid vertices and then with respect
to the local $\Bu$- and $d$-degrees of freedom leading to
a decomposition
\begin{align}
  \label{eq:subspace_splitting}
  \R^{M,\Dim} \times \R^M
  = (\R^{\Dim} \times \R)^M 
  = \sum_{j=1}^M \Bigl(V_{j,\Bu} + V_{j,d}\Bigr).
\end{align}
Here the $\Dim$-dimensional subspace $V_{j,\bu}$ represents
the displacement components at the $j$-th grid vertex,
while the one-dimensional subspace $V_{j,d}$ represents
the $d$-component at this vertex.
All other components are set to zero in these spaces
such that the subspace decomposition can be written as
a direct sum.
For simplicity we use a plain enumeration of these
subspaces in alternating order
\begin{align}
  \label{eq:subspace_enumeration}
  V_{2j-1} &= V_{j,\Bu}, &
  V_{2j} &= V_{j,d}, &
  &j=1,\dots,M.
\end{align}
Notice that with this splitting none of the
nonsmooth terms $\varphi_i$ in \eqref{eq:nonsmooth_decomposition}
couples across different subspaces.
Furthermore, by Lemma~\ref{lemma:fully_discrete_properties}
the restriction of $\Jcal$ to any
affine subspace $(\Bu,d) + V_i$, $i=1,\dots,2M$ is convex.

We will now introduce the TNNMG method.
For simplicity we first assume that
$\Gamma_{D,\Bu}$ and $\Gamma_{D,d}$ are empty
and that $\Jcal_0$ is $C^2$.
Let $\nu \in \mathbb{N}_0$ denote the iteration number.
Given a previous iterate $\BU^\nu = (\Bu,d)^\nu \in \R^{M,\Dim} \times \R^M$,
one iteration of the TNNMG method consists of the following four steps:
\begin{enumerate}
    \item\label{enum:presmoothing}
        \textbf{Nonlinear presmoothing}
        \begin{enumerate}
            \item
              Set $\BW^0 = \BU^\nu$
            \item
                For $i=1,\dots,2M$ compute $\BW^i \in \BW^{i-1} + V_i$ as
                        \begin{equation}
                        \label{eq:general_local_problem}
                            \BW^i \approx \argmin_{\BW \in \BW^{i-1} + V_i} \Jcal(\BW)
                        \end{equation}
            \item
              Set $\BU^{\nu+\frac12} = \BW^{2M}$
        \end{enumerate}
    \item \label{enum:linear_correction}
        \textbf{Inexact linear correction}
        \begin{enumerate}
            \item
                Determine the maximal subspace $W_\nu \subset \R^{M,\Dim} \times \R^M$
                such that the restriction
                $\Jcal|_{W_\nu}$ is $C^2$ at $\BU^{\nu + \frac{1}{2}}$
            \item
                Compute $c^\nu \in W_\nu$ as an inexact Newton step on $W_\nu$
                \begin{align}
                \label{eq:inexact_linear_correction}
                    c^\nu \approx
                        -\big(\Jcal''(\BU^{\nu+\frac12})|_{W_\nu \times W_\nu}\big)^{-1}
                        \big(\Jcal'(\BU^{\nu+\frac12})|_{W_\nu} \big)
                \end{align}
        \end{enumerate}
    \item \label{enum:projection}
        \textbf{Projection}\\
                Compute the Euclidean projection
                $c_\text{pr}^\nu = P_{\operatorname{dom}\Jcal - \BU^{\nu+1/2}}(c^\nu)$,
                i.e., choose $c_\text{pr}^\nu$ such that $\BU^{\nu + \frac{1}{2}} + c_\text{pr}^\nu$ is closest
                to $\BU^{\nu + \frac{1}{2}} + c^\nu$ in $\operatorname{dom} \Jcal$
    \item  \label{enum:line_search}
        \textbf{Damped update}
        \begin{enumerate}
            \item
                Compute a $\rho_\nu \in [0,\infty)$ such that
                $\Jcal(\BU^{\nu+\frac12} + \rho_\nu c_\text{pr}^\nu) \leq \Jcal(\BU^{\nu+\frac12})$
              \item Set $(\Bu,d)^{\nu+1} = \BU^{\nu+1}= \BU^{\nu+\frac12} + \rho_\nu c_\text{pr}^\nu$
        \end{enumerate}
\end{enumerate}

The algorithm is easily generalized to non-trivial Dirichlet
boundary conditions by leaving out all subspaces associated to
Dirichlet vertices during the nonlinear smoothing,
and by additionally requiring $W_\nu \subset H^\text{alg}_0$
for the linear correction subspace.
Then, if the initial iterate satisfies the boundary conditions,
i.e., if $(\Bu_d)^0 \in H^\text{alg}$, the method will only
iterate within this affine subspace, which preserves the Dirichlet boundary conditions for
all iterates.

The canonical choice for the linear correction step~\eqref{eq:inexact_linear_correction}
is a single linear multigrid step, which explains
why the overall method is classified as a multigrid method. If a grid hierarchy is available, then
a geometric multigrid method is preferable. Otherwise, a suitable constructed
algebraic multigrid step for small-strain elasticity problems will work just
as well.
In the case of the phase-field brittle-fracture increment functional considered here,
the nonlinear presmoothing step takes a large part of the run-time of a single iteration.
The convergence speed can therefore be improved considerably by doing a small fixed number
(larger than~1) of multigrid steps, without appreciably increasing the time per iteration.

Section~\ref{sec:convergence} will discuss convergence
of the method based on an abstract convergence theory.
The abstract theory will be used as a guideline
for the discussion of nonlinear smoothers in Sections~\ref{sec:smoothers}.
Finally \ref{sec:multigrid_corrections} will
discuss the linear correction in more detail.

\subsection{Convergence results}
\label{sec:convergence}

The TNNMG method was originally introduced for convex
problems where global convergence to global minimizers
can be shown~\cite{graeser_sack_sander:2009,GraeserKornhuber2009b,GraeserSander2014_preprint}.
These classical results cannot be applied here, due to
the non-convexity of~$\Jcal$.
As a generalization of previous results, \cite{GraeserSander2019}
introduced an abstract convergence theory that also
covers non-convex problems.
In the following we will summarize some results from
this work.
These will later be used as a guideline
for specifying how to solve the local subproblems
\eqref{eq:general_local_problem} and the
linear correction problem \eqref{eq:inexact_linear_correction}.
In order to simplify the presentation some of the terminology and notation
used in \cite{GraeserSander2019} is avoided in favor
of a more specific notation adjusted to the algorithm
as introduced above.

\begin{theorem}
  \label{thm:convergence}
  Let $\Jcal: \R^L \to \R \cup \{\infty\}$
  be coercive, proper, lower semicontinuous,
  and continuous on its domain, and assume that
  $\Jcal(\BV+ (\cdot))$ has a unique global minimizer
  in $V_i$ for all $i$ and each $\BV \in \operatorname{dom} \Jcal$.
  Assume that the inexact local corrections $\BW^i$ are given by
  $\BW^i = \Mcal_i(\BW^{i-1})$
  for local correction operators
  \begin{align*}
    \Mcal_i : \operatorname{dom}\Jcal \to \operatorname{dom} \Jcal,
    \qquad \Mcal_i - \operatorname{Id}: \operatorname{dom}\Jcal \to V_i
  \end{align*}
  having the properties:
  \begin{enumerate}
    \item
      Monotonicity: $\Jcal(\Mcal_i(\BV)) \leq \Jcal(\BV)$
      for all $\BV \in \operatorname{dom}\Jcal$.
    \item
      Continuity: $\Jcal \circ \Mcal_i$ is continuous.
    \item\label{enum:stability}
      Stability: $\Jcal(\Mcal_i(\BV)) < \Jcal(\BV)$ if $\Jcal(\BV)$ is not minimal in $\BV + V_i$.
  \end{enumerate}
  Furthermore assume that the initial iterate is feasible, i.e., $\BU^0 \in \operatorname{dom} \Jcal$,
  and that the linear correction is monotone, i.e., $\Jcal(\BU^{\nu+1}) \leq \Jcal(\BU^{\nu+\frac12})$.
  Then any accumulation point $\BU$ of $(\BU^\nu)$ is stationary in the sense that
  \begin{align}\label{eq:stationary_point}
    \Jcal(\BU) \leq \Jcal(\BU + \BV) \qquad \forall \BV \in V_i, \quad \forall i.
  \end{align}
\end{theorem}
\begin{proof}
  This is Theorem~4.1 in~\cite{GraeserSander2019}.
\end{proof}

\begin{remark}
  Note that the statement requires both global lower semicontinuity of~$\Jcal$
  and continuity on its domain, because neither property implies the other.
  In the proof given in~\cite{GraeserSander2019},
  lower semicontinuity is needed to show that accumulation points
  of the iteration have finite energy and are thus feasible,
  while continuity of~$\Jcal$ on its domain is needed to show
  that they are stationary.
\end{remark}

Now we discuss the application of this theorem to the phase-field brittle-fracture problem.
First we interpret the stationarity result.

\begin{proposition}
  Let $\Jcal$ be given by \eqref{eq:lumped_algebraic_increment_functional}
  and the subspaces $V_i$ by \eqref{eq:subspace_splitting} and \eqref{eq:subspace_enumeration}.
  Then any stationary point $\BU$ in the sense of \eqref{eq:stationary_point}
  is first-order optimal in the sense of
  \begin{align*}
    \inner{\nabla \Jcal_0(\BU)}{\BW-\BU} \leq 0 \qquad \forall \BW \in \operatorname{dom}{\Jcal}.
  \end{align*}
\end{proposition}
\begin{proof}
  The stationarity \eqref{eq:stationary_point} implies the
  variational inequalities
  \begin{align*}
    \inner{\nabla \Jcal_0(\BU)}{\BW_i-\BU} \leq 0 \qquad \forall \BW_i \in \operatorname{dom}{\Jcal} \cap (\BU+V_i)
  \end{align*}
  for each subspace $V_i$.
  Now let $\BW \in \operatorname{dom}\Jcal$.
  Since the splitting \eqref{eq:subspace_splitting}
  is direct one can splitting $\BW$ uniquely into
  \begin{align*}
    \BW = \BU + \sum_{i=1}^{2M} \BV_i, \qquad \BV_i \in V_i.
  \end{align*}
  Using the product structure of $\operatorname{dom} \Jcal$
  we find that $\BW_i \colonequals \BU + \BV_i \in \operatorname{dom}{\Jcal} \cap (\BU+V_i)$.
  Summing up the variational inequalities for those $\BW_i$
  we obtain the assertion.
\end{proof}

Next we investigate the assumptions of the theorem.
First we note that $\Jcal$ as given in \eqref{eq:lumped_algebraic_increment_functional}
is coercive, proper, and lower semicontinuous by Lemma~\ref{lemma:fully_discrete_coercivity}.
Furthermore, $\Jcal_0$ is continuous and the indicator function $\varphi$
is continuous on its domain. Hence the latter is also true for $\Jcal = \Jcal_0+\varphi$.
Subspaces $V_i$ with odd index $i$ only vary in $\Bu(p_i)$
such that existence of a unique minimizer of $\Jcal(\BW+ (\cdot))|_{V_i}$
follows from the strong convexity of $\Jcal(\cdot,d)$ shown in
Lemma~\ref{lemma:fully_discrete_properties}.
For even $i$ these subspaces are associated to nodal damage degrees
of freedom $d(p_i)$.
Although $\Jcal(\Bu,\cdot)$ is in general only convex,
but not strictly convex, the restriction
$\Jcal(\BW+ (\cdot))|_{V_i}$ to a single node
is a strictly convex quadratic functional,
which again implies existence of a unique minimizer.
Finally the monotonicity $\Jcal(\BU^{\nu+1}) \leq \Jcal(\BU^{\nu+\frac12})$
of the linear correction is a direct consequence of the damped update.

It remains to identify proper local correction operators
$\Mcal_i$ satisfying the above assumptions.
As a first result we show that solving the local minimization
problems~\eqref{eq:general_local_problem} exactly
leads to a convergent algorithm in the sense given above.

\begin{lemma}
  \label{lem:exact_local_solve}
  Let $\Jcal$ be given by \eqref{eq:lumped_algebraic_increment_functional}
  and the subspaces $V_i$ by \eqref{eq:subspace_splitting} and \eqref{eq:subspace_enumeration}.
  Then the exact local solution operators
  \begin{align*}
    \Mcal_i(\BW) \colonequals \argmin_{\BV \in \BW + V_i} \Jcal(\BV)
  \end{align*}
  satisfy the assumptions of Theorem~\ref{thm:convergence}.
\end{lemma}
\begin{proof}
  This is Lemma~5.1 in~\cite{GraeserSander2019}.
\end{proof}

Depending on the damaged energy density $\psi$, solving
the restricted problems exactly may not be practical.
As a remedy, it is also shown in \cite{GraeserSander2019}
that inexact minimization is sufficient as long
as it guarantees sufficient decrease of the energy.
In fact we do not need $\BW^i = \Mcal_i(\BW^{i-1})$
exactly but may relax this to $\Jcal(\BW^i) \leq \Jcal(\Mcal_i(\BW^{i-1}))$
for a suitable continuous $\Mcal_i$.
However, sufficient decrease is in general hard to check rigorously.
In the following we cite one inexact variant from~\cite{GraeserSander2019}
where sufficient descent
is guaranteed a priori.

\begin{lemma}
  \label{lem:model_local_solve}
  Let $\Jcal$ be given by \eqref{eq:lumped_algebraic_increment_functional}
  and the subspaces $V_i$ by \eqref{eq:subspace_splitting} and \eqref{eq:subspace_enumeration}.
  For each subspace $V_i$ let $C_i$ be a symmetric positive definite matrix
  that satisfies
  \begin{align*}
    \inner{\nabla \Jcal_0(\BW+\BV) - \nabla \Jcal_0(\BW)}{\BV} \leq \inner{C_i\BV}{\BV} \qquad \forall \BW \in \operatorname{dom}\Jcal, \BV \in V_i.
  \end{align*} 
  Then the correction operators
  \begin{multline*}
    \Mcal_i(\BW)
       \colonequals \argmin_{\BV \in (\BW + V_i) \cap \operatorname{dom}\Jcal}
        \Jcal(\BW) + \inner{\nabla \Jcal_0(\BW)}{\BV-\BW} \\ + \frac{1}{2}\inner{C_i(\BV-\BW)}{\BV-\BW}
  \end{multline*}
  satisfy the assumptions of Theorem~\ref{thm:convergence}.
\end{lemma}
\begin{proof}
  This is Lemma~5.8 in~\cite{GraeserSander2019}.
\end{proof}

\subsection{Smoothers for brittle fracture problems}
\label{sec:smoothers}

The smoother of the TNNMG method performs a sequence of (inexact) minimization problems in
low-dimensional subspaces $V_i$. Different approaches are possible here, implementing
different compromises between convergence speed, wall-time per iteration, and ease of programming.
As there are two types of degrees of freedom, two types of solvers are needed as well.

\subsubsection{Subspaces of displacement degrees of freedom}

We first consider the subspaces
$V_i = V_{(i+1)/2,\bu}$
for odd $i$ spanned by the $\Dim$ displacement degrees of freedom
at the vertex $p_{(i+1)/2}$.
Noting that elements of this subspace only vary
in the displacement component,
the minimization problem~\eqref{eq:general_local_problem}
is equivalent to
\begin{equation}
\label{eq:displacement_local_problem}
  \argmin_{(\bv,0) \in V_i} L_i(\bv).
\end{equation}
Here, the restricted functional
$L_i(\bv) \colonequals \Jcal(\BW^{i-1} + (\bv,0))$
takes the form
\begin{align}
  \label{eq:local_displacement_increment_energy}
  L_i(\bv) =
      \int_{\Omega,h}
        (g(d)+k)  \psi^+_{0}(\bve(\bu+\bv)) + (1+k)\psi^-_0(\bve(\bu+\bv)) \,dV
      + P_\text{ext}(t_{n+1},\bv) + \operatorname{const},
\end{align}
where we have used $(\bu,d)=\BW^{i-1}$.
The precise nature of $L_i$ depends on the type of energy splitting used by the model.
If the isotropic splittings~\eqref{eq:symmetric_damaged_energy} or~\eqref{eq:isotropic_volumetric_split}
are used, \eqref{eq:local_displacement_increment_energy}
is a strictly convex quadratic functional on a vector space, and can be minimized
exactly by solving an $\Dim \times \Dim$ system of linear equations.

For the anisotropic splittings~\eqref{eq:anisotropic_volumetric_split}
and~\eqref{eq:spectral_energy_split}, the functional is still strictly
convex and once continuously differentiable. The classical Hesse matrix,
however, is not guaranteed to exist.
However, by Lemma~\ref{lemma:fully_discrete_properties},
the increment functional is semismooth.  This suggests various natural choices
for local solvers, such as steepest-descent methods or
nonsmooth Newton methods~\cite{Ulbrich2002}.
When these are used to solve the local problems~\eqref{eq:general_local_problem} exactly,
global convergence of the overall TNNMG method follows from
Theorem~\ref{thm:convergence} and Lemma~\ref{lem:exact_local_solve} above.

However, as mentioned in the previous section, Theorem~\ref{thm:convergence}
is more general, and also shows convergence for certain types
of inexact local solvers (such as the one in Lemma~\ref{lem:model_local_solve}).
Such a setup can make iterations much faster, while keeping
the corresponding deterioration of the convergence rate within acceptable limits.
Possible approaches are:
\begin{itemize}
  \item One Newton-type step
  where the $\psi''$-term in the Hessian $\Jcal_0''$ of the differentiable part $\Jcal_0$ of $\Jcal$ is
  replaced by the
    quadratic upper bound $(1+k) \psi''_0$,
    that is, the undamaged St.\,Venant--Kirchhoff energy density (scaled with $(1+k)$), which is
    strictly convex and quadratic.  By the construction of the splittings in
    Section~\ref{sec:elastic_energy}, this term bounds the degraded energy density~\eqref{eq:damaged_energy_density}
    for any admissible value of $d$. We call this approach a preconditioned smoother.
  \item One (or another fixed number of) semismooth Newton steps.
  \item One gradient step with exact line search.
\end{itemize}
For the first variant, global convergence of the TNNMG solver
follows from Lemma~\ref{lem:model_local_solve}.
For the other two, the problem of showing convergence is open.
Section~\ref{sec:numerical_examples} will show
how the first two choices perform in practice.

\subsubsection{Subspaces spanned by damage degrees of freedom}

For subspaces
$V_i = V_{i/2,d}$
with even $i$, i.e., subspaces spanned by the damage degrees of freedom
at the vertices $p_{i/2}$,
the minimization problem~\eqref{eq:general_local_problem}
is equivalent to
\begin{equation*}
  \argmin_{(0,v) \in V_i} L_i(v),
\end{equation*}
with a restricted functional $L_i(v) \colonequals \Jcal(\BW^{i-1} + (0,v))$.
For all choices of damage functions $g$ described in Section~\ref{sec:phase_field_models}
this is a strictly convex quadratic functional on a closed interval,
whose minimizer can be computed directly
by computing the unconstrained minimizer and projecting onto the admissible interval.
We therefore always assume that these problems are solved exactly.

\subsection{Linear multigrid corrections}
\label{sec:multigrid_corrections}

For the linear correction step~\eqref{eq:inexact_linear_correction}
we need to compute a constrained Newton-type correction
\begin{align}
\label{eq:inexact_linear_correction_repeated}
  c^\nu \approx
  -\big(\Jcal''(\BU^{\nu+\frac12})|_{W_\nu \times W_\nu}\big)^{-1}
  \big(\Jcal'(\BU^{\nu+\frac12})|_{W_\nu} \big)
\end{align}
at least inexactly.
This requires to determine the subspace $W_\nu$,
to compute the constrained first- and second-order derivatives
$\Jcal'(\BU^{\nu+\frac12})|_{W_\nu}$ and
$\Jcal''(\BU^{\nu+\frac12})|_{W_\nu \times W_\nu}$
on this subspace, and finally to solve the system
inexactly.

It is easy to see that the largest subspace $W_\nu$
where $\Jcal$ is differentiable
in a neighborhood of $\BU^{\nu+\frac12}=(\Bu^{\nu+\frac12},d^{\nu+\frac12})$
is given by
\begin{align*}
  W_\nu = \Big\{\BU=(\Bu, d)\;\big|\; d_i=0 \text{ if } (d^{\nu+\frac12})_i \notin (d_{n}(p_i),1)\Big\}.
\end{align*}
In this subspace the nonsmooth indicator functional
$\varphi$ is identical to zero such that we only
need to compute first- and second-order derivatives of the
smooth part~$\Jcal_0$, which are then restricted
to
the degrees of freedom that are allowed to be nonzero
in $W_\nu$. This can easily be achieved by setting
rows and columns of $\Jcal'$ and $\Jcal''$ to zero for degrees of freedom not
contained in $W_\nu$.
For all splittings where the degraded
density $\psi$ is not $C^2$, it is at least locally
Lipschitz and semismooth.
In this case a generalized second-order derivative
$\Jcal_0''(\BU^{\nu+\frac12})$ can be used as a replacement
of the classical Hesse matrix making
\eqref{eq:inexact_linear_correction} a semismooth Newton step.
Such a generalized Hessian can be obtained by the following procedure:
The density $\psi$ is piecewise $C^2$, and every point $(\bve,d)$
where it is not $C^2$
is at the boundary of several subdomains on which $\psi$ is $C^2$.
Whenever the second derivative~$\psi''$ needs to be computed at such a point,
one uses instead the second derivative from any of these adjacent subdomains.

To practically compute the first and second derivatives of the degraded elastic energy
with the splitting~\eqref{eq:spectral_energy_split} based on eigenvalues, recall that
the positive and negative parts~$\psi_0^\pm$ of that splitting are spectral functions
\begin{equation*}
 \psi_0^\pm
 =
 \widehat{\psi}_0^\pm \circ \Eigfunc : \mathbb{S}^m \to \R,
\end{equation*}
where $\Eigfunc : \mathbb{S}^m \to \R^m$ is the ordered eigenvalue function,
and the $\widehat{\psi}^\pm_0 : \R^m \to \R$ are invariant under permutations of their arguments.
Such functions are called \emph{spectral}, and they are once or twice differentiable
if and only if the functions $\widehat{\psi}^\pm_0$ are once or twice differentiable, respectively.
The first and second derivatives of $\psi^\pm_0$ can be expressed in closed form
in terms of the derivatives of $\widehat{\psi}^\pm_0$.  This is shown, for example,
in \cite[Lemma~3.1]{lewis_sendov:2001} for first-order derivatives, and
in \cite[Theorem~3.3]{lewis_sendov:2001} for second-order ones. The expressions involve
computing the eigenvector decomposition of the argument of $\psi^\pm_0$.

For the inexact solution~\eqref{eq:inexact_linear_correction_repeated} one step of a classical linear multigrid method
can be used.
Here we only need to take care that the linear smoother
can deal with the non-trivial kernel resulting from constraining
the linearization. For a linear Gau\ss--Seidel smoother this amounts
to omitting corrections for rows with zero diagonal entry.
When using the TNNMG method for problems with a nonquadratic smooth part
such as the phase-field brittle-fracture problem considered here, the smoother
is relatively expensive. One can then improve the overall convergence rate
by doing more than one multigrid iteration, without increasing the time
per iteration.

%compute the first and second derivatives of the
%smooth part of the increment functional $\Pi^\tau_n$. The first derivative is
%\begin{align}
%  \delta\Pi^\tau & = \int_\Omega \bigg\{
% [g(d) \bsigma_0^+ + \bsigma_0^-] : \nabla\delta\Bu
%  + [\frac{g_c}{l}d  +g'(d) \psi_0^+ + \partial_d I_{[d_n,\infty)}]\delta d
%  \nonumber \\
%  &  + g_c l \nabla d \cdot\nabla\delta d
%  \bigg\} \,\text{d}V
%\ .
%\end{align}
%where $g'(d)=-2(1-d)$ is the derivative of the degradation
%function and the relation $\bsigma_0 \colonequals\partial_\bvarepsilon\psi_0$ is
%used.
%Note that only the tensile part $(*)^+$ of the energy $\psi_0$ and the
%stress $\bsigma_0$ are degraded. The linearization then reads
%\begin{align}
%  \Delta\delta\Pi^\tau & = \int_\Omega\bigg\{
%  \nabla\Delta\Bu :[ g(d) \partial_{\bvarepsilon} \bsigma_0^+
%                   + \partial_{\bvarepsilon} \bsigma_0^- ]:
%                         \nabla\delta\Bu
%   + \Delta d [2\psi + \frac{g_c}{l} + \partial^2_{dd}I_{[d_n,\infty)}] \delta d
%  \nonumber \\
%  & + \nabla\Delta d g_c l \nabla\delta d
%  + \Delta d [g'(d)\bsigma_0^+)] : \nabla\delta\Bu
%  + \Delta\nabla\Bu : [g'(d)\bsigma_0^+] \delta d
%   \bigg\} \,\text{d}V
%\ .
%\label{eq:linpot2}
%\end{align}

\section{Numerical examples}
\label{sec:numerical_examples}

% \TNNMGEX entspricht:
% \texttt{parameterSet.localSolver.type = 'exact'}
% \texttt{parameterSet.outerIterations = 1}
% \texttt{parameterSet.displacementIterations = 1}
%
% \TNNMGPRE entspricht:
% \texttt{parameterSet.localSolver.type = 'preconditioned'}

In this last section we demonstrate the speed and robustness of the TNNMG solver with two
numerical examples. For both of them we study four instances from the family of models discussed in
this manuscript: the isotropic and the spectral
splittings (\eqref{eq:symmetric_damaged_energy} and \eqref{eq:spectral_energy_split}, respectively),
combined with the AT-1 and AT-2 crack density functionals ($w(d) = d$ and $w(d) = d^2$, respectively).

In the experiments we will consider two variants
of the TNNMG method with two different nonlinear smoothers,
both based on the splitting~\eqref{eq:subspace_enumeration}
which alternates displacement and damage degrees of freedom.
For the first variant---denoted \TNNMGEX{} in the following---%
the smoother will solve the local displacement problems \eqref{eq:displacement_local_problem}
inexactly in the sense that one damped (generalized) Newton step
is applied.
In the second variant---denoted \TNNMGPRE{}---%
the smoother will solve approximate local displacement problems exactly.
To this end it employs a single Newton-like step
where the $\psi''$-term in the Hessian $\Jcal_0''$ of the differentiable part $\Jcal_0$ of $\Jcal$ is
replaced
by the quadratic upper bound $(1+k) \psi''_0$.

The constrained quadratic local damage problems are solved exactly
by both smoothers.
By Lemma~\ref{lem:model_local_solve}, the \TNNMGPRE{} smoother satisfies the
assumptions of the convergence Theorem~\ref{thm:convergence}, but the \TNNMGEX{} smoother
does so only for the degraded elastic energy density without a splitting
(otherwise its local solution operator is not continuous).
For the linear correction step~\eqref{eq:inexact_linear_correction} we use three standard $V(3,3)$
linear geometric multigrid steps with a block Gauß--Seidel smoother
operating on the canonical $(d+1)\times(d+1)$ blocks.
The coarse linear problems of the multigrid iteration are solved
using the UMFPack sparse direct solver~\cite{davis:2004}.
Damage degrees of freedom are truncated when they are less than $10^{-10}$
away from their lower bound.
The TNNMG iteration is set to run until the relative degraded energy norm of the
correction drops below $10^{-7}$.  The degraded energy norm for the coupled problem combines
the energy norm of the degraded linear elasticity problem with the energy norm
of the AT-2 crack surface energy density.%
\footnote{Because the AT-1 model does not induce a norm.}

We measure iteration numbers, wall-time, and memory consumption.
The TNNMG algorithm and the operator-splitting algorithm used for comparison
are implemented in C++ using the \textsc{Dune} libraries%
\footnote{\url{www.dune-project.org}}
\cite{sander:2020,bastian_et_al:2021}.

\subsection{An operator splitting method for phase-field brittle fracture models}
\label{sec:operator_splitting}

We compare the performance of the TNNMG method to an operator-splitting
method from the literature.  This method is described here in some detail
to allow readers to reproduce the results. We choose an operator-splitting method
because such methods are widely used,
and reported to be robust~\cite{ambati_gerasimov_delorenzis:2015,wu_huang_nguyen:2020}.

Starting from an initial displacement $\bu^0$ and damage field $d^0$
the operator-splitting method alternately repeats the following steps:
\begin{align}
 \label{eq:os_displacement_step}
 \bu^{\nu+1}
 & =
 \argmin_{\bu \in H^\text{alg}_{\bu_0}} \Jcal(\bu, d^\nu) \\
 \label{eq:os_damage_step}
 d^{\nu+1}
 & =
 \argmin_{d \in H^\text{alg}_{d_0}} \Jcal(\bu^{\nu+1}, d).
\end{align}
The iteration is terminated under the same condition as the
TNNMG method above.

If the degraded elastic energy $\psi(\bve,d)$ is of the
type~\eqref{eq:symmetric_damaged_energy}, which does not distinguish between compressive and
tensile strains, then~\eqref{eq:os_displacement_step} is a linear
problem with a symmetric positive definite matrix.
Our implementation solves these problems using the CHOLMOD
direct sparse solver~\cite{chen_davis_hager_rajamanickam:2008}.%
\footnote{It is well-known that sparse direct
solvers do not scale well to larger problems both run-time and memory-wise.
Therefore, an alternative method such as linear multigrid may be a better
choice here. However, direct solvers are widely used in practice,
and we therefore consider a comparison with such a solver worthwhile.}

If the degraded elastic energy incorporates a compressive--tensile split
such as~\eqref{eq:spectral_energy_split}, then $\Jcal$ is not
a quadratic functional in $\bu$.  In this case, as suggested by
\citet[Equation~(66)]{miehe+hofacker+welschinger10}, we perform a single
Newton step at $(\bu^\nu, d^\nu)$, viz.
\begin{equation*}
 \bu^{\nu+1}
 =
 \bu^\nu - \Big( \nabla_{\bu\bu}^2 \Jcal (\bu^\nu, d^\nu) \Big)^{-1} \nabla_\bu \Jcal(\bu^\nu, d^\nu),
\end{equation*}
where $\nabla_\bu \Jcal$ und $\nabla^2_{\bu\bu} \Jcal$ are the first and
second derivatives of $\Jcal$, respectively, with respect to $\bu$ only.
Because of the splitting, the Hesse matrix $\nabla_{\bu\bu}^2 \Jcal$ does not depend
continuously on $\bu$.  At points where an entry of $\nabla_{\bu\bu}^2 \Jcal$
jumps, we simply pick one of the one-sided limits,
which are elements of the generalized Jacobian
of $\nabla_{\bu} \Jcal$ in the sense of Clarke (cf.\ Remark~\ref{rem:generalized_jacobian}).

For the AT-1 and AT-2 models, the damage subproblem~\eqref{eq:os_damage_step}
is a quadratic minimization problem subject to lower bound constraints
\begin{equation*}
 d_i^\nu \le d_i \le 1
 \qquad
 i = 1,\dots, M.
\end{equation*}
The Hesse matrix is
\begin{equation*}
 \Big(\nabla_{dd}^2 \Jcal(\bu^{\nu+1},d) \Big)_{ij}
 =
 \int_\Omega \Big[ g''(d)\theta_i \theta_j \psi_0^+(\bve(\bu^{\nu+1}))
 +
 \underbrace{\frac{g_c}{2 c_w l} \theta_i \theta_j}_{\text{AT-2 only}}
 +
 \frac{g_c l}{2 c_w} \nabla \theta_i \nabla \theta_j \Big]\,dx,
\end{equation*}
with $g''(d) = 2$ for our choice $g(d) = (1-d)^2$.
For the AT-2 model this Hesse matrix is always positive definite. For the AT-1
model, it is only positive definite if the tensile elastic energy $\psi_0^+$
is positive everywhere, and positive semidefinite otherwise.
This did not lead to any problems in the numerical tests done for this article.
If necessary, a small amount of $L^2$-regularization can be added to
increase the robustness.

Following a suggestion by \citet{burke_ortner_sueli:2013}, we use a projected Newton method
to solve the constrained damage problems.  We use the variant proposed
by~\citet{bertsekas:1982}, because it is well described and, for the quadratic problems
considered here, converges locally quadratically~\cite[Proposition~4]{bertsekas:1982}.
The method combines projected gradient steps for degrees of freedom in the vicinity
of the constraints with Newton-type steps for the other degrees of freedom,
and an Armijo-style line search.  Consult the original article~\cite{bertsekas:1982}
of \citeauthor{bertsekas:1982} for a detailed description.  In the notation
of that article, we let $M$ be the identity matrix and $D_k$ ($k$ being the
iteration number) the truncated Hesse matrix. This matrix results from the true
Hesse matrix $\nabla^2_{dd}\Jcal(\bu^{\nu+1},d)$ by replacing the rows and columns
for which the degrees of freedom are close to or at the constraint
by the corresponding rows and columns of the identity matrix.%
\footnote{Note that this way of truncating degrees of freedom differs from
the one employed by the TNNMG method (the construction of the space $W_\nu$)
in Section~\ref{sec:multigrid_corrections}.}
The maximum truncation distance is $\varepsilon = 10^{-5}$.
For the line search parameters we set $\beta = 0.5$ and $\sigma = 0.49$.
We solve the linear subproblems with the CHOLMOD solver.

Note that the projected Newton method is closely related to the TNNMG method.
Indeed, the TNNMG method could also be applied to solve only the damage
problem~\eqref{eq:os_damage_step} within an operator splitting loop.
The main differences are that TNNMG has a presmoother, and that it does not solve
the linear correction problems exactly.

\subsection{Pure tension test of a notched, symmetric specimen}

\begin{figure}
  \begin{center}
    \begin{overpic}[width = 0.75\textwidth]{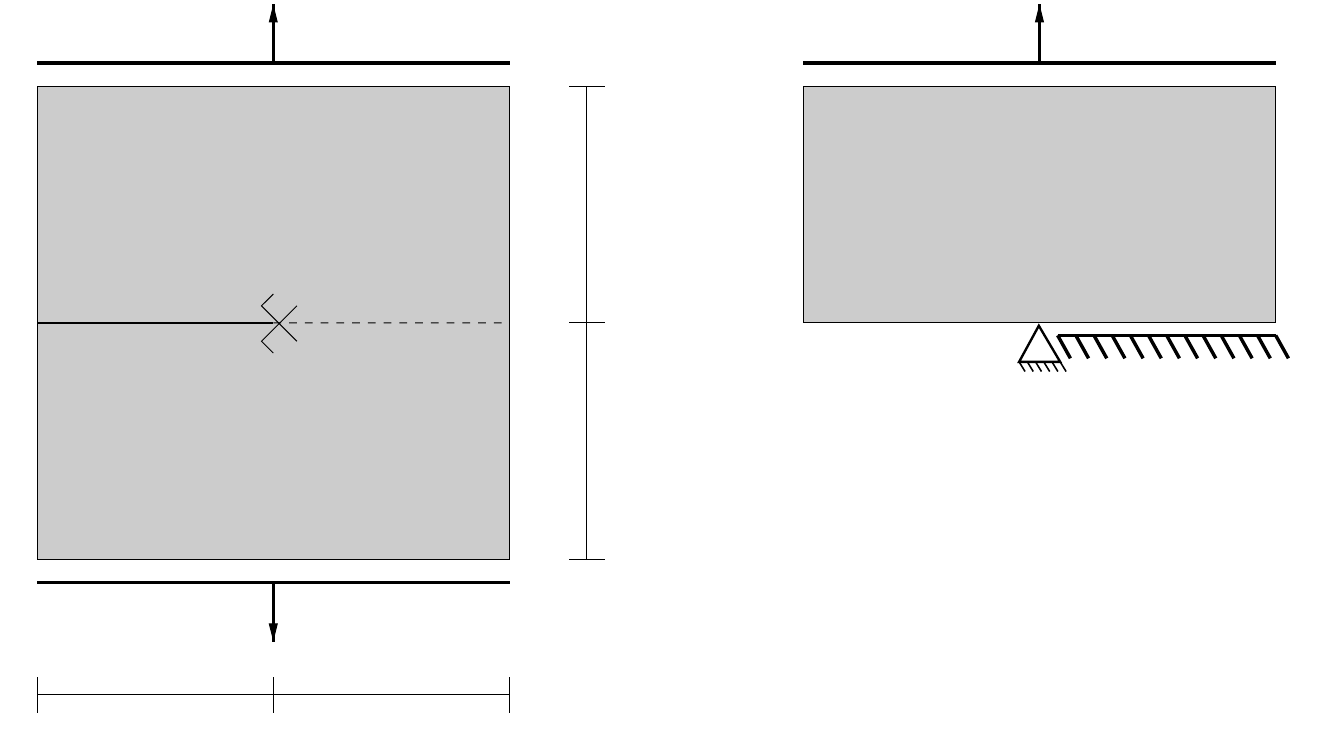}
     \footnotesize
     \put(-5, 31){$(0,0)$}
     \put(53,31){$(0,0)$}
     \put(22, 9){$-\bar\Bu$}
     \put(22,53){$\bar\Bu$}
     \put(81,53){$\bar\Bu$}

     \put( 9, 1){$L/2$}
     \put(27, 1){$L/2$}
     \put(46,23){$L/2$}
     \put(46,40){$L/2$}
    \end{overpic}
    \caption{Notched square with a vertical displacement load (left).
    Exploiting symmetry we only simulate on the upper half of the domain (right).}
    \label{fig:singlenotchbvp}
  \end{center}
\end{figure}

The first numerical example is a two-dimensional, square-shaped notched specimen
of size $L \times L$ under tension.
Due to symmetry we simulate only its upper half. Geometry and boundary conditions
are shown in Figure \ref{fig:singlenotchbvp}. On the top edge a
time-dependent normal displacement $\bar \bu$ is prescribed, while the horizontal displacement
is left free.  The bottom edge of the upper half is clamped vertically for all $x>L/2$,
and fixed vertically and horizontally at the single point $(L/2, 0)$,
which is where the initial crack tip is.
With increasing normal displacement $\bar{\bu}$, the preexisting crack opens,
and the specimen ruptures suddenly, when a limit load is exceeded.
By symmetry, the crack energy is accounted for correctly, even though
only one half of the crack profile appears in the simulation result.

The simulations are performed with
parameters taken from the corresponding experiment in~\cite{miehe+welschinger+hofacker10a},
viz.\ $L=1$\,mm, Lamé parameters $\lambda=121$\,kN/mm$^2$ and $\mu=80$\,kN/mm$^2$,
critical energy release rate $g_c=2.7$\,N/mm, and residual stiffness $k=10^{-5}$.
The phase-field regularization parameter is set to $l=0.03125$\,mm.
We apply the loading in 160 steps, and set the displacement load~$\bar{\bu}$
at step $i$ to
\begin{equation*}
 \bar{\bu}_i = i \cdot 2\cdot 10^{-5}\,\text{mm} \cdot \be_2,
 \qquad
 i = 1,\dots, 160,
\end{equation*}
where $\be_2$ is the canonical basis vector point upwards.
We start the evolution with no displacement and no damage anywhere.

For the spatial discretization we use three different uniform grids
with $256 \times 128$~($h_1$), $512 \times 256$ ($h_2$), and
$1024\times 512$ ($h_3$) quadrilateral elements, respectively.
These were all constructed by uniform refinement of a grid with $32 \times 16$
elements, and hence the grid hierarchy for the multigrid solver consists of~4, 5, and 6 levels,
respectively.
To separate the effect of the discretization parameter $h$
from the modeling parameter $l$, we explicitly study different
element sizes $h$ to illustrate the performance of the algorithm,
instead of choosing a single element size proportional to the crack width.
Relating the grid resolution to the average fracture with $l$,
the average element edge length corresponds to $l/8$, $l/16$, and $l/32$, respectively.

\subsubsection{Isotropic splitting}

\begin{figure}

 \begin{tikzpicture}[scale=0.8]
  % AT-1 model
  \node at (0,10) {\includegraphics[width = 0.32\textwidth]{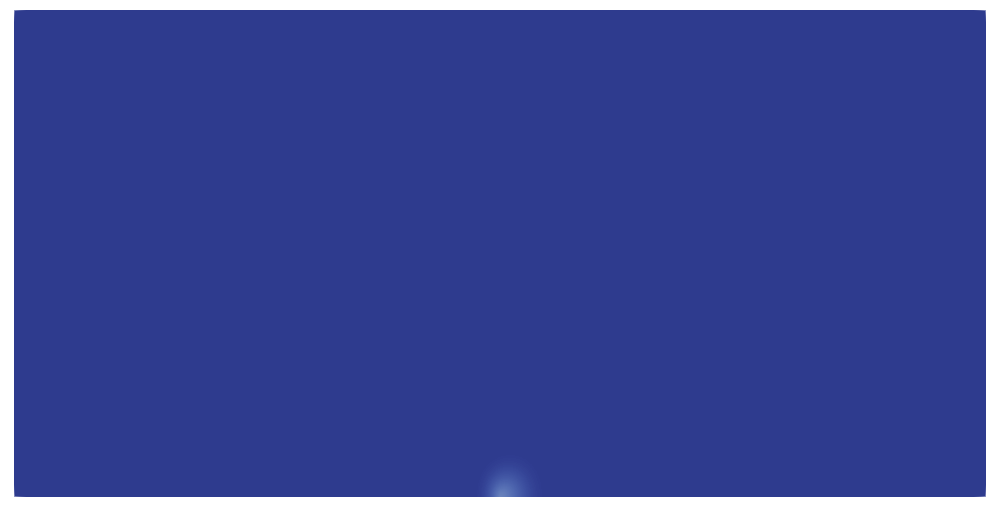}};
  \node at (0, 7) {\includegraphics[width = 0.32\textwidth]{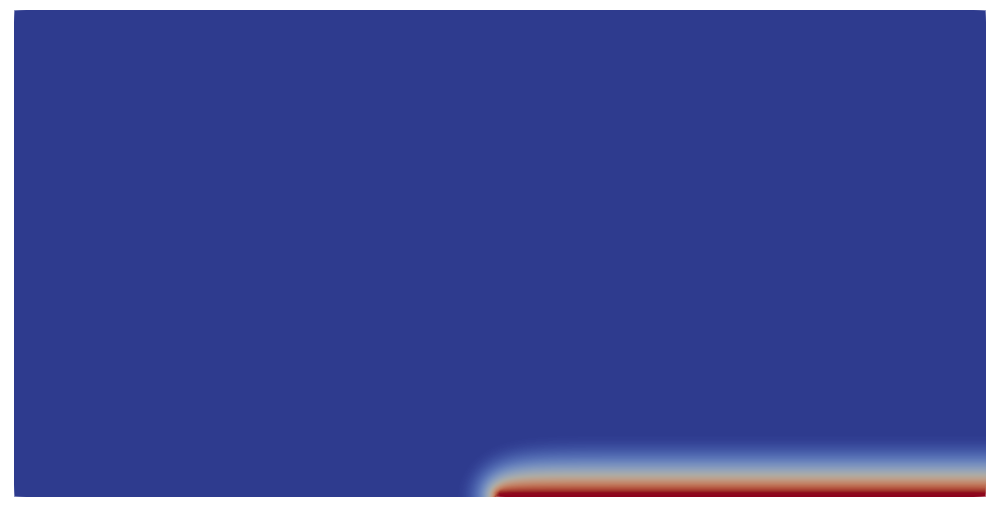}};

  \node at (8,8.5) {\begin{minipage}{0.5\textwidth}
                    \input{2d_reac_at1_iso.pgf}
                   \end{minipage}};

  \node at (4,6) {AT-1};

  % AT-2 model
  \node at (0, 3) {\includegraphics[width = 0.32\textwidth]{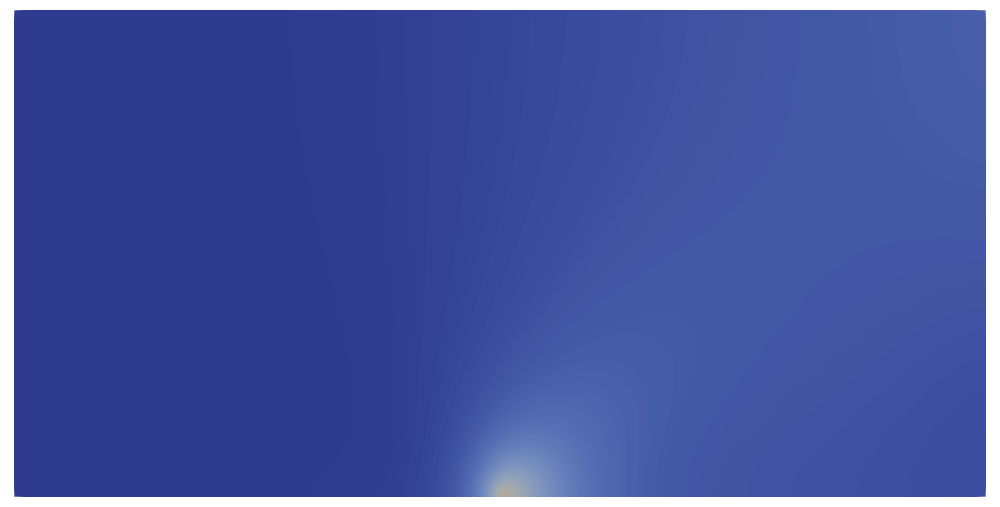}};
  \node at (0, 0) {\includegraphics[width = 0.32\textwidth]{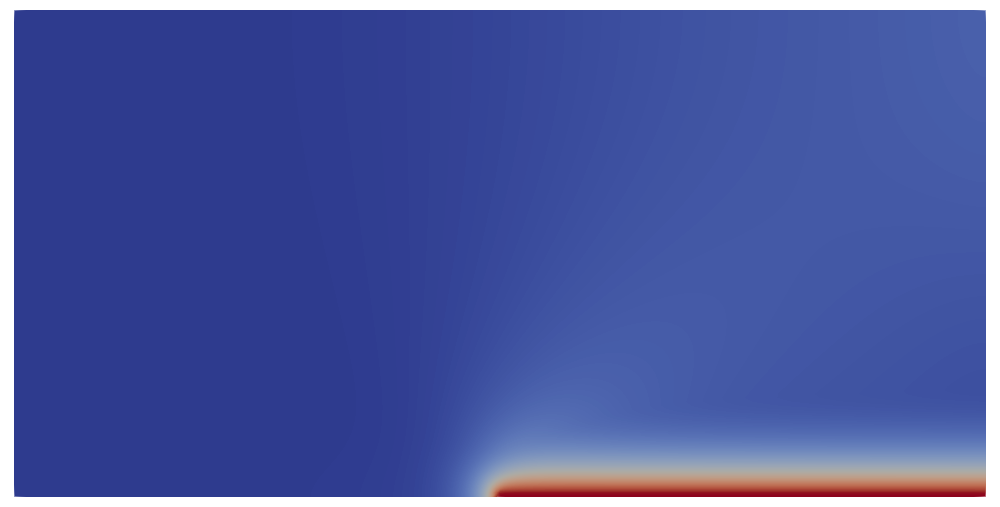}};

  \node at (8, 1.5) {\begin{minipage}{0.5\textwidth}
                    \input{2d_reac_at2_iso.pgf}
                   \end{minipage}};

  \node at (4,-1) {AT-2};
 \end{tikzpicture}

 \caption{2d example, isotropic energy split: Evolution at time steps~145 and~160
  computed with the TNNMG algorithm, and displacement--force curves for the AT-1 (top)
  and AT-2 (bottom) models}
 \label{fig:evolution2diso}
\end{figure}

We first consider the model with the isotropic splitting~\eqref{eq:symmetric_damaged_energy}
of the elastic energy density, where all elastic strains contribute to the degradation
of the material.   Figure~\ref{fig:evolution2diso}
shows the evolution and displacement--force curves, both for the AT-1 and AT-2 functionals.
The force plotted here is the total normal force on the top edge of the specimen.

\begin{figure}
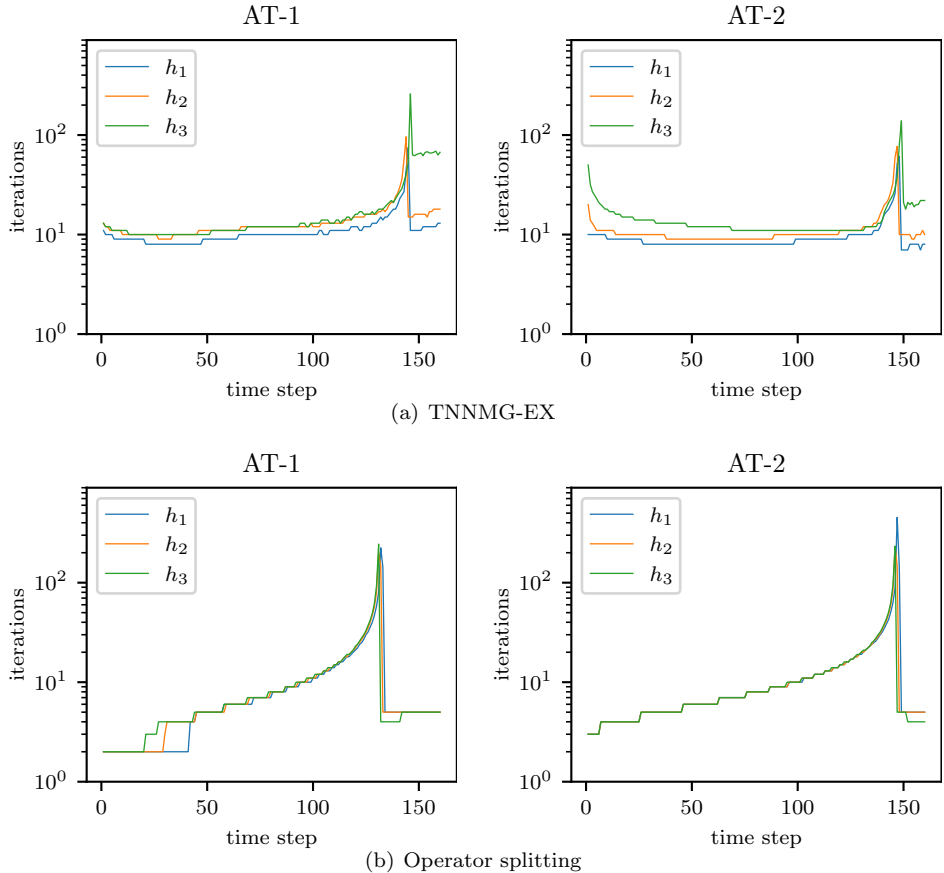


 \subfigure[\TNNMGEX]
 {
  \begin{minipage}{0.5\textwidth}
     \begin{center}
     \input{iter_2d_AT1_iso_ex.pgf}
    \end{center}
  \end{minipage}
  \begin{minipage}{0.5\textwidth}
     \begin{center}
     \input{iter_2d_AT2_iso_ex.pgf}
    \end{center}
  \end{minipage}
 }

 \subfigure[Operator splitting]
 {
  \begin{minipage}{0.5\textwidth}
     \begin{center}
     \input{iter_2d_AT1_iso_feap.pgf}
    \end{center}
  \end{minipage}
  \begin{minipage}{0.5\textwidth}
     \begin{center}
     \input{iter_2d_AT2_iso_feap.pgf}
    \end{center}
  \end{minipage}
 }
 \caption{2d example, isotropic energy split: Iterations per time step, for grid sizes $h_1$, $h_2$, $h_3$}
 \label{fig:iterations_2d_isotropic}
\end{figure}

We first compare iteration numbers.  At each loading step, the increment problem is solved starting
from the solution of the previous time step
until the energy norm of the correction normalized by the energy norm of the previous
iterate drops below~$10^{-7}$.
The upper row of Figure~\ref{fig:iterations_2d_isotropic} shows the number of iterations for the TNNMG solver
with the exact smoother (\TNNMGEX), with logarithmically scaled vertical axis.
As can be seen, the iteration numbers remain essentially bounded independently
from the grid resolution.
The peak shortly before the 150th load step is where the material ruptures.
In this situation, the system becomes highly unstable, which results in a
higher number of iterations.

\begin{figure}
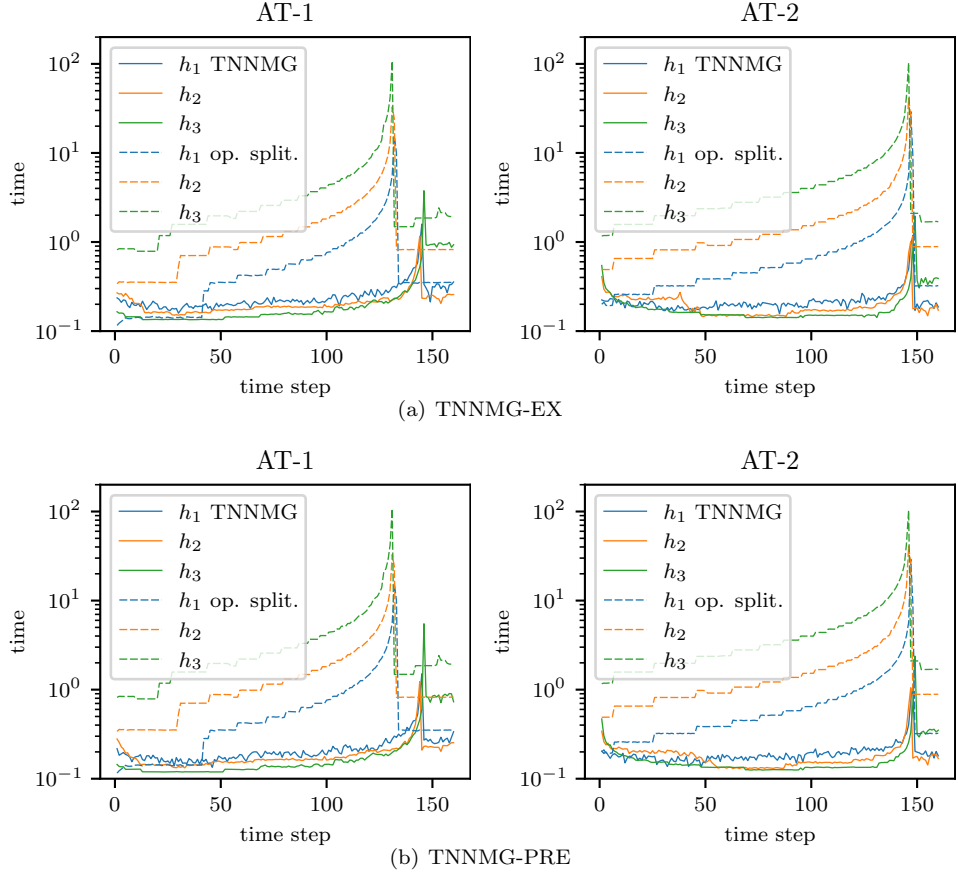


 \subfigure[\TNNMGEX]
 {
  \begin{minipage}{0.5\textwidth}
  \begin{center}
     \input{time_2d_AT1_iso_ex.pgf}
  \end{center}
  \end{minipage}
  \begin{minipage}{0.5\textwidth}
  \begin{center}
     \input{time_2d_AT2_iso_ex.pgf}
  \end{center}
  \end{minipage}
 }

 \subfigure[\TNNMGPRE]
 {
  \begin{minipage}{0.5\textwidth}
  \begin{center}
     \input{time_2d_AT1_iso_pre.pgf}
  \end{center}
  \end{minipage}
  \begin{minipage}{0.5\textwidth}
  \begin{center}
     \input{time_2d_AT2_iso_pre.pgf}
  \end{center}
  \end{minipage}
 }
 \caption{2d example, isotropic energy split: Wall-time per degree of freedom
   per time step, for grid sizes $h_1$, $h_2$, $h_3$}
 \label{fig:walltime_2d_isotropic}
\end{figure}

After the rupture, the iteration numbers do depend on the grid resolution.
This is atypical for a multigrid solver---presumably it is caused by the fact that,
given the particular boundary conditions, the completely ruptured specimen
is essentially an ill-posed problem.

For comparison, the lower row of Figure~\ref{fig:iterations_2d_isotropic} shows the
iteration numbers of the operator splitting method.
One can see that for the first two thirds of the loading history, this method
needs less iterations than the TNNMG algorithm, and that iteration numbers
are independent from the grid resolution.
Recall, however, that TNNMG iterations are much cheaper than operator-splitting
iterations, because they are essentially a low fixed number multigrid iterations,
whereas each operator-splitting iteration involves solving two global
linear systems.  In contrast to the TNNMG method, the number of iterations
increases with increasing load, and shortly before the peak iteration numbers
are higher.

We do not show the iteration numbers of TNNMG with the inexact
smoother (\TNNMGPRE), because they coincide with the results for \TNNMGEX. This is not
surprising:
As the model uses the isotropic splitting of the elastic energy, the local
displacement problems are quadratic, and a single Newton step solves them exactly.
\TNNMGPRE uses a preconditioner for those quadratic problems, which
in this particular situation is very similar to the actual problem. Therefore, preconditioning here
has only a limited impact on the smoothing and thus on the speed
of convergence.
However, since the preconditioner is independent of the current iterate,
the corresponding local matrices can be precomputed.

Wall-time behavior is discussed next.
We plot wall-time per time step for the two multigrid variants
\TNNMGEX and \TNNMGPRE, and for the operator-splitting method (Figure~\ref{fig:walltime_2d_isotropic},
again with logarithmic vertical axes).  The plots show
the time normalized by the number of degrees of freedom.

We see that TNNMG is about 2~to 3~times faster than operator-splitting.
The time per degree of freedom stays roughly constant for
both methods, independent of the grid resolution.
Presumably, the superlinear complexity of the direct solver used in the
operator-splitting method only shows for larger grids.
Table~\ref{tbl:2d_iso_times} shows the accumulated normalized run-times for the
entire load history.  It shows that the speed difference for the $h_3$ grid
accumulates to a factor of about
3.5 for the AT-1 model.  For the AT-2 model the difference is only a factor
of about~2.  This is a bit surprising: The factor is larger for the smaller grids,
but TNNMG, for the AT-2 model, slows down considerably when going from the $h_2$~grid
to the $h_3$~grid. This is caused by higher iteration numbers throughout the load
history, as can be seen in Figure~\ref{fig:iterations_2d_isotropic}.
The reason for this behavior is unclear.

\begin{table}
% TNNMG: time_2d_AT1_iso_ex.dat [5.89357132 9.55973812 8.12692169]
% TNNMG: time_2d_AT1_iso_pre.dat [5.67322153 6.57834235 7.81396231]
% FEAP:  time_2d_AT1_iso_feap.dat [20.47279784 20.5061114  28.84459621]
%
% TNNMG: time_2d_AT2_iso_ex.dat [ 6.73590394  9.9434536  20.77936975]
% TNNMG: time_2d_AT2_iso_pre.dat [ 5.38512     7.18139805 18.28866632]
% FEAP:  time_2d_AT2_iso_feap.dat [29.18803376 29.42741459 40.9208234 ]

\small
 \begin{tabular}{c|ccc|ccc}
 \hline
& \multicolumn{3}{c}{AT-1} & \multicolumn{3}{c}{AT-2} \\
&  \TNNMGEX & \TNNMGPRE & OS & \TNNMGEX & \TNNMGPRE & OS \\
\hline
$h_1$ &  5.89 &  5.67 & 20.47  &       6.74 &  5.39 & 29.19 \\
$h_2$ &  9.56 &  6.58 & 20.51  &       9.94 &  7.18 & 29.43 \\
$h_3$ &  8.13 &  7.81 & 28.84  &      20.78 & 18.29 & 40.92
\end{tabular}

\bigskip

\caption{2d example, isotropic energy split: Total wall time per degree of freedom (in milliseconds)}
\label{tbl:2d_iso_times}
\end{table}

Figure~\ref{fig:walltime_2d_isotropic} and Table~\ref{tbl:2d_iso_times} also show
the wall-times for \TNNMGPRE, the TNNMG variant smoothing with an inexact local solver.
It can be observed that using that smoother decreases the computation time
by roughly further 5\% to 30\%.  This is the effect of precomputing the
local preconditioned Hessians in \TNNMGPRE{}.

Finally, we point out that the AT-1 model is solved more quickly than
the AT-2 one, even though the threshold for damage formation
makes it more challenging.

\subsubsection{Spectral splitting of the elastic energy density}
\begin{figure}

 \begin{tikzpicture}[scale=0.8]
  % AT-1 model
  \node at (0,10) {\includegraphics[width = 0.32\textwidth]{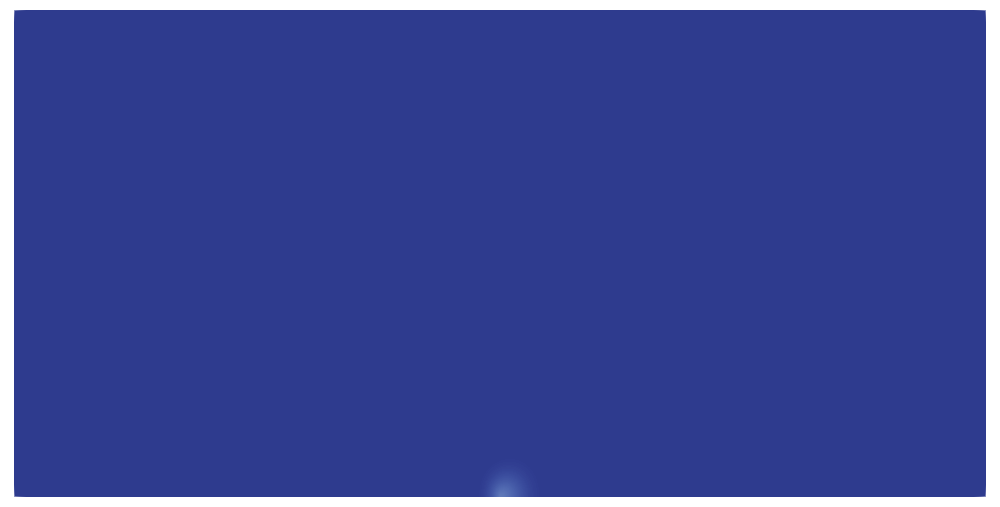}};
  \node at (0, 7) {\includegraphics[width = 0.32\textwidth]{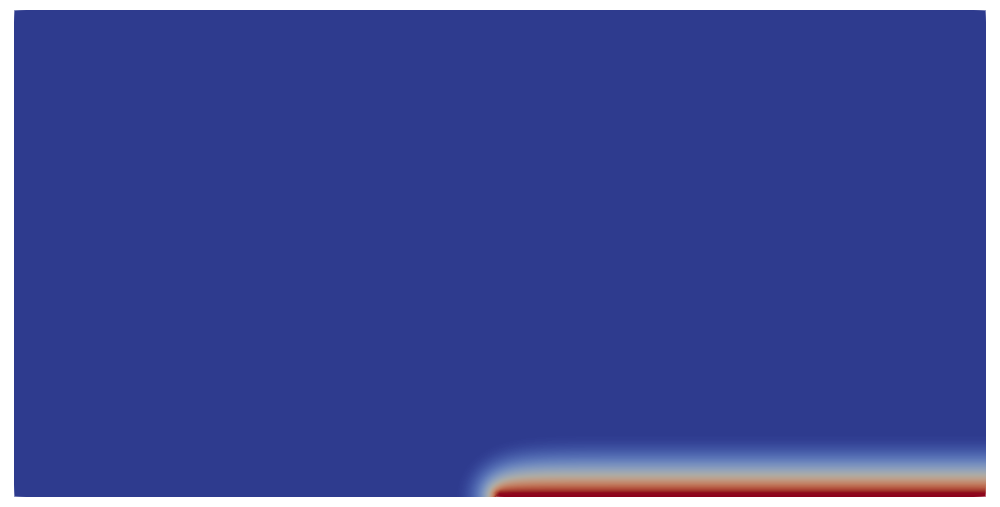}};

  \node at (8,8.5) {\begin{minipage}{0.5\textwidth}
                    \input{2d_reac_at1_aniso.pgf}
                   \end{minipage}};

  \node at (4,6) {AT-1};

  % AT-2 model
  \node at (0, 3) {\includegraphics[width = 0.32\textwidth]{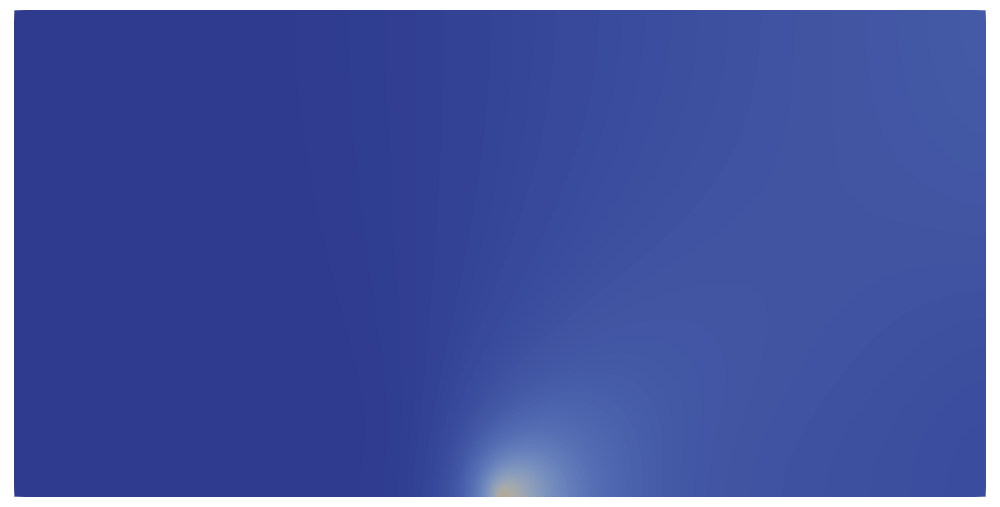}};
  \node at (0, 0) {\includegraphics[width = 0.32\textwidth]{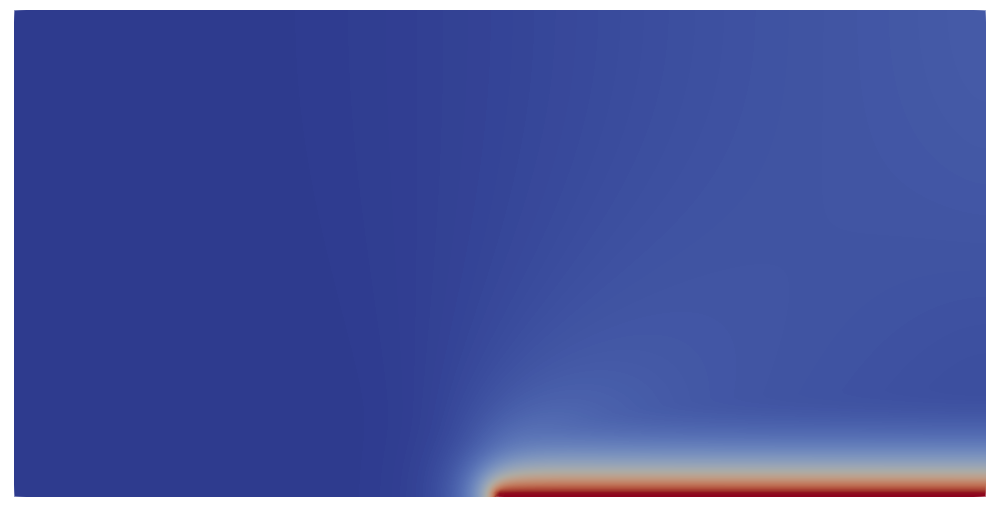}};

  \node at (8, 1.5) {\begin{minipage}{0.5\textwidth}
                    \input{2d_reac_at2_aniso.pgf}
                   \end{minipage}};

  \node at (4,-1) {AT-2};
 \end{tikzpicture}

 \caption{2d example, spectral energy split: Evolution at time steps~145 and~160, and displacement--force curves for the AT-1 (top) and AT-2 (bottom) models}
 \label{fig:evolution2daniso}
\end{figure}

For the next experiment we exchange the isotropic energy splitting~\eqref{eq:symmetric_damaged_energy}
by the spectral one~\eqref{eq:spectral_energy_split}.  As the example specimen
is loaded in tension we expect few differences to the previous experiments,
and indeed, the displacement--load curves (in Figure~\ref{fig:evolution2daniso})
show only minor differences compared to the ones of the isotropic splitting
(Figure~\ref{fig:evolution2diso}).

\begin{figure}
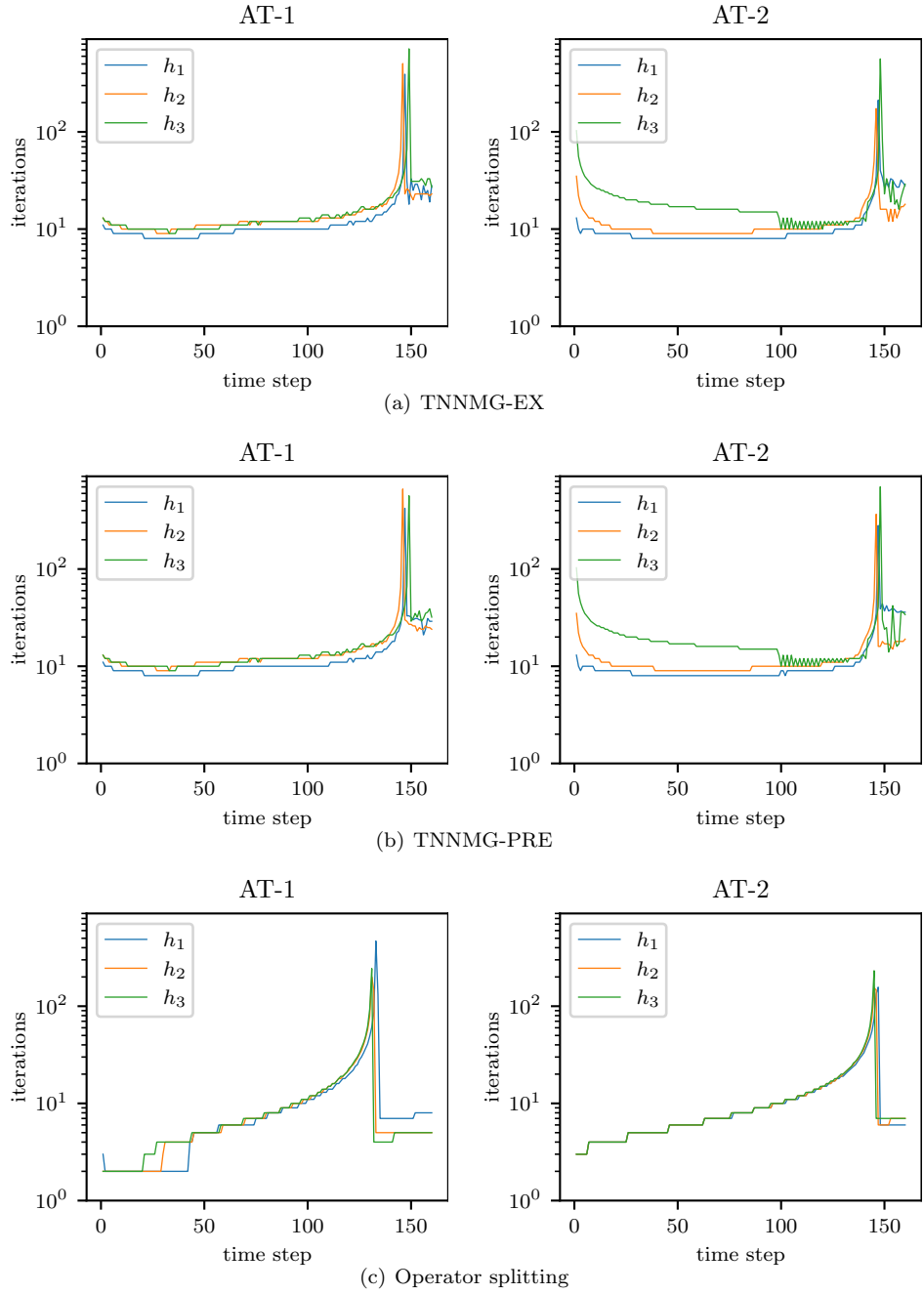


 \subfigure[\TNNMGEX]
 {
  \begin{minipage}{0.5\textwidth}
     \begin{center}
     \input{iter_2d_AT1_aniso_ex.pgf}
    \end{center}
  \end{minipage}
  \begin{minipage}{0.5\textwidth}
     \begin{center}
     \input{iter_2d_AT2_aniso_ex.pgf}
    \end{center}
  \end{minipage}
 }

 \subfigure[\TNNMGPRE]
 {
  \begin{minipage}{0.5\textwidth}
     \begin{center}
     \input{iter_2d_AT1_aniso_pre.pgf}
    \end{center}
  \end{minipage}
  \begin{minipage}{0.5\textwidth}
     \begin{center}
     \input{iter_2d_AT2_aniso_pre.pgf}
    \end{center}
  \end{minipage}
 }

 \subfigure[Operator splitting]
 {
  \begin{minipage}{0.5\textwidth}
     \begin{center}
     \input{iter_2d_AT1_aniso_feap.pgf}
    \end{center}
  \end{minipage}
  \begin{minipage}{0.5\textwidth}
     \begin{center}
     \input{iter_2d_AT2_aniso_feap.pgf}
    \end{center}
  \end{minipage}
 } \caption{2d example, spectral energy split: Iterations per
  time step, for grid sizes $h_1$, $h_2$, $h_3$}
 \label{fig:iterations_2d_anisotropic}
\end{figure}

The purpose of this test is primarily to assess the cost of the two-dimensional eigenvalue
decomposition and its derivatives required for the spectral splitting.
Figure~\ref{fig:iterations_2d_anisotropic} shows the iteration numbers per time step
for the three methods.  As the model is not quadratic in the displacement anymore,
we now distinguish between the \TNNMGEX and the \TNNMGPRE smoothers.
Not surprisingly though, the iteration numbers are virtually identical.
The iterations needed by the operator splitting method,
shown in the lowest row of Figure~\ref{fig:iterations_2d_anisotropic}, have not changed
appreciably either.
As an exception to this, the \TNNMGPRE method on the $h_3$ grid needs a much larger number
of iterations at the actual rupture, where the problem is very ill-conditioned.
While this is only a single load step, it markedly influences the accumulated
wall-times.

\begin{figure}
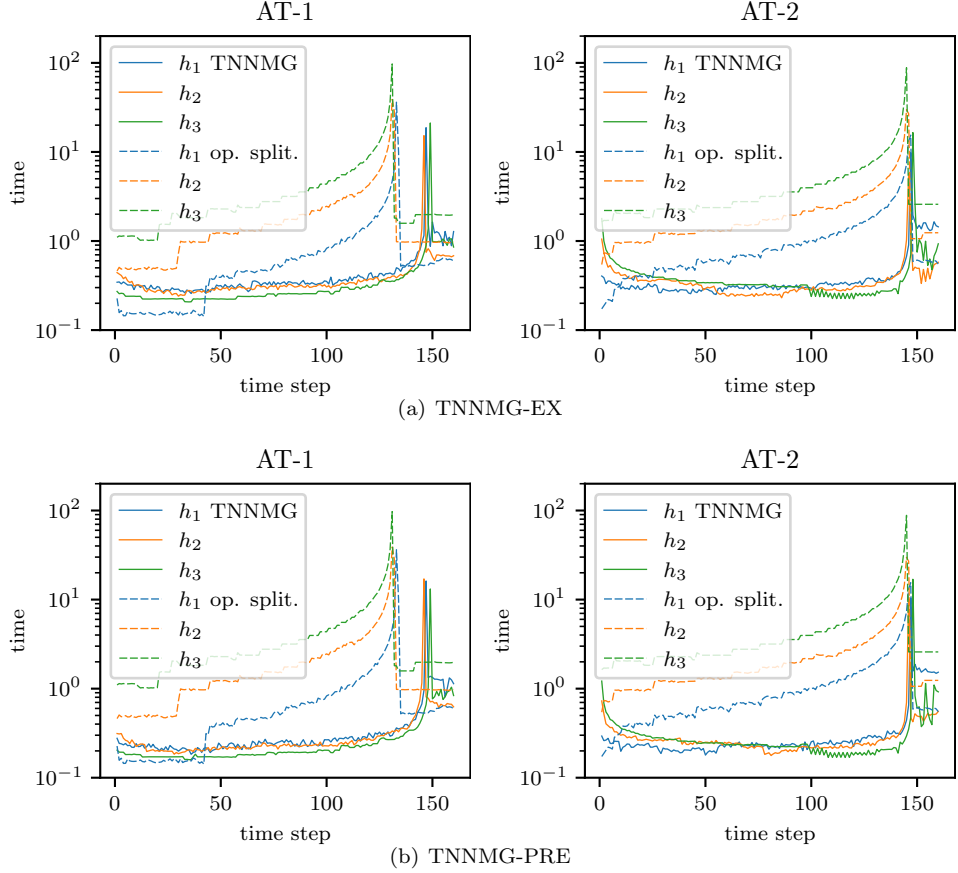


 \subfigure[\TNNMGEX]
 {
  \begin{minipage}{0.5\textwidth}
  \begin{center}
     \input{time_2d_AT1_aniso_ex.pgf}
  \end{center}
  \end{minipage}
  \begin{minipage}{0.5\textwidth}
  \begin{center}
     \input{time_2d_AT2_aniso_ex.pgf}
  \end{center}
  \end{minipage}
 }

 \subfigure[\TNNMGPRE]
 {
  \begin{minipage}{0.5\textwidth}
  \begin{center}
     \input{time_2d_AT1_aniso_pre.pgf}
  \end{center}
  \end{minipage}
  \begin{minipage}{0.5\textwidth}
  \begin{center}
     \input{time_2d_AT2_aniso_pre.pgf}
  \end{center}
  \end{minipage}
 }
 \caption{2d example, spectral energy split: Wall-time per degree of freedom
  per time step, for grid sizes $h_1$, $h_2$, $h_3$}
 \label{fig:walltime_2d_anisotropic}
\end{figure}

Figure~\ref{fig:walltime_2d_anisotropic} shows the time
needed to solve the increment problems, again normalized by the number of
degrees of freedom. One can see that the overall behavior
remains unchanged, but that the time needed by TNNMG goes up a
bit,
in particular for the AT-2 model.  As the number of iterations per load step
remains largely unchanged, the observed cost increase is caused by the
spectral decompositions performed by the smoother.
Interestingly, the operator-splitting solver for the AT-2 problem is a bit
faster during the rupturing of the specimen than it was for the
isotropically split energy.

\begin{table}
% TNNMG: time_2d_AT1_aniso_ex.dat [10.20155562 19.58246474 16.25969667]
% TNNMG: time_2d_AT1_aniso_pre.dat [ 9.09193507 11.696225   21.78786003]
% FEAP:  time_2d_AT1_aniso_feap.dat [20.99192552 21.75960033 32.83605376]
%
% TNNMG: time_2d_AT2_aniso_ex.dat [14.78129269 23.55445145 39.63698575]
% TNNMG: time_2d_AT2_aniso_pre.dat [11.16711469 16.78417387 28.64799727]
% FEAP:  time_2d_AT2_aniso_feap.dat [29.58519068 29.6798879  40.95810678]
\small
 \begin{tabular}{c|ccc|ccc}
 \hline
& \multicolumn{3}{c}{AT-1} & \multicolumn{3}{c}{AT-2} \\
&  \TNNMGEX & \TNNMGPRE & OS & \TNNMGEX & \TNNMGPRE & OS \\
\hline
$h_1$ & 10.20 &  9.09 & 20.99  &      14.78 & 11.17 & 29.59 \\
$h_2$ & 19.58 & 11.70 & 21.76  &      23.55 & 16.78 & 29.68 \\
$h_3$ & 16.26 & 21.79 & 32.84  &      39.64 & 28.65 & 40.96
\end{tabular}

\bigskip

\caption{2d example, spectral energy split: Total wall time per degree of freedom (in milliseconds)}
\label{tbl:2d_spectral_times}
\end{table}

When looking at the accumulated run-times shown in Table~\ref{tbl:2d_spectral_times}
(still normalized by the number of degrees of freedom), we see that the time needed
by the operator-splitting method remains virtually unchanged.
The TNNMG methods have gotten slower, though, and have only a small
wall-time lead over operator-splitting.  As in the case of the isotropic splitting,
moving from the $h_2$ grid to the $h_3$ grid increases the time per degree of freedom
considerably for the AT-2 model.  The reason is unclear.

Using the preconditioned smoother instead of the exact one now leads to a significant
decrease in the accumulated wall-time, with time reductions between 10\% and 40\%.
This is caused by the fact that the preconditioned smoother computes much fewer
eigenvector decompositions.  No such speedup can be seen for the AT-1 model
on the $h_3$ grid, where the preconditioned smoother is actually slower
than the exact one. As mentioned, this is because the algorithm needs many
more iterations at the load step with the rupture in this case.
The reason is unknown.

\subsubsection{Comparison to an interior-point solver from the literature}

To allow further assessment of the solver wall-times, we compare them with
results published by~\citet[Chapter~4.1]{wambacq_ulloa_lombaert_francois:2021}, that use
an interior-point solver for a very
similar example.  There, the authors consider the same problem geometry
and loading (quantitatively), and similar material parameters.
They model the material with an AT-2 functional and an elastic energy using
Amor's volumetric--deviatoric splitting~\cite{amor_marigo_maurini:2009},
see Section~\ref{sec:volumetric_decompositions}.
Note that this splitting is cheaper to evaluate than the eigenvalue-based one
whose times are given in Table~\ref{tbl:2d_spectral_times}.

As we do, the authors simulate a linearly increasing displacement load
until a bit after complete rupture, but they use only 60 load steps
compared to our 160 ones.  Apparently they simulate the entire square,
but with 17\,421 vertices, their grid is a bit coarser than our $h_1$ grid
(which has 33\,153 vertices).  They give cumulative simulation times for
four types of solvers, of which three are of operator-splitting type
(with different local solvers), and one is a monolithic interior-point solver.
For the operator-splitting solvers they report wall-times per degree of freedom
between 6.9\,ms and 27.6\,ms.  For a better comparison with our simulation
with 160 load steps we multiply these times with $\frac{160}{60}$ and obtain
values between 18.4\,ms and 73.5\,ms.  Likewise, for the monolithic solver
they report a cumulative time per degree of freedom of 82.7\,ms (for 60 load steps).
Scaled to 160 load steps this gives a value of about 220\,ms.
These times should be compared to the ones of the AT-2 column of
Tables~\ref{tbl:2d_iso_times} and~\ref{tbl:2d_spectral_times}.
One can see that our implementations are quite a bit faster.

Note, however, that the two example problems are not exactly the same,
that termination criteria differ,
and that both hardware and implementation differ as well.  Also, the grid used
by~\citet{wambacq_ulloa_lombaert_francois:2021} is coarser than ours.  Therefore,
the comparison should not be over-interpreted.
On the other hand, the interior-point implementation involves sparse matrix
factorizations, and it is therefore to be expected that the run-times
deteriorate rapidly for increasing grid resolution,
in particular for three-dimensional problems.

\subsection{Bending test of a notched bar in three dimensions}

The second example uses a three-dimensional object.  In three dimensions, stiffness
matrices get denser, and hence direct solvers for global linear systems
get more expensive and need considerably more memory.  In contrast,
the TNNMG memory consumption remains linear in the number of unknowns.
Also, eigenvalue decompositions are more expensive for $3 \times 3$ matrices,
and we therefore expect larger run-time differences between the isotropic
and the spectrally split model, and between the two TNNMG smoothers.

\begin{figure}
  \begin{center}
    \begin{overpic}[width = 0.8\textwidth]{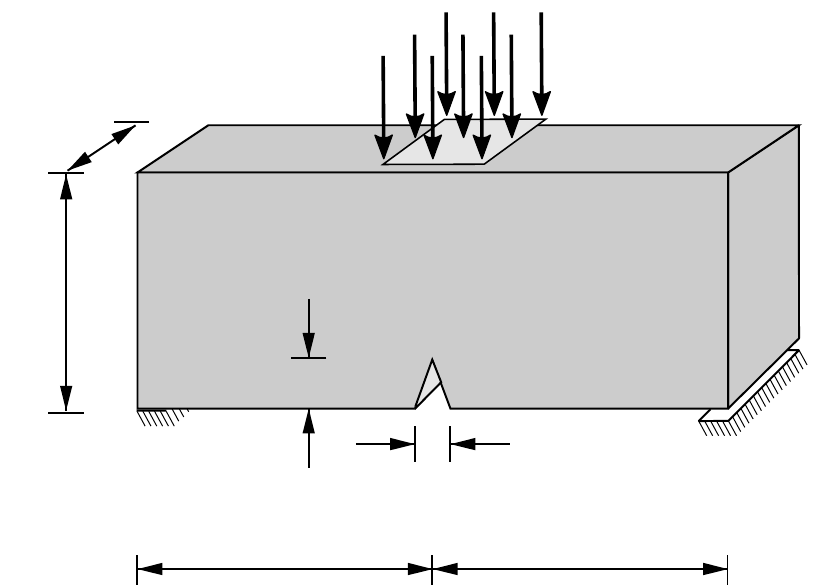}
     \put( 6,53){$L_z$}
     \put( 3,35){$L_y$}
     \put(68,64){$\bar \Bu$}
     \put(33, 5){$L_x/2$}
     \put(68, 5){$L_x/2$}
     \put(51,12){$l_1$}
     \put(32,23.5){$l_2$}
    \end{overpic}

    \caption{Boundary value problem of bending test in three dimensions}
    \label{fig:bending3d1}
  \end{center}
\end{figure}

We consider a bending test for a rectangular bar with a triangular notch.
In this setting, the decomposition of the elastic energy density plays a
crucial role.  Under the given loading, parts of the specimen undergo severe compression, and material
models that degrade under such compression will therefore show unphysical results.
We test the solvers with the isotropic splitting~\eqref{eq:symmetric_damaged_energy}
nevertheless, to obtain an idea of the cost of the spectral splitting.

The example setting is again taken from~\cite{miehe+welschinger+hofacker10a}.
The geometry and the boundary conditions are
visualized in Figure \ref{fig:bending3d1}. The dimensions of the
specimen are $L_x=8$\,mm, $L_y=2$\,mm and $L_z=1$\,mm.
Width and height of the triangular notch are $l_1=0.4$\,mm and $l_2=0.2$\,mm, respectively.
We use three unstructured hexahedral grids to discretize the domain,%
with 1\,920, 15\,360, and 122\,880 elements, respectively.
In the figures these are denoted by $h_1$, $h_2$, and $h_3$, respectively.
The grids were constructed by 2, 3, and~4 steps of uniform refinement
of a coarse grid with 30 elements, respectively,
and this refinement
history is used by the linear geometric multigrid step of the TNNMG algorithm.
The grids are graded a bit towards the expected fracture,
and have an average edge length of about $l/2.3$ (for $h_1$), $l/4.6$ (for $h_2$)
and $l/9.2$ (for $h_3$) there.

As in~\cite{miehe+welschinger+hofacker10a},
the Lamé parameters are set to $\lambda=121$\,kN/mm$^2$ and $\mu=80$\,kN/mm$^2$.
For the further parameters we set
$g_c=2.7$\,N/mm, $l=0.2$\,mm, and $k=10^{-5}$.
The object is loaded with a time-dependent displacement load in downward
direction in a strip of width 1.2\,mm in the center of the top surface.
In that strip, the surface is fixed in $z$-direction, but free to move
in $x$-direction. The displacement load is set to $\bar{\bu}_i = -i \cdot 5 \cdot 10^{-3}\,\mathrm{mm} \cdot \be_3$
for the load steps $i=1,\dots,13$.  No adaptive load-stepping as in~\cite{miehe+welschinger+hofacker10a}
is necessary, because time discretization and solvers are stable enough for this
range of load step sizes. The object is fully clamped at a strip of
width 0.4\,mm at the left end of the lower surface.  At the right end
of the same surface, a strip of width 0.4\,mm is fixed in the $y$ and $z$ directions.

For each time step, we solve the increment problem with the same solver settings as for the
two-dimensional example: Each iteration starts at the solution of the previous load step,
and we terminate when the degraded energy norm of the correction, scaled by the
degraded energy norm of the previous iterate, drops below $10^{-7}$.
The linear correction step~\eqref{eq:inexact_linear_correction} consists of three geometric
multigrid $V(3,3)$-cycle steps with a block-Gauß--Seidel smoother,
and damage degrees
of freedom are truncated when their current value is less than $10^{-10}$
away from the lower obstacle.

\subsubsection{Isotropic splitting}
\begin{figure}
  \begin{tikzpicture}[scale=0.8]
  % AT-1 model
  \node at (0,10) {\includegraphics[width = 0.35\textwidth]{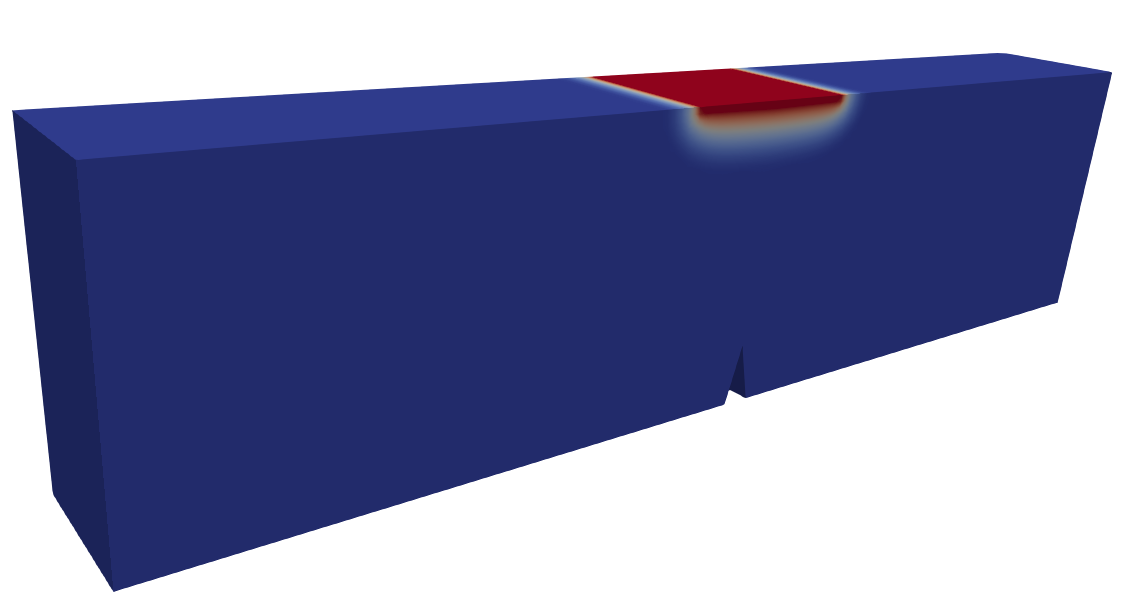}};
  \node at (0, 7) {\includegraphics[width = 0.35\textwidth]{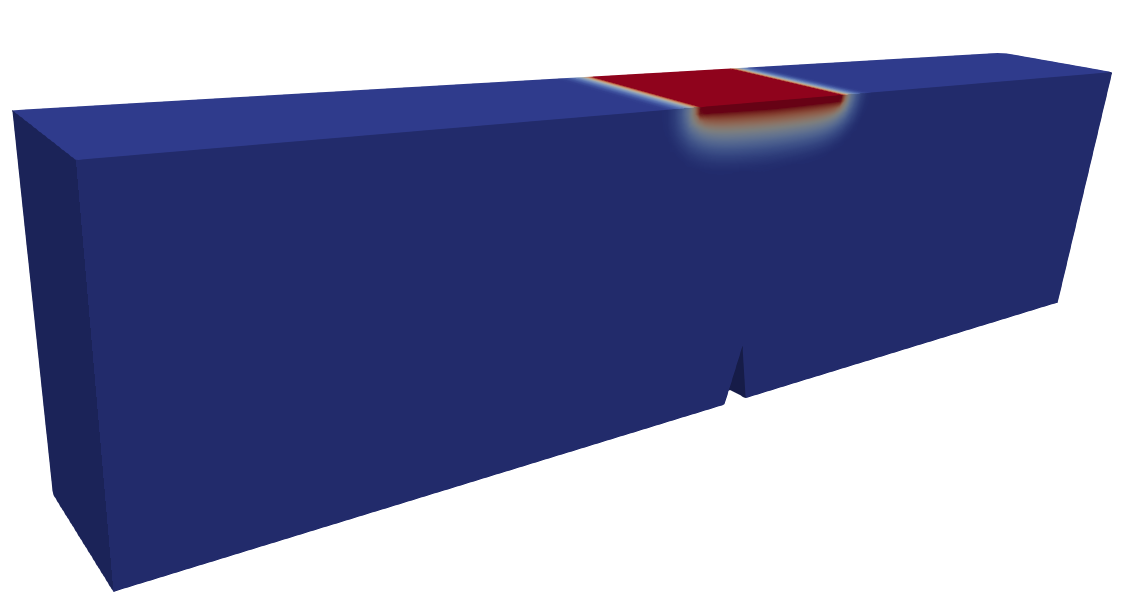}};

  \node at (8,8.5) {\begin{minipage}{0.5\textwidth}
                    \input{3d_reac_at1_iso.pgf}
                   \end{minipage}};

  \node at (4,6) {AT-1};

  % AT-2 model
  \node at (0, 3) {\includegraphics[width = 0.35\textwidth]{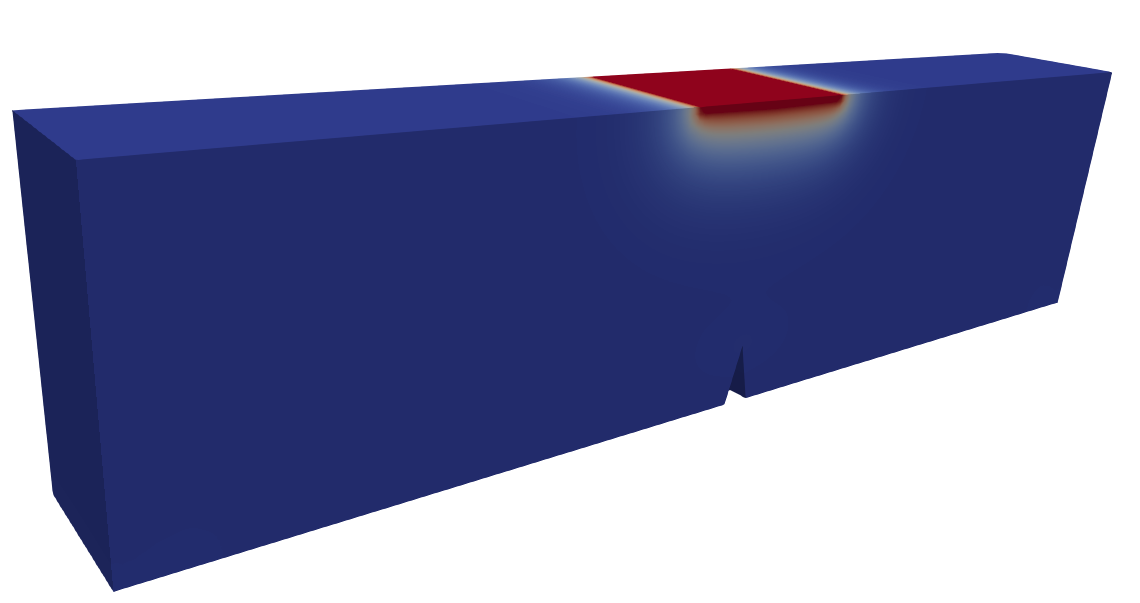}};
  \node at (0, 0) {\includegraphics[width = 0.35\textwidth]{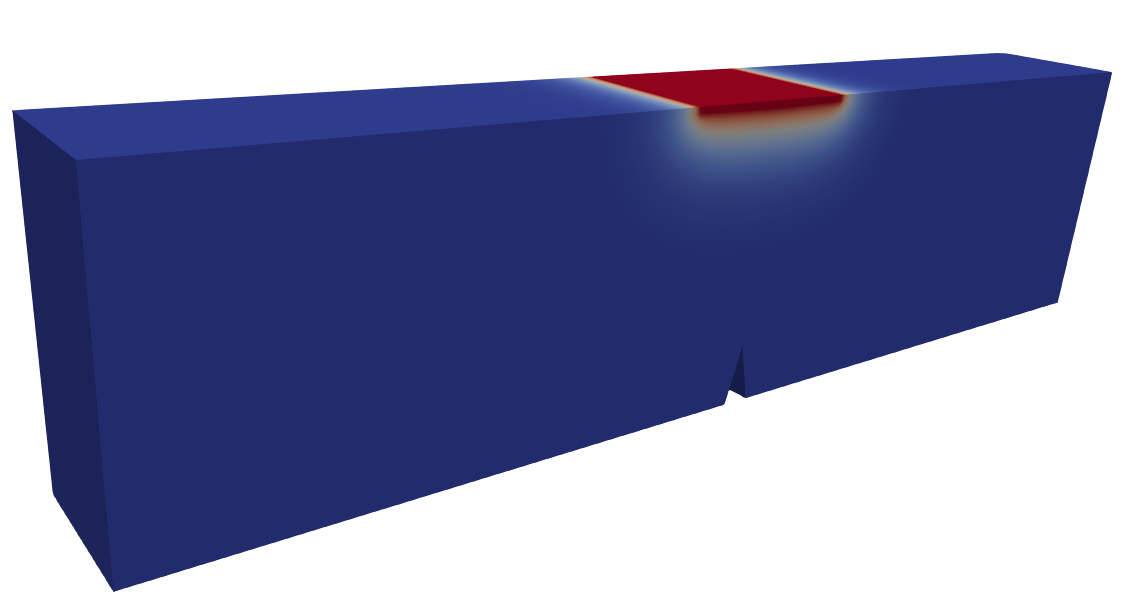}};

  \node at (8, 1.5) {\begin{minipage}{0.5\textwidth}
                    \input{3d_reac_at2_iso.pgf}
                   \end{minipage}};

  \node at (4,-1) {AT-2};
 \end{tikzpicture}

 \caption{3d example, isotropic energy split: Configurations at load steps~7 and~13,
   and displacement--force curves.}
 \label{fig:evolution_3d_isotropic}
\end{figure}

We first consider the model with the isotropic splitting~\eqref{eq:symmetric_damaged_energy}
of the elastic energy density, i.e., the model where all elastic strains contribute to the degradation
of the material.  We expect unphysical results: Virtually all damage will happen
in the vicinity of the load, whereas the region around the notch will remain intact.
Indeed, the simulation results in Figure~\ref{fig:evolution_3d_isotropic}
show that this is exactly what is happening.
We are nevertheless interested in the solver behavior for this model, primarily
to assess the cost of the spectral splitting in the
next section. No costly eigenvalue decompositions are necessary for the isotropic
splitting considered here, and
the local displacement problems~\eqref{eq:local_displacement_increment_energy}
are quadratic. Consequently, the two smoother variants \TNNMGEX and \TNNMGPRE solve almost
the same problems, and we expect them to behave more or less the same, too,
with a small run-time advantage for \TNNMGPRE.

To assess solver performance,
we again first compare iteration numbers.  Figure~\ref{fig:iterations_3d_isotropic} shows the
iteration numbers per time step needed by the two TNNMG variants, and by the
operator-splitting iteration. Note again that the vertical axis is scaled logarithmically.
In contrast to the two-dimensional experiment we did notice different iteration
numbers for the exact and preconditioned smoothers, and we therefore show both plots.

We see that TNNMG needs about 10~iterations for each load step for the AT-1 model.
For the AT-2 model, for which the damage variable changes throughout the domain
immediately, the solver also needs about 10 iterations on the $h_2$ grid,
but almost 10 times as many for the other two grids for the first 5 load steps.
Operator-splitting, on the other hand, behaves roughly the same for both models.
Unlike TNNMG, it needs less iterations initially, but many more later on.
More specifically, the method needs around 10~iterations or even less for the first 5 load steps,
but after that it consistently needs about 100~iterations per load step
for all grids, with one outlier even going up to 1000~iterations. For an attempt
of an explanation, recall that in this example the specimen never breaks into
two parts (Figure~\ref{fig:evolution_3d_isotropic}), and the problem therefore
remains much better conditioned than the previous ones. In this situation,
the monolithic TNNMG method seems to have a clear advantage.

The difference in iteration numbers together with the lower cost of TNNMG iterations
compared to operator-splitting iterations leads to a tremendous speed difference.
As shown in Figure~\ref{fig:walltime_3d_isotropic}, for all load steps beyond
the first few, the TNNMG method needs about two orders of magnitude less
wall-time than the operator-splitting method. This difference can also be seen
in Table~\ref{tbl:3d_isotropic_times}, which shows the accumulated wall-times
for the entire load history, normalized by the number of degrees of freedom.
We see that \TNNMGEX is about 25 to 35 times faster than operator-splitting
for the AT-1 model, and 25 to 30 times faster for the AT-2 model.
Using the preconditioned smoother makes the wall-time decrease further.
In this three-dimensional situation, not having to recompute the matrix
block diagonal entries does lead to large savings. The speed advantage of \TNNMGPRE
over operator-splitting rises to about 42--58 for the AT-1 model, and to
33--57 for the AT-2 one.

\begin{figure}
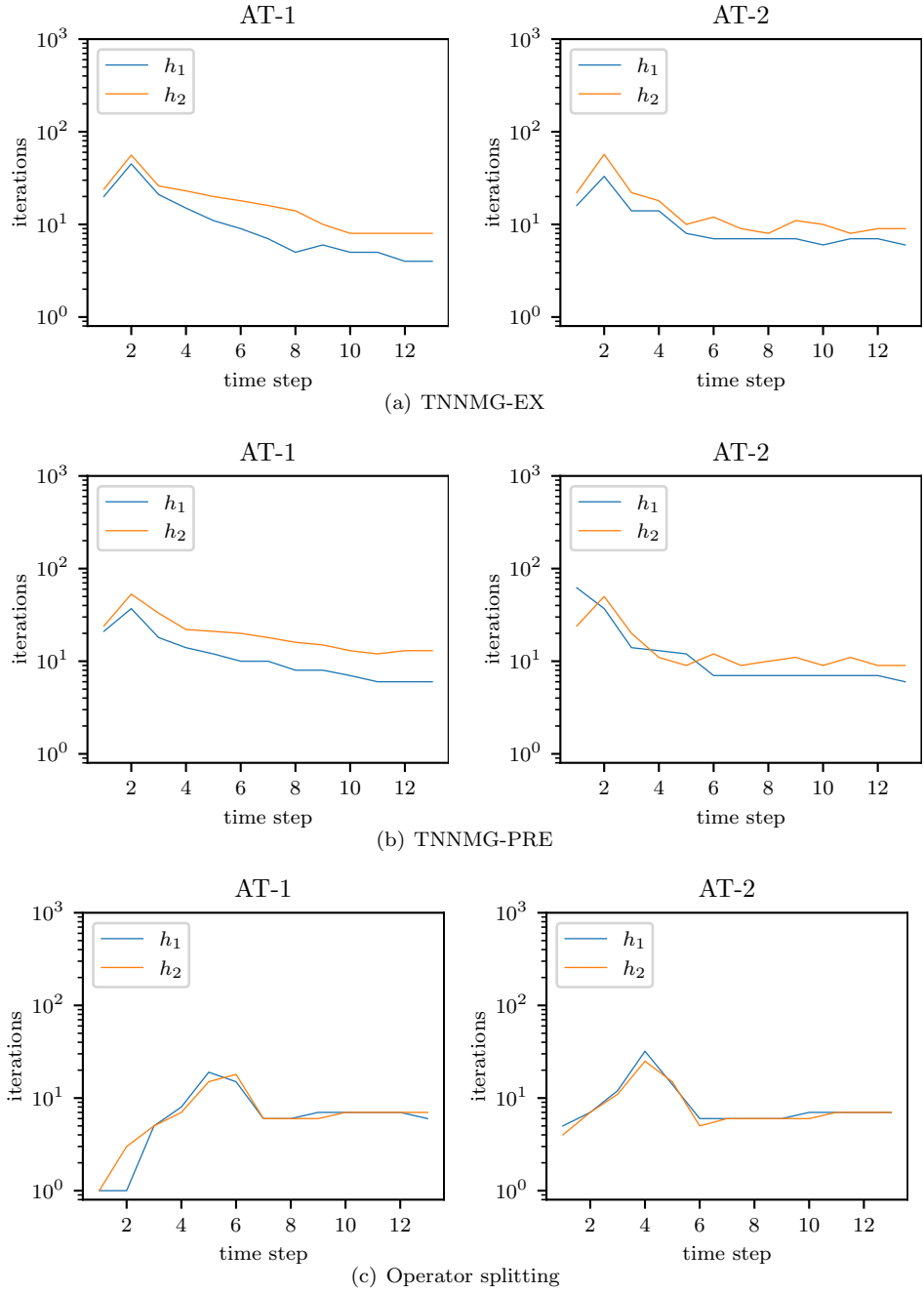


 \subfigure[\TNNMGEX]
 {
  \begin{minipage}{0.5\textwidth}
    \begin{center}
    \input{iter_3d_AT1_iso_ex.pgf}
    \end{center}
  \end{minipage}
  \begin{minipage}{0.5\textwidth}
    \begin{center}
    \input{iter_3d_AT2_iso_ex.pgf}
    \end{center}
  \end{minipage}
 }

 \subfigure[\TNNMGPRE]
 {
  \begin{minipage}{0.5\textwidth}
    \begin{center}
    \input{iter_3d_AT1_iso_pre.pgf}
    \end{center}
  \end{minipage}
  \begin{minipage}{0.5\textwidth}
    \begin{center}
    \input{iter_3d_AT2_iso_pre.pgf}
    \end{center}
  \end{minipage}
 }

 \subfigure[Operator splitting]
 {
  \begin{minipage}{0.49\textwidth}
  \begin{center}
     \input{iter_3d_AT1_iso_feap.pgf}
  \end{center}
  \end{minipage}
  \begin{minipage}{0.49\textwidth}
  \begin{center}
     \input{iter_3d_AT2_iso_feap.pgf}
  \end{center}
  \end{minipage}
 }
 \caption{3d example, isotropic energy split: Iterations per time step, for grid sizes $h_1$, $h_2$, $h_3$}
 \label{fig:iterations_3d_isotropic}
\end{figure}

\begin{figure}
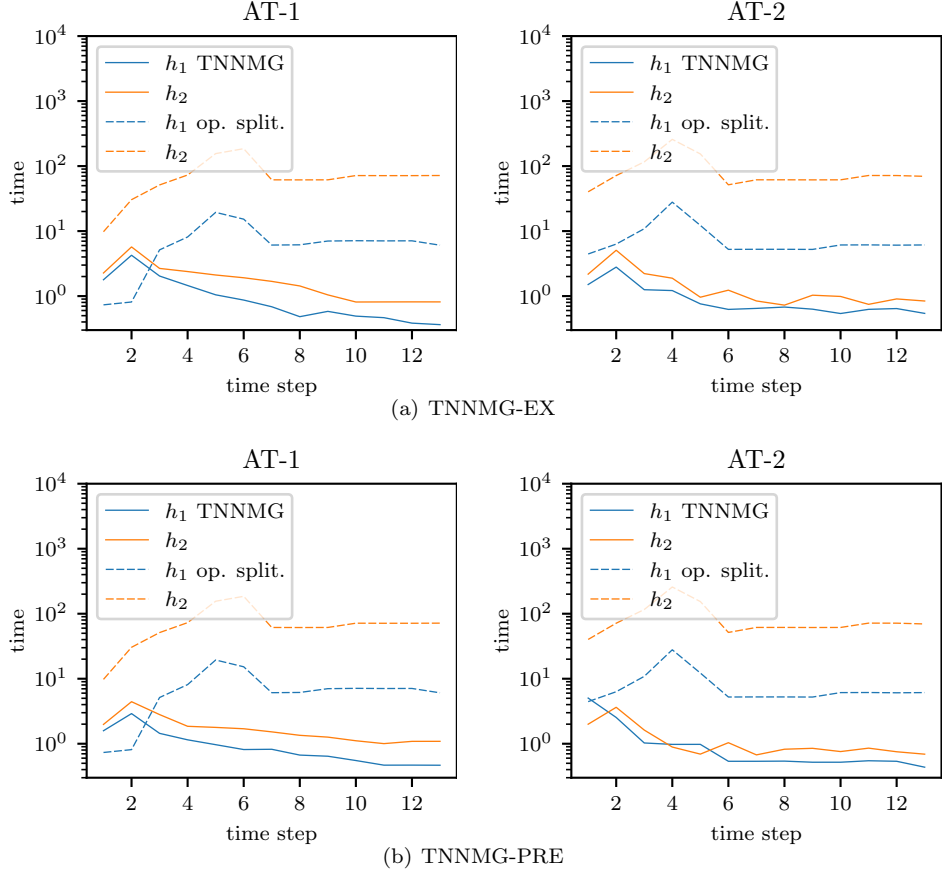


 \subfigure[\TNNMGEX]
 {
  \begin{minipage}{0.5\textwidth}
    \begin{center}
    \input{time_3d_AT1_iso_ex.pgf}
    \end{center}
  \end{minipage}
  \begin{minipage}{0.5\textwidth}
    \begin{center}
    \input{time_3d_AT2_iso_ex.pgf}
    \end{center}
  \end{minipage}
 }

 \subfigure[\TNNMGPRE]
 {
  \begin{minipage}{0.5\textwidth}
    \begin{center}
    \input{time_3d_AT1_iso_pre.pgf}
    \end{center}
  \end{minipage}
  \begin{minipage}{0.5\textwidth}
    \begin{center}
    \input{time_3d_AT2_iso_pre.pgf}
    \end{center}
  \end{minipage}
 }
 \caption{3d example, isotropic energy split: Wall-time per degree of freedom over
 time step, for grid sizes $h_1$, $h_2$, $h_3$}
 \label{fig:walltime_3d_isotropic}
\end{figure}

\begin{table}
% TNNMG: time_3d_AT1_iso_ex.dat [ 4.14798848  6.11304487 10.34157558]
% TNNMG: time_3d_AT1_iso_pre.dat [2.63129846 3.73710048 6.2753344 ]
% FEAP:  time_3d_AT1_iso_feap.dat [109.77075528 216.41605447 357.35397458]
%
% TNNMG: time_3d_AT2_iso_ex.dat [13.30851344  7.12382975 14.35879652]
% TNNMG: time_3d_AT2_iso_pre.dat [11.93951919  4.07215559  9.40214127]
% FEAP:  time_3d_AT2_iso_feap.dat [393.91832054 223.96226887 349.97077954]
\small
 \begin{tabular}{c|ccc|ccc}
 \hline
& \multicolumn{3}{c}{AT-1} & \multicolumn{3}{c}{AT-2} \\
&  \TNNMGEX & \TNNMGPRE & OS & \TNNMGEX & \TNNMGPRE & OS \\
\hline
$h_1$ &  4.15 &  2.63 & 109.77  &      13.31 & 11.94 & 393.92 \\
$h_2$ &  6.11 &  3.74 & 216.42  &       7.12 &  4.07 & 223.96 \\
$h_3$ & 10.34 &  6.28 & 357.35  &      14.36 &  9.40 & 349.97
\end{tabular}

\bigskip

\caption{3d example, isotropic energy split: Total wall time per degree of freedom (in milliseconds),
 for grid sizes $h_1$, $h_2$, $h_3$}
\label{tbl:3d_isotropic_times}
\end{table}

\subsubsection{Spectral splitting of the elastic energy density}
\begin{figure}
 \begin{tikzpicture}[scale=0.8]
  % AT-1 model
  \node at (0,10) {\includegraphics[width = 0.35\textwidth]{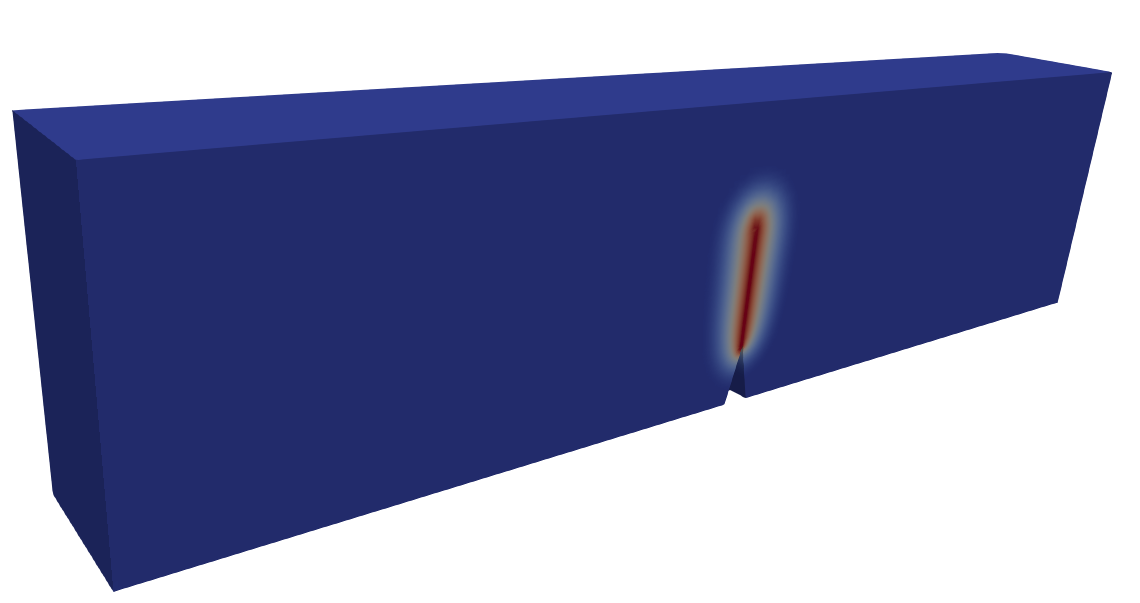}};
  \node at (0, 7) {\includegraphics[width = 0.35\textwidth]{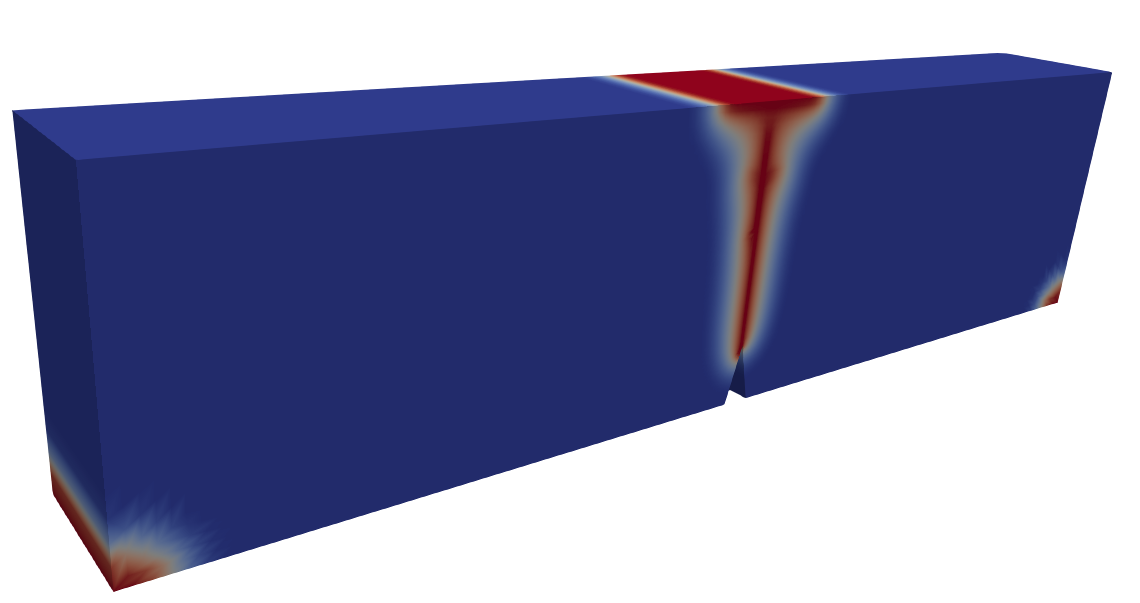}};

  \node at (8,8.5) {\begin{minipage}{0.5\textwidth}
                    \input{3d_reac_at1_aniso.pgf}
                   \end{minipage}};

  \node at (4,6) {AT-1};

  % AT-2 model
  \node at (0, 3) {\includegraphics[width = 0.35\textwidth]{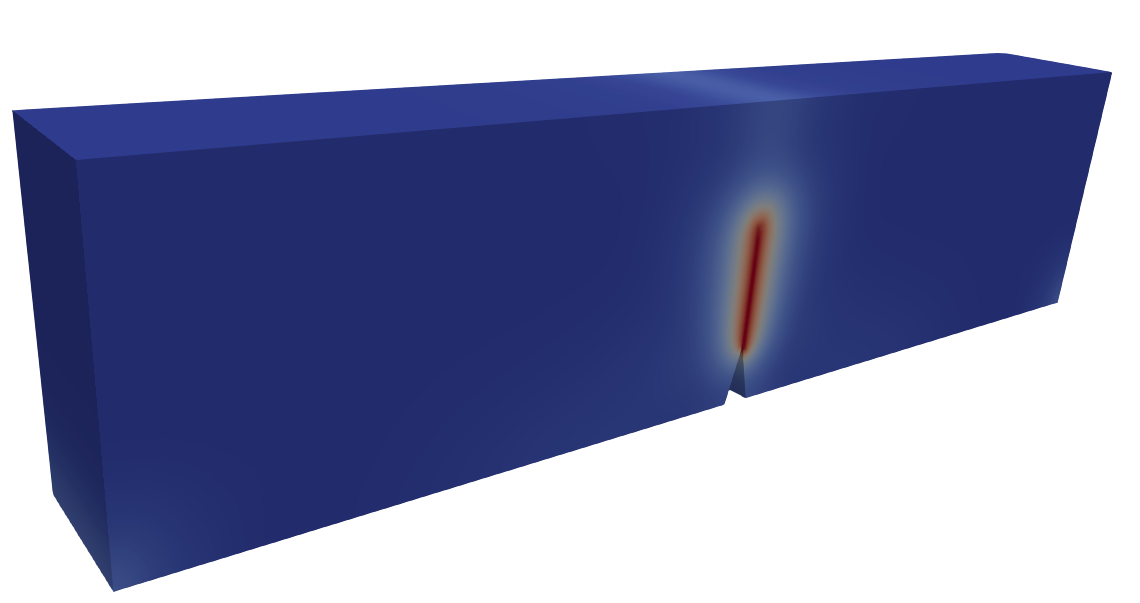}};
  \node at (0, 0) {\includegraphics[width = 0.35\textwidth]{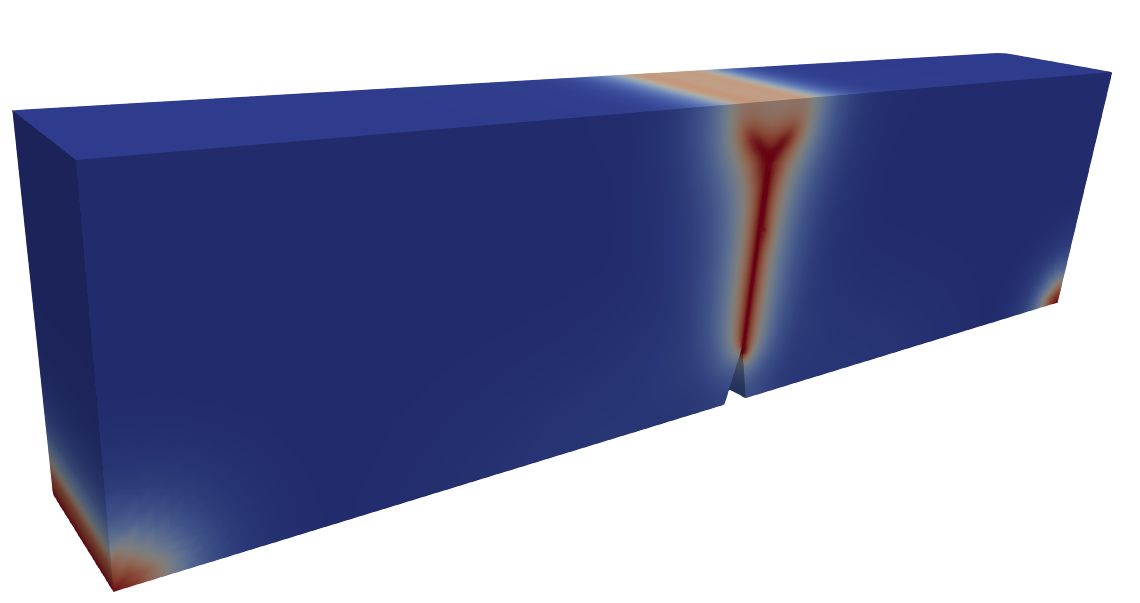}};

  \node at (8, 1.5) {\begin{minipage}{0.5\textwidth}
                    \input{3d_reac_at2_aniso.pgf}
                   \end{minipage}};

  \node at (4,-1) {AT-2};
 \end{tikzpicture}

 \caption{3d example, spectral energy split: Evolution at time steps~7 and~13, and displacement--force curves}
 \label{fig:3d_evolution_spectral}
\end{figure}

\begin{figure}
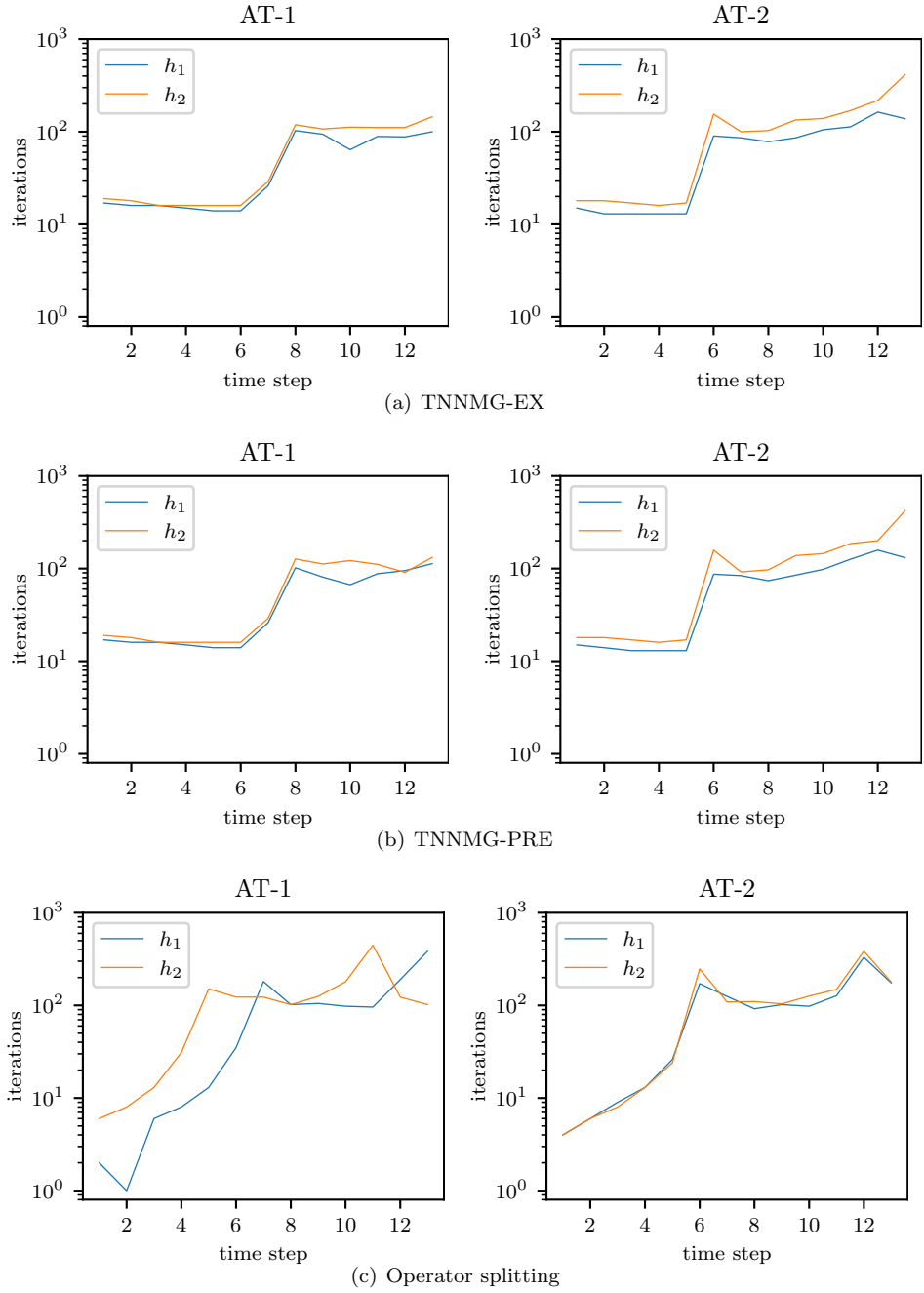


 \subfigure[\TNNMGEX]
 {
  \begin{minipage}{0.5\textwidth}
    \begin{center}
    \input{iter_3d_AT1_aniso_ex.pgf}
    \end{center}
  \end{minipage}
  \begin{minipage}{0.5\textwidth}
    \begin{center}
    \input{iter_3d_AT2_aniso_ex.pgf}
    \end{center}
  \end{minipage}
 }

 \subfigure[\TNNMGPRE]
 {
  \begin{minipage}{0.5\textwidth}
    \begin{center}
    \input{iter_3d_AT1_aniso_pre.pgf}
    \end{center}
  \end{minipage}
  \begin{minipage}{0.5\textwidth}
    \begin{center}
    \input{iter_3d_AT2_aniso_pre.pgf}
    \end{center}
  \end{minipage}
 }

 \subfigure[Operator splitting]
 {
  \begin{minipage}{0.49\textwidth}
  \begin{center}
     \input{iter_3d_AT1_aniso_feap.pgf}
  \end{center}
  \end{minipage}
  \begin{minipage}{0.49\textwidth}
  \begin{center}
     \input{iter_3d_AT2_aniso_feap.pgf}
  \end{center}
  \end{minipage}
 }
 \caption{3d example, spectral energy split: Iterations per time step, for grid sizes $h_1$, $h_2$, $h_3$}
 \label{fig:iterations_3d_spectral}
\end{figure}

\begin{figure}
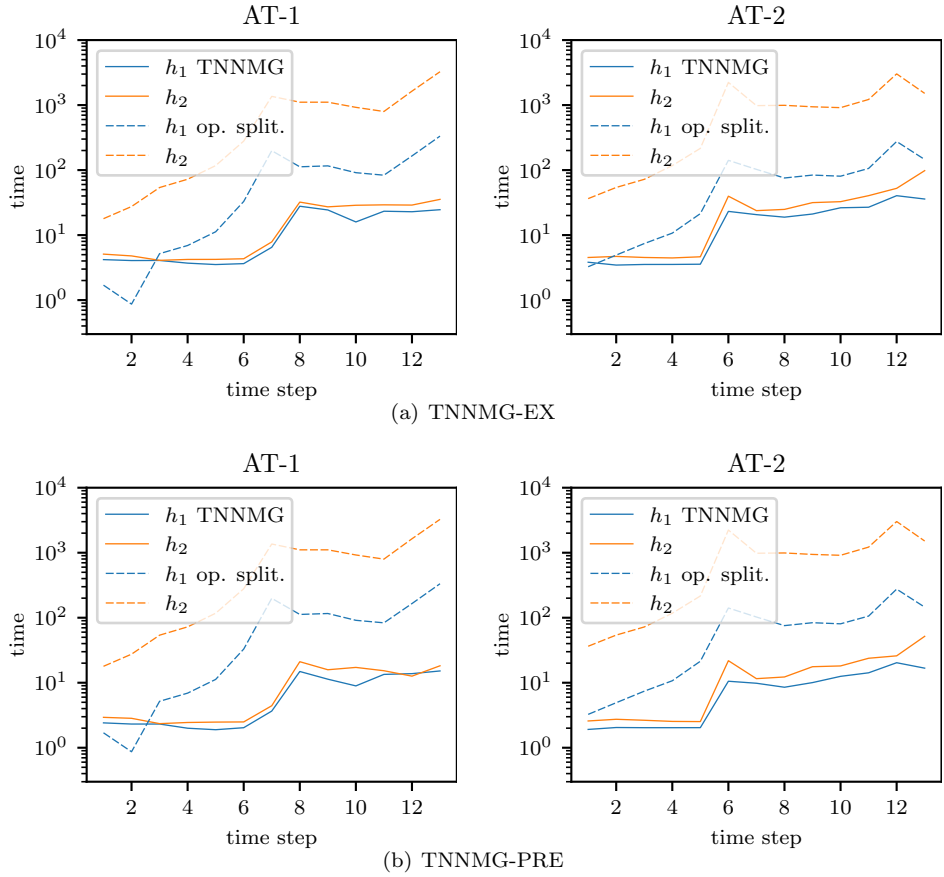


 \subfigure[\TNNMGEX]
 {
  \begin{minipage}{0.5\textwidth}
    \begin{center}
    \input{time_3d_AT1_aniso_ex.pgf}
    \end{center}
  \end{minipage}
  \begin{minipage}{0.5\textwidth}
    \begin{center}
    \input{time_3d_AT2_aniso_ex.pgf}
    \end{center}
  \end{minipage}
 }

 \subfigure[\TNNMGPRE]
 {
  \begin{minipage}{0.5\textwidth}
    \begin{center}
    \input{time_3d_AT1_aniso_pre.pgf}
    \end{center}
  \end{minipage}
  \begin{minipage}{0.5\textwidth}
    \begin{center}
    \input{time_3d_AT2_aniso_pre.pgf}
    \end{center}
  \end{minipage}
 }
 \caption{3d example, spectral energy split: Wall-time per degree of freedom over
  time step, for grid sizes $h_1$, $h_2$, $h_3$
  }
 \label{fig:walltime_3d_anisotropic}
\end{figure}

In the second set of experiments we replace the unsplit degraded energy
density~\eqref{eq:symmetric_damaged_energy} by the one with the spectral splitting according
to~\eqref{eq:spectral_energy_split}.
Figure~\ref{fig:3d_evolution_spectral} shows two snapshots from the problem evolution,
and the reaction force as a function of the applied displacement.
One can clearly see the differences to the simulation with the unsplit energy.  As one would expect,
the damage now happens primarily near the notch, and the
specimen does now break into two parts.

Figure~\ref{fig:iterations_3d_spectral}
shows the iteration numbers for TNNMG with both smoothers and the operator-splitting
method again.  Iteration numbers for the multigrid
algorithms and both AT models are roughly the same.  In particular,
\TNNMGPRE, the multigrid method with the inexact smoother avoiding the costly
$3 \times 3$ eigenvalue decomposition, does not need more iterations than \TNNMGEX.
The operator-splitting method needs slightly less iterations than TNNMG for the
first few load steps, and a few more after that.

When looking at wall-times, the
TNNMG algorithm is consistently faster than the operator-splitting method.
Figure~\ref{fig:walltime_3d_anisotropic} shows the
normalized times in the same way as in the previous sections.
The \TNNMGEX method is roughly five times faster in each load step.
This shows the effect of the cheaper iteration times: In three spaces dimensions,
tangent matrices are denser than in two dimensions, and direct solving
of linear systems with such matrices gets more expensive.
Consequently, operator-splitting iterations get relatively more expensive
than multigrid ones. This is also
reflected in Table~\ref{tbl:3d_spectral_times}, which shows the accumulated
run-times for the entire load history per degree of freedom.
Here, \TNNMGEX is about 4 to 7 times faster than operator-splitting for the
AT-1 model, and roughly 5 times faster for the AT-2 model. % 4.7 to 5.5

\begin{table}
% TNNMG: time_3d_AT1_aniso_ex.dat [23.52491459 38.25566939 64.75672933]
% TNNMG: time_3d_AT1_aniso_pre.dat [14.68976008 21.36459826 36.55489808]
% FEAP:  time_3d_AT1_aniso_feap.dat [124.63426104 275.53968571 266.59330105]
%
% TNNMG: time_3d_AT2_aniso_ex.dat [23.42464875 46.01720529 78.29074923]
% TNNMG: time_3d_AT2_aniso_pre.dat [12.87371305 25.07024404 38.31279024]
% FEAP:  time_3d_AT2_aniso_feap.dat [128.93370441 217.46916819 393.44561393]
\small
 \begin{tabular}{c|ccc|ccc}
 \hline
& \multicolumn{3}{c}{AT-1} & \multicolumn{3}{c}{AT-2} \\
&  \TNNMGEX & \TNNMGPRE & OS & \TNNMGEX & \TNNMGPRE & OS \\
\hline
$h_1$ & 23.52 & 14.69 & 124.63 &       23.42 & 12.87 & 128.93 \\
$h_2$ & 38.26 & 21.36 & 275.54 &       46.02 & 25.07 & 217.47 \\
$h_3$ & 64.76 & 36.55 & 266.59 &       78.29 & 38.31 & 393.45
\end{tabular}

\bigskip

\caption{3d example, spectral energy split: Total wall time per degree of freedom (in milliseconds),
 for grid sizes $h_1$, $h_2$, $h_3$}
\label{tbl:3d_spectral_times}
\end{table}

The speed difference gets even larger when using the preconditioned smoother.
Computing eigenvector decompositions (and the derivative transformation formulas
for spectral functions; Section~\ref{sec:multigrid_corrections})
is much more expensive for $3 \times 3$ matrices than for $2 \times 2$ ones.
Therefore, using the preconditioned smoother, which avoids many of these
computations, saves a lot of time. As can be seen from Figure~\ref{fig:walltime_3d_anisotropic}
and Table~\ref{tbl:3d_spectral_times}, \TNNMGPRE is about 7 to 12 times
faster for the AT-1 model, and 8 to 10 times faster for the AT-2 model.

\bigskip

Finally, we investigate the memory consumption of the two solver implementations.
Figure~\ref{fig:3d_spectral_memory} shows the maximum amount of memory used by the two algorithms
for the problem of this section, as a function of the grid size.
Memory consumption was measured using the \texttt{valgrind-massif} tool%
\footnote{\url{valgrind.org}}
with the \texttt{--pages-as-heap} option.
One can see that TNNMG clearly requires a linear amount
of memory. On the other hand, the proportionality factor is rather large, because
TNNMG needs the full Newton tangent matrix, restriction operators, and approximate
tangent matrices on coarser grid levels.
Also, in order to obtain first- and second-order derivatives efficiently during
the local solves, our implementation precomputes and stores the values and derivatives
of the shape functions at all quadrature points (what Miehe et al.\
call the \emph{global interpolation matrix} \cite[Chap.\,4.2]{miehe+hofacker+welschinger10}).

In contrast, operator splitting does not use coarser grid levels and only assembles
tangent matrices for the displacement and damage problems separately. This leads
to a smaller memory footprint for small grids. Once the grid gets finer, though,
the superlinear memory complexity of the direct solver begins to show.
To highlight this effect we measured the memory consumption for one further step
of grid refinement. The resulting grid has 983\,040 elements, and the simulation
using TNNMG needs about 47\,GB of memory.  The corresponding operator-splitting
implementation, however, would exhaust even a machine with 512\,GB of main memory
(Figure~\ref{fig:3d_spectral_memory}).

\begin{figure}
 \begin{minipage}{0.49\textwidth}
 \begin{center}
  \input{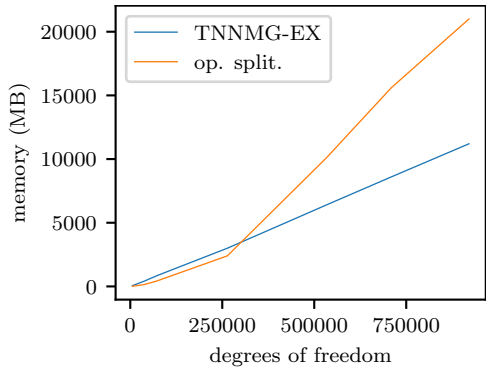}
 \end{center}
 \end{minipage}
 \caption{3d example: Memory consumption for solving
   one increment problem of the AT-2 functional}
 \label{fig:3d_spectral_memory}
\end{figure}

\bibliographystyle{abbrvnat}
\bibliography{paper-graeser-kienle-sander-tnnmg-for-brittle-fracture}

\begin{thebibliography}{54}
\providecommand{\natexlab}[1]{#1}
\providecommand{\url}[1]{\texttt{#1}}
\expandafter\ifx\csname urlstyle\endcsname\relax
  \providecommand{\doi}[1]{doi: #1}\else
  \providecommand{\doi}{doi: \begingroup \urlstyle{rm}\Url}\fi

\bibitem[Alberti(2000)]{Alberti2000}
G.~Alberti.
\newblock Variational models for phase transitions, an approach via
  $\gamma$-convergence.
\newblock In G.~Buttazzo, A.~Marino, and M.~K.~V. Murthy, editors,
  \emph{Calculus of Variations and Partial Differential Equations: Topics on
  Geometrical Evolution Problems and Degree Theory}, pages 95--114. Springer,
  Berlin, Heidelberg, 2000.
\newblock \doi{10.1007/978-3-642-57186-2_3}.

\bibitem[Ambati et~al.(2015)Ambati, Gerasimov, and
  Lorenzis]{ambati_gerasimov_delorenzis:2015}
M.~Ambati, T.~Gerasimov, and L.~D. Lorenzis.
\newblock A review on phase-field models of brittle fracture and a new fast
  hybrid formulation.
\newblock \emph{Comput. Mech.}, 55:\penalty0 383--405, 2015.
\newblock \doi{10.1007/s00466-014-1109-y}.

\bibitem[Amor et~al.(2009)Amor, Marigo, and Maurini]{amor_marigo_maurini:2009}
H.~Amor, J.-J. Marigo, and C.~Maurini.
\newblock Regularized formulation of the variational brittle fracture with
  unilateral contact: {N}umerical experiments.
\newblock \emph{J. Mech. Phys. Solids}, 57\penalty0 (8):\penalty0 1209--1229,
  2009.
\newblock \doi{10.1016/j.jmps.2009.04.011}.

\bibitem[Baes(2007)]{baes:2007}
M.~Baes.
\newblock Convexity and differentiability properties of spectral functions and
  spectral mappings on {E}uclidean {J}ordan algebras.
\newblock \emph{Linear Algebra and its Applications}, 422\penalty0
  (2):\penalty0 664--700, 2007.
\newblock \doi{10.1016/j.laa.2006.11.025}.

\bibitem[Bastian et~al.(2021)Bastian, Blatt, Dedner, Dreier, Engwer, Fritze,
  Gr\"aser, Gr\"uninger, Kempf, Kl\"ofkorn, Ohlberger, and
  Sander]{bastian_et_al:2021}
P.~Bastian, M.~Blatt, A.~Dedner, N.-A. Dreier, C.~Engwer, R.~Fritze,
  C.~Gr\"aser, C.~Gr\"uninger, D.~Kempf, R.~Kl\"ofkorn, M.~Ohlberger, and
  O.~Sander.
\newblock The {DUNE} framework: Basic concepts and recent developments.
\newblock \emph{Computers \& Mathematics with Applications}, 81:\penalty0
  75--112, 2021.
\newblock \doi{10.1016/j.camwa.2020.06.007}.

\bibitem[Bertsekas(1982)]{bertsekas:1982}
D.~P. Bertsekas.
\newblock Projected {N}ewton methods for optimization problems with simple
  constraints.
\newblock \emph{SIAM J. Control Optim.}, 20\penalty0 (2):\penalty0 221--246,
  1982.
\newblock \doi{10.1137/0320018}.

\bibitem[Bourdin(2007)]{bourdin07}
B.~Bourdin.
\newblock Numerical implementation of the variational formulation for
  quasi-static brittle fracture.
\newblock \emph{Interfaces and free boundaries}, 9\penalty0 (3):\penalty0
  411--430, 2007.

\bibitem[Bourdin et~al.(2000)Bourdin, Francfort, and
  Marigo]{bourdin_francfort_marigo:2000}
B.~Bourdin, G.~Francfort, and J.~J. Marigo.
\newblock Numerical experiments in revisited brittle fracture.
\newblock \emph{J. Mech. Phys. Solids}, 48:\penalty0 797--826, 2000.

\bibitem[Burke et~al.(2010)Burke, Ortner, and S{\"u}li]{burke+ortner+sueli10}
S.~Burke, C.~Ortner, and E.~S{\"u}li.
\newblock An adaptive finite element approximation of a variational model of
  brittle fracture.
\newblock \emph{SIAM Journal on Numerical Analysis}, 48\penalty0 (3):\penalty0
  980--1012, 2010.

\bibitem[Burke et~al.(2013)Burke, Ortner, and
  S{\"u}li]{burke_ortner_sueli:2013}
S.~Burke, C.~Ortner, and E.~S{\"u}li.
\newblock An adaptive finite element approximation of a generalized
  {A}mbrosio--{T}ortorelli functional.
\newblock \emph{Mathematical Models and Methods in Applied Sciences},
  23\penalty0 (9):\penalty0 1663--1697, 2013.

\bibitem[Chen et~al.(2008)Chen, Davis, Hager, and
  Rajamanickam]{chen_davis_hager_rajamanickam:2008}
Y.~Chen, T.~A. Davis, W.~W. Hager, and S.~Rajamanickam.
\newblock Algorithm 887: {CHOLMOD}, supernodal sparse {C}holesky factorization
  and update/downdate.
\newblock \emph{ACM Trans. Math. Softw.}, 35\penalty0 (3):\penalty0
  22:1--22:14, 2008.
\newblock \doi{10.1145/1391989.1391995}.

\bibitem[Dacorogna(1989)]{Dacorogna1989}
B.~Dacorogna.
\newblock \emph{Direct methods in the calculus of variations}.
\newblock Springer, first edition, 1989.

\bibitem[Davis(2004)]{davis:2004}
T.~A. Davis.
\newblock Algorithm 832: {UMFPACK} {V4.3}---an unsymmetric-pattern multifrontal
  method.
\newblock \emph{ACM Trans. Math. Softw.}, 30\penalty0 (2):\penalty0 196--199,
  2004.
\newblock \doi{10.1145/992200.992206}.

\bibitem[Farrell and Maurini(2017)]{farrell+maurini17}
P.~Farrell and C.~Maurini.
\newblock Linear and nonlinear solvers for variational phase-field models of
  brittle fracture.
\newblock \emph{Int. J. Numer. Meth. Engng.}, 109\penalty0 (5):\penalty0
  648--667, 2017.

\bibitem[Gerasimov and {De Lorenzis}(2016)]{gerasimov+delorenzis16}
T.~Gerasimov and L.~{De Lorenzis}.
\newblock A line search assisted monolithic approach for phase-field computing
  of brittle fracture.
\newblock \emph{Comput. Methods Appl. Mech. Engrg.}, 312:\penalty0 276--303,
  2016.

\bibitem[Gerasimov and {De Lorenzis}(2019)]{gerasimov_delorenzis:2019}
T.~Gerasimov and L.~{De Lorenzis}.
\newblock On penalization in variational phase-field models of brittle
  fracture.
\newblock \emph{Comput. Methods Appl. Mech. Engrg.}, 354:\penalty0 990--1026,
  2019.
\newblock \doi{10.1016/j.cma.2019.05.038}.

\bibitem[Glowinski(1984)]{Glowinski1984_NumNonlinVar}
R.~Glowinski.
\newblock \emph{Numerical Methods for Nonlinear Variational Problems}.
\newblock Springer Series in Computational Physics. Springer, third edition,
  1984.

\bibitem[Gr{\"a}ser and Kornhuber(2009)]{GraeserKornhuber2009b}
C.~Gr{\"a}ser and R.~Kornhuber.
\newblock Multigrid methods for obstacle problems.
\newblock \emph{J. Comp. Math.}, 27\penalty0 (1):\penalty0 1--44, 2009.

\bibitem[Gr{\"a}ser and Sander(2014)]{GraeserSander2014_preprint}
C.~Gr{\"a}ser and O.~Sander.
\newblock Truncated nonsmooth {N}ewton multigrid methods for
  simplex-constrained minimization problems.
\newblock Preprint 384, IGPM Aachen, 2014.

\bibitem[Gr{\"a}ser and Sander(2019)]{GraeserSander2019}
C.~Gr{\"a}ser and O.~Sander.
\newblock Truncated nonsmooth {N}ewton multigrid methods for block-separable
  minimization problems.
\newblock \emph{IMA J. Numer. Anal.}, 39:\penalty0 454--481, 2019.
\newblock \doi{10.1093/imanum/dry073}.

\bibitem[Gr\"aser et~al.(2009)Gr\"aser, Sack, and
  Sander]{graeser_sack_sander:2009}
C.~Gr\"aser, U.~Sack, and O.~Sander.
\newblock Truncated nonsmooth {N}ewton multigrid methods for convex
  minimization problems.
\newblock In M.~Bercovier, M.~Gander, R.~Kornhuber, and O.~Widlund, editors,
  \emph{Domain Decomposition Methods in Science and Engineering XVIII},
  volume~70 of \emph{Lecture Notes in Computational Science and Engineering}.
  Springer, 2009.

\bibitem[Heister et~al.(2015)Heister, Wheeler, and
  Wick]{heister_wheeler_wick:2015}
T.~Heister, M.~F. Wheeler, and T.~Wick.
\newblock A primal--dual active set method and predictor--corrector mesh
  adaptivity for computing fracture propagation using a phase-field approach.
\newblock \emph{Comput. Methods Appl. Mech. Engrg.}, 290:\penalty0 466--495,
  2015.
\newblock \doi{10.1016/j.cma.2015.03.009}.

\bibitem[Jodlbauer et~al.(2020)Jodlbauer, Langer, and
  Wick]{jodlbauer_langer_wick:2020}
D.~Jodlbauer, U.~Langer, and T.~Wick.
\newblock Matrix-free multigrid solvers for phase-field fracture problems.
\newblock \emph{Comput. Methods Appl. Mech. Engrg.}, 372\penalty0 (113431),
  2020.
\newblock \doi{10.1016/j.cma.2020.113431}.

\bibitem[Kopani\v{c}\'{a}kov\'{a} and Krause(2020)]{KopanicakovaKrause2020}
A.~Kopani\v{c}\'{a}kov\'{a} and R.~Krause.
\newblock A recursive multilevel trust region method with application to fully
  monolithic phase-field models of brittle fracture.
\newblock \emph{Comput. Methods Appl. Mech. Engrg.}, 360\penalty0 (112720),
  2020.
\newblock \doi{10.1016/j.cma.2019.112720}.

\bibitem[Kopaničáková et~al.(2023)Kopaničáková, Kothari, and
  Krause]{kopanicakova_kothari_krause:2022}
A.~Kopaničáková, H.~Kothari, and R.~Krause.
\newblock Nonlinear field-split preconditioners for solving monolithic
  phase-field models of brittle fracture.
\newblock \emph{Comput. Methods Appl. Mech. Engrg.}, 403\penalty0 (115733),
  2023.
\newblock \doi{10.1016/j.cma.2022.115733}.

\bibitem[Kornhuber(1997)]{kornhuber:1997}
R.~Kornhuber.
\newblock \emph{Adaptive monotone multigrid methods for nonlinear variational
  problems}.
\newblock Vieweg + Teubner Verlag, 1997.
\newblock ISBN 3519027224.

\bibitem[Kristensen and Martínez-Pañeda(2020)]{kristensen_martinez:2020}
P.~K. Kristensen and E.~Martínez-Pañeda.
\newblock Phase field fracture modelling using quasi-{N}ewton methods and a new
  adaptive step scheme.
\newblock \emph{Theoretical and Applied Fracture Mechanics}, 107:\penalty0
  102446, 2020.
\newblock \doi{10.1016/j.tafmec.2019.102446}.

\bibitem[Kuhn et~al.(2015)Kuhn, Schl{\"u}ter, and
  M{\"u}ller]{kuhn+schuelter+mueller15}
C.~Kuhn, A.~Schl{\"u}ter, and R.~M{\"u}ller.
\newblock On degradation functions in phase field fracture models.
\newblock \emph{Computational Materials Science}, 108:\penalty0 374--384, 2015.

\bibitem[Lancioni and Royer-Carfagni(2009)]{lancioni_royer_carfagni:2009}
G.~Lancioni and G.~Royer-Carfagni.
\newblock The variational approach to fracture mechanics. {A} practical
  application to the french {P}anth\'eon in {P}aris.
\newblock \emph{J Elast}, 95:\penalty0 1--30, 2009.
\newblock \doi{10.1007/s10659-009-9189-1}.

\bibitem[Lewis(1996)]{lewis:1996}
A.~S. Lewis.
\newblock Convex analysis on the {H}ermitian matrices.
\newblock \emph{SIAM J. Optimization}, 6\penalty0 (1):\penalty0 164--177, 1996.

\bibitem[Lewis and Sendov(2001)]{lewis_sendov:2001}
A.~S. Lewis and H.~S. Sendov.
\newblock Twice differentiable spectral functions.
\newblock \emph{SIAM J. Matrix Anal. Appl.}, 23\penalty0 (2):\penalty0
  368--386, 2001.
\newblock \doi{10.1137/S089547980036838X}.

\bibitem[Mang et~al.(2020)Mang, Wick, and Wollner]{MangWickWollner2020}
K.~Mang, T.~Wick, and W.~Wollner.
\newblock A phase-field model for fractures in nearly incompressible solids.
\newblock \emph{Comput. Mech.}, 65\penalty0 (1):\penalty0 61--78, 2020.

\bibitem[Marigo et~al.(2016)Marigo, Maurini, and
  Pham]{marigo_maurini_pham:2016}
J.-J. Marigo, C.~Maurini, and K.~Pham.
\newblock An overview of the modelling of fracture by gradient damage models.
\newblock \emph{Meccanica}, 51:\penalty0 3107--3128, 2016.
\newblock \doi{10.1007/s11012-016-0538-4}.

\bibitem[May et~al.(2016)May, Vignollet, and {de
  Borst}]{may_vignollet_deborst:2016}
S.~May, J.~Vignollet, and R.~{de Borst}.
\newblock A new arc-length control method based on the rates of the internal
  and the dissipated energy.
\newblock \emph{Engineering Computations}, 33\penalty0 (1):\penalty0 100--115,
  2016.

\bibitem[Miehe et~al.(2010{\natexlab{a}})Miehe, Hofacker, and
  Welschinger]{miehe+hofacker+welschinger10}
C.~Miehe, M.~Hofacker, and F.~Welschinger.
\newblock A phase field model for rate-independent crack propagation: Robust
  algorithmic implementation based on operator splits.
\newblock \emph{Comput. Methods Appl. Mech. Engrg.}, 199:\penalty0 2765--2778,
  2010{\natexlab{a}}.

\bibitem[Miehe et~al.(2010{\natexlab{b}})Miehe, Welschinger, and
  Hofacker]{miehe+welschinger+hofacker10a}
C.~Miehe, F.~Welschinger, and M.~Hofacker.
\newblock Thermodynamically consistent phase-field models of fracture:
  Variational principles and multi-field {FE} implementations.
\newblock \emph{Int. J. Numer. Meth. Engng.}, 83:\penalty0 1273--1311,
  2010{\natexlab{b}}.

\bibitem[Mielke and {Roub\'\i\v cek}(2015)]{mielke_roubicek:2015}
A.~Mielke and T.~{Roub\'\i\v cek}.
\newblock \emph{Rate-Independent Systems}.
\newblock Springer, 2015.

\bibitem[Modica and Mortola(1977{\natexlab{a}})]{ModicaMortola77a}
L.~Modica and S.~Mortola.
\newblock Un esempio di {$\Gamma$}-convergenza.
\newblock \emph{Boll. Un. Mat. Ital. B}, 14:\penalty0 285--299,
  1977{\natexlab{a}}.

\bibitem[Modica and Mortola(1977{\natexlab{b}})]{ModicaMortola77b}
L.~Modica and S.~Mortola.
\newblock The {$\Gamma$}-convergence of some functionals.
\newblock preprint 77-7, Istituto Matematico `Leonida Tonelli', Università di
  Pisa, 1977{\natexlab{b}}.

\bibitem[Pham et~al.(2011)Pham, Amor, Marigo, and
  Maurini]{pham_amor_marigo_maurini:2011}
K.~Pham, H.~Amor, J.-J. Marigo, and C.~Maurini.
\newblock Gradient damage models and their use to approximate brittle fracture.
\newblock \emph{International Journal of Damage Mechanics}, 20\penalty0
  (4):\penalty0 618--652, 2011.
\newblock \doi{10.1177/1056789510386852}.

\bibitem[Qi and Yang(2003)]{QiYang2003}
H.~Qi and X.~Yang.
\newblock Semismoothness of spectral functions.
\newblock \emph{SIAM J. Matrix Anal. Appl.}, 25\penalty0 (3):\penalty0
  766--783, 2003.

\bibitem[Qi and Sun(1993)]{QiSun:1993}
L.~Qi and J.~Sun.
\newblock A nonsmooth version of {N}ewtons's method.
\newblock \emph{Math. Prog.}, 58:\penalty0 353--367, 1993.

\bibitem[Sander(2020)]{sander:2020}
O.~Sander.
\newblock \emph{DUNE --- The Distributed and Unified Numerics Environment}.
\newblock Springer, 2020.
\newblock \doi{10.1007/978-3-030-59702-3}.

\bibitem[Sander and Jaap(2020)]{sander_jaap:2020}
O.~Sander and P.~Jaap.
\newblock Solving primal plasticity increment problems in the time of a single
  predictor--corrector iteration.
\newblock \emph{Comput. Mech.}, 65:\penalty0 663--685, 2020.
\newblock \doi{10.1007/s00466-019-01788-y}.

\bibitem[Singh et~al.(2016)Singh, Verhoosel, {de Borst}, and {van
  Brummelen}]{singh_verhoosel_deborst_vanbrummelen:2016}
N.~Singh, C.~Verhoosel, R.~{de Borst}, and E.~{van Brummelen}.
\newblock A fracture-controlled path-following technique for phase-field
  modeling of brittle fracture.
\newblock \emph{Finite Elements in Analysis and Design}, 113:\penalty0 14--29,
  2016.
\newblock \doi{10.1016/j.finel.2015.12.005}.

\bibitem[Steinke and Kaliske(2018)]{steinke+kaliske18}
C.~Steinke and M.~Kaliske.
\newblock A phase-field crack model based on directional stress decomposition.
\newblock \emph{Comput. Mech.}, pages 1--28, 2018.

\bibitem[Thomas(2010)]{Thomas_phd:2010}
M.~Thomas.
\newblock \emph{Rate-independent Damage Processes in Nonlinearly Elastic
  Materials}.
\newblock PhD thesis, Humboldt-Universität zu Berlin, 2010.

\bibitem[Ulbrich(2002)]{Ulbrich2002}
M.~Ulbrich.
\newblock \emph{Nonsmooth {N}ewton-like Methods for Variational Inequalities
  and Constrained Optimization Problems in Function Spaces}.
\newblock Habilitationsschrift, Technische Universit{\"a}t M{\"u}nchen, 2002.

\bibitem[Wambacq et~al.(2021)Wambacq, Ulloa, Lombaert, and
  François]{wambacq_ulloa_lombaert_francois:2021}
J.~Wambacq, J.~Ulloa, G.~Lombaert, and S.~François.
\newblock Interior-point methods for the phase-field approach to brittle and
  ductile fracture.
\newblock \emph{Comput. Methods Appl. Mech. Engrg.}, 375\penalty0 (113612),
  2021.
\newblock \doi{10.1016/j.cma.2020.113612}.

\bibitem[Wheeler et~al.(2014)Wheeler, Wick, and
  Wollner]{wheeler+wick+wollner14}
M.~Wheeler, T.~Wick, and W.~Wollner.
\newblock An augmented-{L}agrangian method for the phase-field approach for
  pressurized fractures.
\newblock \emph{Comput. Methods Appl. Mech. Engrg.}, 271:\penalty0 69--85,
  2014.

\bibitem[Wick(2017{\natexlab{a}})]{wick17error}
T.~Wick.
\newblock An error-oriented {N}ewton/inexact augmented {L}agrangian approach
  for fully monolithic phase-field fracture propagation.
\newblock \emph{SIAM Journal on Scientific Computing}, 39\penalty0
  (4):\penalty0 B589--B617, 2017{\natexlab{a}}.

\bibitem[Wick(2017{\natexlab{b}})]{wick17modified}
T.~Wick.
\newblock Modified {N}ewton methods for solving fully monolithic phase-field
  quasi-static brittle fracture propagation.
\newblock \emph{Comput. Methods Appl. Mech. Engrg.}, 325:\penalty0 577--611,
  2017{\natexlab{b}}.
\newblock \doi{10.1016/j.cma.2017.07.026}.

\bibitem[Wu(2018)]{wu:2018}
J.~Wu.
\newblock Numerical implementation of non-standard phase-field damage models.
\newblock \emph{Comput. Methods Appl. Mech. Engrg.}, 340:\penalty0 767--797,
  2018.

\bibitem[Wu et~al.(2020)Wu, Huang, and Nguyen]{wu_huang_nguyen:2020}
J.-Y. Wu, Y.~Huang, and V.~P. Nguyen.
\newblock On the {BFGS} monolithic algorithm for the unified phase field damage
  theory.
\newblock \emph{Comput. Methods Appl. Mech. Engrg.}, 360:\penalty0 112704,
  2020.
\newblock \doi{10.1016/j.cma.2019.112704}.

\end{thebibliography}

\end{document}